\documentclass[sn-mathphys-num]{sn-jnl}

\usepackage{graphicx}%
\usepackage{multirow}%
\usepackage{amsmath,amssymb,amsfonts}%
\usepackage{amsthm}%
\usepackage{mathrsfs}%
\usepackage[title]{appendix}%
\usepackage{xcolor}%
\usepackage{textcomp}%
\usepackage{manyfoot}%
\usepackage{booktabs}%
\usepackage{algorithm}%
\usepackage{algorithmicx}%
\usepackage{algpseudocode}%
\usepackage{listings}%

\theoremstyle{thmstyleone}%
\theoremstyle{thmstyletwo}%

\theoremstyle{thmstylethree}%

\usepackage{amssymb}
\usepackage{latexsym}
\usepackage{lscape}
\usepackage{comment}
\usepackage{amscd}
\usepackage{wasysym}  
\usepackage{tikz}
\usetikzlibrary{matrix,arrows}
\usepackage{tikz-cd}
\usepackage{appendix}
\usepackage[nice]{nicefrac}
\usepackage{dsfont, stmaryrd}
\usepackage{xparse}
\usepackage{appendix}

\usepackage{color}

\usetikzlibrary{decorations.markings}

\makeatletter
\tikzcdset{
  open/.code     = {\tikzcdset{hook, circled};},
  closed/.code   = {\tikzcdset{hook, slashed};},
  open'/.code    = {\tikzcdset{hook', circled};},
  closed'/.code  = {\tikzcdset{hook', slashed};},
  circled/.code  = {\tikzcdset{markwith = {\draw (0,0) circle (.375ex);}};},
  slashed/.code  = {\tikzcdset{markwith = {\draw[-] (-.4ex,-.4ex) -- (.4ex,.4ex);}};},
  markwith/.code ={
    \pgfutil@ifundefined%
    {tikz@library@decorations.markings@loaded}%
    {\pgfutil@packageerror{tikz-cd}{You need to say %
      \string\usetikzlibrary{decorations.markings} to use arrows with markings}{}}{}%
    \pgfkeysalso{/tikz/postaction = {
      /tikz/decorate,
      /tikz/decoration={markings, mark = at position 0.5 with {#1}}}
    }
  },
}
\makeatother

\renewcommand{\geq}{\geqslant}
\renewcommand{\leq}{\leqslant}

\newtheorem{thm}{Theorem}[section]

\newtheorem{propo}[thm]{Proposition}

\newtheorem{propodef}[thm]{Definition and Proposition}
\newtheorem{lem}[thm]{Lemma}
\newtheorem{sublem}[thm]{Sublemma}
\newtheorem{lem-def}[thm]{Lemma-Definition}
\newtheorem{cor}[thm]{Corollary}
\newtheorem{conject}[thm]{Conjecture}
\newtheorem{propert}[thm]{Properties}
\newtheorem{observ}[thm]{Observation}

\newtheorem{fac}[thm]{Fact}

\newtheorem{notat}[thm]{Notation}
\theoremstyle{definition}
\newtheorem*{ack}{Acknowledgement}

\newtheorem{ex}[thm]{Example}

\newtheorem{rmk}[thm]{Remark}
\newtheorem{dfn}[thm]{Definition}
\newtheorem{quest}[thm]{Question}

\newtheorem*{abs}{Abstract}
\numberwithin{equation}{section}

\newcommand{\nc}{\newcommand}

\nc{\abst}{\begin{abs}} \nc{\xabst}{\end{abs}}
\nc{\theo}{\begin{thm}} \nc{\xtheo}{\end{thm}}
\nc{\prop}{\begin{propo}} \nc{\xprop}{\end{propo}}

\nc{\nota}{\begin{notat}} \nc{\xnota}{\end{notat}}
\nc{\depr}{\begin{propodef}} \nc{\xdepr}{\end{propodef}}
\nc{\lemm}{\begin{lem}} \nc{\xlemm}{\end{lem}}
\nc{\sublemm}{\begin{sublem}} \nc{\xsublemm}{\end{sublem}}
\nc{\lemmdefi}{\begin{lem-def}} \nc{\xlemmdefi}{\end{lem-def}}
\nc{\coro}{\begin{cor}} \nc{\xcoro}{\end{cor}}
\nc{\conj}{\begin{conject}} \nc{\xconj}{\end{conject}}
\nc{\proper}{\begin{propert}} \nc{\xproper}{\end{propert}}
\nc{\obse}{\begin{observ}} \nc{\xobse}{\end{observ}}
\nc{\ques}{\begin{quest}} \nc{\xques}{\end{quest}}

\nc{\fact}{\begin{fac}} \nc{\xfact}{\end{fac}}
\nc{\ackn}{\begin{ack}} \nc{\xackn}{\end{ack}}
\nc{\exam}{\begin{ex}} \nc{\xexam}{\end{ex}}
\nc{\rema}{\begin{rmk}} \nc{\xrema}{\end{rmk}}
\nc{\defi}{\begin{dfn}} \nc{\xdefi}{\end{dfn}}

\nc{\pf}{\begin{proof}} \nc{\xpf}{\end{proof}}

\nc{\on}{\operatorname}
\nc{\fraka}{{\mathfrak a}} \nc{\bba}{{\mathbf a}}
\nc{\frakb}{{\mathfrak b}}
\nc{\frakc}{{\mathfrak c}}
\nc{\frakd}{{\mathfrak d}}
\nc{\frake}{{\mathfrak e}}
\nc{\frakf}{{\mathfrak f}}
\nc{\frakg}{{\mathfrak g}}
\nc{\frakh}{{\mathfrak h}}
\nc{\fraki}{{\mathfrak i}}
\nc{\frakj}{{\mathfrak j}}
\nc{\frakk}{{\mathfrak k}}
\nc{\frakl}{{\mathfrak l}}
\nc{\frakm}{{\mathfrak m}}
\nc{\frakn}{{\mathfrak n}}
\nc{\frako}{{\mathfrak o}}
\nc{\frakp}{{\mathfrak p}}
\nc{\frakq}{{\mathfrak q}}
\nc{\frakr}{{\mathfrak r}}
\nc{\fraks}{{\mathfrak s}}
\nc{\frakt}{{\mathfrak t}}
\nc{\fraku}{{\mathfrak u}}
\nc{\frakv}{{\mathfrak v}}
\nc{\frakw}{{\mathfrak w}}
\nc{\frakx}{{\mathfrak x}}
\nc{\fraky}{{\mathfrak y}}
\nc{\frakz}{{\mathfrak z}}
\nc{\frakA}{{\mathfrak A}}
\nc{\frakB}{{\mathfrak B}}
\nc{\frakC}{{\mathfrak C}}
\nc{\frakD}{{\mathfrak D}}
\nc{\frakE}{{\mathfrak E}}
\nc{\frakF}{{\mathfrak F}}
\nc{\frakG}{{\mathfrak G}}
\nc{\frakH}{{\mathfrak H}}
\nc{\frakI}{{\mathfrak I}}
\nc{\frakJ}{{\mathfrak J}}
\nc{\frakK}{{\mathfrak K}}
\nc{\frakL}{{\mathfrak L}}
\nc{\frakM}{{\mathfrak M}}
\nc{\frakN}{{\mathfrak N}}
\nc{\frakO}{{\mathfrak O}}
\nc{\frakP}{{\mathfrak P}}
\nc{\frakQ}{{\mathfrak Q}}
\nc{\frakR}{{\mathfrak R}}
\nc{\frakS}{{\mathfrak S}}
\nc{\frakT}{{\mathfrak T}}
\nc{\frakU}{{\mathfrak U}}
\nc{\frakV}{{\mathfrak V}}
\nc{\frakW}{{\mathfrak W}}
\nc{\frakX}{{\mathfrak X}}
\nc{\frakY}{{\mathfrak Y}}
\nc{\frakZ}{{\mathfrak Z}}
\nc{\bbA}{{\mathbb A}}
\nc{\bbB}{{\mathbb B}}
\nc{\bbC}{{\mathbb C}}
\nc{\bbD}{{\mathbb D}}
\nc{\bbE}{{\mathbb E}}
\nc{\bbF}{{\mathbb F}} \nc{\bbf}{{\mathbf f}}
\nc{\bbG}{{\mathbb G}}
\nc{\bbH}{{\mathbb H}}
\nc{\bbI}{{\mathbb I}}
\nc{\bbJ}{{\mathbb J}}
\nc{\bbK}{{\mathbb K}}
\nc{\bbL}{{\mathbb L}}
\nc{\bbM}{{\mathbb M}}
\nc{\bbN}{{\mathbb N}}
\nc{\bbO}{{\mathbb O}}
\nc{\bbP}{{\mathbb P}}
\nc{\bbQ}{{\mathbb Q}}
\nc{\bbR}{{\mathbb R}}
\nc{\bbS}{{\mathbb S}}
\nc{\bbT}{{\mathbb T}}
\nc{\bbU}{{\mathbb U}}
\nc{\bbV}{{\mathbb V}}
\nc{\bbW}{{\mathbb W}}
\nc{\bbX}{{\mathbb X}}
\nc{\bbY}{{\mathbb Y}}
\nc{\bbZ}{{\mathbb Z}}
\nc{\calA}{{\mathcal A}}
\nc{\calB}{{\mathcal B}}
\nc{\calC}{{\mathcal C}}
\nc{\calD}{{\mathcal D}}
\nc{\calE}{{\mathcal E}}
\nc{\calF}{{\mathcal F}}
\nc{\calG}{{\mathcal G}}
\nc{\calH}{{\mathcal H}}
\nc{\calI}{{\mathcal I}}
\nc{\calJ}{{\mathcal J}}
\nc{\calK}{{\mathcal K}}
\nc{\calL}{{\mathcal L}}
\nc{\calM}{{\mathcal M}}
\nc{\calN}{{\mathcal N}}
\nc{\calO}{{\mathcal O}}
\nc{\calP}{{\mathcal P}}
\nc{\calQ}{{\mathcal Q}}
\nc{\calR}{{\mathcal R}}
\nc{\calS}{{\mathcal S}}
\nc{\calT}{{\mathcal T}}
\nc{\calU}{{\mathcal U}}
\nc{\calV}{{\mathcal V}}
\nc{\calW}{{\mathcal W}}
\nc{\calX}{{\mathcal X}}
\nc{\calY}{{\mathcal Y}}
\nc{\calZ}{{\mathcal Z}}

\nc{\Spac}{{{Spaces}}}
\nc{\scrA}{{\mathscr A}}
\nc{\scrE}{{\mathscr E}}
\nc{\scrR}{{\mathscr R}}

\nc{\Bmu}{\mbox{$\raisebox{-0.59ex}{$l$}\hspace{-0.18em}\mu\hspace{-0.88em}\raisebox{-0.98ex}{\scalebox{2}{$\color{white}.$}}\hspace{-0.416em}\raisebox{+0.88ex}{$\color{white}.$}\hspace{0.46em}$}{}}

\nc{\bnu}{{\bar{ \nu}}}

\nc{\olO}{\bar{\calO}}

\nc{\al}{{\alpha}} 
\nc{\be}{{\beta}}
\nc{\ga}{{\gamma}} \nc{\Ga}{{\Gamma}}
 \nc{\hGa}{\hat{\Gamma}}
\nc{\ve}{{\varepsilon}} 
\nc{\la}{{\lambda}} \nc{\La}{{\Lambda}}
\nc{\om}{\omega} \nc{\Om}{\Omega} 
\nc{\sig}{{\sigma}} \nc{\Sig}{{\Sigma}}

\nc{\tnb}{\psi_{\rm tame}}
\nc{\oM}{\overline{{M}}}
\nc{\op}{{\on{op}}}
\nc{\ad}{{\on{ad}}}
\nc{\alg}{{\on{alg}}}
\nc{\Ad}{{\on{Ad}}}
\nc{\Adm}{{\on{Adm}}} \nc{\aff}{{\on{aff}}}
\nc{\Aut}{{\on{Aut}}}
\nc{\Bun}{{\on{Bun}}}
\nc{\cha}{{\on{char}}}
\nc{\der}{{\on{der}}}
\nc{\Der}{{\on{Der}}}
\nc{\diag}{{\on{diag}}}
\nc{\End}{{\on{End}}}
\nc{\Fl}{{\calF\!\ell}}
\nc{\Tr}{{\on{Transp}}}
\nc{\TR}{{\calT\!\calR}}
\nc{\Gal}{{\on{Gal}}}
\nc{\Gr}{{\on{Gr}}}
\nc{\rH}{{\on{H}}}
\nc{\Hom}{{\on{Hom}}}
\nc{\IC}{{\on{IC}}}
\nc{\id}{{\on{id}}}
\nc{\Id}{{\on{Id}}}
\nc{\ind}{{\on{ind}}}
\nc{\Ind}{{\on{Ind}}}
\nc{\Lie}{{\on{Lie}}}
\nc{\Pic}{{\on{Pic}}}
\nc{\pr}{{\on{pr}}}
\nc{\Res}{{\on{Res}}}
\nc{\res}{{\on{res}}} \nc{\Sat}{{\on{Sat}}}
\nc{\s}{{\on{sc}}}
\nc{\drv}{{\on{der}}}
\nc{\sgn}{{\on{sgn}}}
\nc{\Spec}{{\on{Spec}}}\nc{\Spf}{\on{Spf}} 
\nc{\Sph}{\on{Sph}}
\nc{\St}{{\on{St}}}
\nc{\tr}{{\on{tr}}}
\nc{\Mod}{{\mathrm{-Mod}}}
\nc{\Hilb}{{\on{Hilb}}} 
\nc{\Ext}{{\on{Ext}}} 
\nc{\vs}{{\on{Vec}}}
\nc{\ev}{{\on{ev}}}
\nc{\nO}{{\breve{\calO}}}
\nc{\tS}{{\tilde{S}}}
\nc{\spe}{{\on{sp}}}
\nc{\loc}{{\on{loc}}}
\nc{\Sym}{{\on{Sym}}}
\nc{\Cone}{{\on{C}}}
\nc{\syn}{{\on{syn}}}
\nc{\reg}{{\on{reg}}}
\nc{\colim}{{\on{colim}}}
\nc{\Norm}{{\on{N}}}

\nc{\nscrR}{{\mathscr{R}^{\on{nr}}}}

\nc{\GL}{{\on{GL}}}
\nc{\U}{{\on{U}}}
\nc{\Gl}{\on{Gl}} 
\nc{\GSp}{{\on{GSp}}}
\nc{\gl}{{\frakg\frakl}}
\nc{\SL}{{\on{SL}}} 
\nc{\SU}{{\on{SU}}} 
\nc{\SO}{{\on{SO}}}
\nc{\PGL}{{\on{PGL}}}

\nc{\Conv}{{\on{Conv}}}
\nc{\Rep}{{\on{Rep}}}
\nc{\Dom}{{\on{Dom}}}
\nc{\red}{{\on{red}}}
\nc{\act}{{\on{act}}}
\nc{\nr}{{\on{nr}}}
\nc{\ctf}{{\on{ctf}}}

\nc{\str}{{\on{-}}} 
\nc{\os}{{\bar{s}}}
\nc{\oeta}{{\bar{\eta}}}

\nc{\hookto}{\hookrightarrow}
\nc{\longto}{\longrightarrow}
\nc{\leftto}{\leftarrow}
\nc{\onto}{\twoheadrightarrow}
\nc{\lonto}{\twoheadleftarrow}

\nc{\uG}{{\underline{G}}}
\nc{\uA}{{\underline{A}}}
\nc{\uS}{{\underline{S}}}
\nc{\uT}{{\underline{T}}}
\nc{\uM}{{\underline{M}}}
\nc{\uP}{{\underline{P}}}
\nc{\uB}{{\underline{B}}}
\nc{\uN}{{\underline{N}}}

\nc{\ucG}{{\underline{\calG}}}
\nc{\ucA}{{\underline{\calA}}}
\nc{\ucS}{{\underline{\calS}}}
\nc{\ucT}{{\underline{\calT}}}
\nc{\ucalM}{{\underline{\calM}}}
\nc{\ucP}{{\underline{\calP}}}
\nc{\ucalN}{{\underline{\calN}}}

\nc{\bF}{{\breve{F}}}

\nc{\oFl}{{\overline{\Fl}}} 
\nc{\bU}{{\overline{U}}}
\nc{\tGr}{{\tilde{\Gr}}}
\nc{\cGr}{\calG\! r}
\nc{\oGr}{\overline{\on{Gr}}} 
\nc{\ocGr}{\overline{\calG\! r}}
\nc{\co}{{\colon}}
\nc{\sch}[1]{(Sch/{#1})}
\nc{\HypLoc}[1]{HypLoc({#1})}

\nc{\ohtimes}{\stackrel{!}{\otimes}}
\nc{\boxtilde}{\widetilde{\boxtimes}}
\nc{\vstar}{{\varhexstar}}

\nc{\Div}{\on{Div}}
\nc{\Sht}{\on{Sht}}
\nc{\Frob}{\on{Frob}}

\nc{\x}{\times}
\nc{\bsl}{\backslash}
\nc{\algQl}{{\bar{\bbQ}_\ell}}
\nc{\sF}{{\bar{F}}}
\nc{\nF}{{\breve{F}}}
\nc{\nW}{{W^{\on{nr}}}}
\nc{\sk}{{\bar{k}}}
\nc{\cont}{\on{c}}
\nc{\Supp}{\on{Supp}}
\nc{\blt}{\bullet}  
\nc{\dom}{\on{dom}}
\nc{\scon}{{\on{sc}}} 
\nc{\Affine}{\on{Aff}} 
\nc{\nscrA}{\mathscr{A}^{\on{nr}}} 
\nc{\nfraka}{{\bbf^{\on{nr}}}}
\nc{\ran}{{\rangle}}
\nc{\lan}{{\langle}}
\nc{\bk}{{\bar{k}}}
\nc{\tF}{{\tilde{F}}}
\nc{\sS}{{\bar{S}}}
\nc{\LG}{{^\text{L}\hspace{-0.04cm}G}}
\nc{\LL}{{^\text{L}\hspace{-0.07cm}L}}
\nc{\et}{{\text{\rm \'et}}}
\nc{\inv}{{\on{inv}}}
\nc{\Hecke}{{\on{Hecke}}}
\nc{\Isom}{{\on{Isom}}}
\nc{\oSht}{{\overline{\on{Sht}}}}
\nc{\umu}{{\underline \mu}}
\nc{\AIJ}{{\calO_X[{\scriptstyle{\calI\over \calJ}}]}}
\nc{\Proj}{{\on{Proj}}}
\nc{\Bl}{{\on{Bl}}}

\nc{\Pos}{{\on{Pos}}}
\nc{\Sets}{{\on{Sets}}}
\nc{\AffSch}{{\on{AffSch}}}
\nc{\Groups}{{\on{Groups}}}
\nc{\Gpds}{{\on{Groupoids}}}
\nc{\Sch}{{\on{Sch}}}
\nc{\fl}{{\on{flat}}}

\nc{\pot}[1]{ [\hspace{-0,5mm}[ {#1} ]\hspace{-0,5mm}] }
\nc{\rpot}[1]{ (\hspace{-0,7mm}( {#1} )\hspace{-0,7mm}) }

\nc{\defined}{\hspace{0.1cm}\stackrel{\text{\tiny \rm def}}{=}\hspace{0.1cm}}

\nc{\nn}{\mathrm{n}}
\begin{document}

\title[\textbf{Multi-centered dilatations, congruent isomorphisms and Rost double deformation space}]{\textbf{Multi-centered dilatations, congruent isomorphisms and Rost double deformation space}}


\author[]{\fnm{Arnaud} \sur{Mayeux}  }\email{arnaud.mayeux@mail.huji.ac.il  
 }

\affil[]{\orgdiv{Einstein Institute of Mathematics}, \orgname{The Hebrew University of Jerusalem}, \orgaddress{\street{Givat Ram}, \city{Jerusalem}, \postcode{9190401},  \country{Israel}}}


\abstract{We introduce multi-centered dilatations of rings, schemes and algebraic spaces, a basic algebraic concept. Dilatations of schemes endowed with a structure (e.g. monoid, group or Lie algebra) are in favorable cases  schemes endowed with the same structure.  As applications, we use our new formalism to contribute to the understanding of mono-centered dilatations, to formulate and deduce some multi-centered congruent isomorphisms and to interpret Rost double deformation space as a double-centered dilatation.}




\maketitle

\begin{center}

$~~$

\end{center}

\ackn
 The idea behind this project appeared while I was preparing a talk at the workshop \textit{Arthur packets and related problems, University of Zagreb, June 2022, Croatia}, I thank Marcela Hanzer for the invitation. I am grateful to Anne-Marie Aubert, Jeffrey Adler, Adrien Dubouloz, Shai Evra, Lucy Moser-Jauslin, Timo Richarz, Matthieu Romagny, João Pedro dos Santos and Simon Riche for useful discussions, support, interest and suggestions. I thank the referees for useful comments and suggestions. Section \ref{Motiv} of the introduction was written in parallel with \cite{DMdS23} and is partly identical with some paragraphs there. \xackn

\tableofcontents

\section{Introduction}

We start with the conventions that we will use, as it will be important to appreciate fully the rest of the introduction. 
\subsection{Notation and conventions}

\begin{enumerate}
\item We use the terminology of \cite{stacks-project} for basic definitions such as ring, scheme, algebraic space, localization, etc; except that the symbol $\bbN $ denotes $ \{0, 1 , 2 ,3 \ldots \}$ in this paper (compare with \cite[\href{https://stacks.math.columbia.edu/tag/055X}{Tag 055X}]{stacks-project}). In particular, a multiplicative subset of a ring is defined as a submonoid of $(A, \times)$ (it always contains $1$) \cite[\href{https://stacks.math.columbia.edu/tag/00CN}{Tag 00CN}]{stacks-project}. For any element $a$ in a ring $A$, $a^0$ equals $1$ and the set $a^\bbN:= \{ a^k | k \in \bbN \}$ is a multiplicative subset. For example $0^\bbN$ is a multiplicative subset. The zero ring is a ring. The spectrum of a ring is empty if and only if it is the zero ring \cite[\href{https://stacks.math.columbia.edu/tag/00E0}{Tag 00E0}]{stacks-project}.
\item A locally principal closed subscheme of a scheme is a closed subscheme whose sheaf of ideals is locally generated by a single element \cite[\href{https://stacks.math.columbia.edu/tag/01WR}{Tag 01WR}]{stacks-project}. Note that this is a very flexible notion: it is stable by pullbacks under our conventions \cite[\href{https://stacks.math.columbia.edu/tag/053P}{Tag 053P}]{stacks-project}.
\item An effective Cartier divisor is a locally principal closed subscheme of a scheme which is locally defined by a non-zero-divisor, cf. \cite[\href{https://stacks.math.columbia.edu/tag/01WQ}{Tag 01WQ}]{stacks-project} and \cite[\href{https://stacks.math.columbia.edu/tag/01WS}{Tag 01WS}]{stacks-project} (historically, cf. \cite{Ca58}). We sometimes omit "effective" in this text as all Cartier divisors that we use are effective. This is a very useful notion, but it is not stable by pullbacks in general \cite[\href{https://stacks.math.columbia.edu/tag/01WV}{Tag 01WV}]{stacks-project}, however the pullback of an effective Cartier divisor is an effective Cartier divisor in several cases, cf. e.g. \cite[\href{https://stacks.math.columbia.edu/tag/02OO}{Tag 02OO}]{stacks-project}.
\item The concepts of locally principal closed subschemes and Cartier divisors also make sense for algebraic spaces, cf. \cite[\href{https://stacks.math.columbia.edu/tag/083B}{Tag 083B}]{stacks-project} and \cite[\href{https://stacks.math.columbia.edu/tag/083C}{Tag 083C}]{stacks-project}.
\item Let $X$ be a scheme (or an algebraic space) over a fixed base scheme $S$. If $f: Y \to X $ is an $S$-morphism and $C \to S$ is an immersion, we denote by $f|_{C} : Y \times _S C \to X \times _S C$ the canonical morphism obtained by base change, and call it the restriction to $C$. 
\item  Let $(M, +)$ be a monoid. A submonoid $F$ is a face of $M$ if whenever $x + y \in  F$, then both $x$
and $y$ belong to $F.$
\end{enumerate}

\subsection{Motivation and goals}\label{Motiv}
 Dilatation of rings is a basic construction of commutative algebra, like localization or tensor product. It can be globalized so that it also makes sense on schemes or algebraic spaces. In fact dilatations generalize localizations.
 
Let $A$ be a ring and let $S$ be a multiplicative subset of $A$. Recall that the localization $S^{-1}A$ is an $A$-algebra such that for any $A$-algebra $ A \to B $ such that the image of $s$ is invertible for any $s \in S$, then $A \to B$ factors through $A \to S^{-1}A$. Intuitively, $S^{-1}A$ is the $A$-algebra obtained from $A$ adding all fractions $\frac{a}{s}$ with $a \in A $ and $s \in S$.
Formally, $S^{-1} A$ is made of classes of fractions $\frac{a}{s}$ where $a \in A$ and $s \in S $ (two representatives $\frac{a}{s}$ and $\frac{b}{t}$ are identified if $atr=bsr$ for some $r \in S$), addition and multiplication are given by usual formulas. Now let us give for any element $s \in S$ an ideal $M_s$ of $A$ containing $s$. The dilatation of $A$ relative to the data $S, \{M_s \}_{s \in S}$ (introduced in this paper) is an $A$-algebra $A'$ obtained intuitively by adding to $A$ only the fractions $\frac{m}{s}$ with $s \in S $ and $m \in M_s$. The dilatation $A'$ satisfies that for any $s \in S$, we have $s A'= M_s A' $ (intuitively any $m \in M_s$ belongs to $s A'$, i.e. becomes a multiple of $s$, so that we have an element $\frac{m}{s} $ such that $m = s  \frac{m}{s} $). As a consequence of the construction, the elements $s \in S$ become a non-zero-divisor in $A'$ so that $\frac{m}{s}$ is well-defined (i.e. unique). It turns out that it is convenient, with dilatations of schemes in mind, to make a bit more flexible the above framework, namely to remove the conditions that $S$ is multiplicative and that $s \in M_s$, so we use the following definition. 

\begin{flushleft}
\textsc{Definition.} Let $A$ be a ring. Let $I$ be an index set. A multi-center in $A$ indexed by $I$ is a set of pairs $\{[M_i , a_i ]\}_{i \in I}$ where for each $i$, $M_i$ is an ideal of $A$ and $a_i$ is an element of $A$.
\end{flushleft}

For each multi-center $\{[M_i , a_i ]\}_{i \in I}$, one has the dilatation $A [\{\frac{M_i}{a_i}\}_{i \in I}]$, it is an $A$-algebra.  We will define and study in details dilatations of rings in Section \ref{sectionalgebr}, in particular we will state formally the universal property they enjoy. We will see that $A [\{\frac{M_i}{a_i}\}_{i \in I}]$ is generated, as $A$-algebra, by $\{\frac{M_i}{a_i}\}_{i \in I}$. We will also see that if $M_i = A$ for all $i$, then $A [\{\frac{M_i}{a_i}\}_{i \in I}]= S^{-1} A$ where $S$ is the multiplicative subset generated by $\{a _i \} _{i \in I }$. Reciprocally, we will see that any sub-$A$-algebra of a localization $S^{-1} A$ for a certain $S$ is isomorphic to a dilatation of $A$. 

Dilatations of schemes and algebraic spaces are obtained from dilatations of rings via glueing. We introduce the following definition (we restrict to schemes in this introduction).

\begin{flushleft}

\textsc{Definition.} Let $X$ be a scheme. Let $I$ be an index set. A multi-center in $X$ indexed by $I$ is a set of pairs $\{[Y_i, D_i]\}_{i \in I } $ such that $Y_i$ and $D_i $ are closed subschemes for each $i$ and such that locally, all $D_i$ are principal for $i \in I$.
\end{flushleft}

Associated to each multi-center, one has the dilatation $\Bl \big\{{}^{D_i}_{Y_i} \big\}_{i \in I } X $, it is a scheme endowed with a canonical affine morphism $f: \Bl \big\{{}^{D_i}_{Y_i} \big\}_{i \in I } X \to X$. It satisfies, in a universal way, that $f^{-1} (D_i)$ is a Cartier divisor (i.e. is locally defined by a non-zero-divisor) and that $f^{-1} (D_i ) \subset  f^{-1} (Y_i)$ for all $i \in I$. If $\#I=1$, we use the terminology mono-centered dilatation.
 We will introduce formally and study several facets of this construction and show that it enjoys many wonderful properties. 

As we explained, dilatations are a basic construction which can be easily encountered in specific situations. As a consequence, the theory of dilatations has deep and distinguished roots. Right from the start, we warn the reader that we could not present a comprehensive historical account.
 As soon as Cremona and Bertini started using quadratic transformations (or blowups) in the framework of algebraic geometry over fields, ``substitutions'' of the form $x'=x$ and $y'=y/x$ started being made by algebraic geometers, see for example equation (8) in \cite[Section 11]{No1884} and Noether's acknowledgement, at the start of \cite[Section 12]{No1884}, that these  manipulations come from Cremona's point of view. Examples of dilatations appear frequently in some works of Zariski and Abhyankar, cf. \cite[Definition, p. 86]{abhyankar-zariski55} and \cite[proof of Th.4, case (b)]{Za43}. Other forerunner examples of dilatations play a central role in several independent and unrelated works later,  cf. \cite{Da67},  \cite[Section 25]{Ner64} and \cite[Section 4]{artin69}. 
As far as we know, the terminology dilatations emerged in \cite[§3.2]{BLR90}, where a section is devoted to study dilatations of schemes over discrete valuation rings systematically. In the context of schemes over a discrete valuation ring, we draw the reader's attention to \cite{Ana73},  \cite{WW80} and \cite{PY06}.
The paper \cite{KZ99} studies dilatations (under the name affine modifications) systematically in the framework of algebraic geometry over fields.
Over two-dimensional base schemes dilatations also appear in \cite[p. 175]{PZ13}.
In recent times, the authors of \cite{Du05} and \cite{MRR20} have set out to accommodate all these mono-centered constructions in a larger and unified frame, namely for arbitrary schemes. The paper \cite{MRR20} introduces dilatations of arbitrary schemes in the mono-centered case and provides a systematic treatment of mono-centered dilatations of general schemes. An equivalent definition of mono-centered dilatations of general schemes, under the name affine modifications, was introduced earlier in \cite[Définition 2.9]{Du05} under few assumptions.
Set aside localizations, mono-centered dilatations have been the main focus of mathematicians in the past. However, in the context of group schemes over discrete valuation rings, examples of multi-centered dilatations of rings and schemes that are not localizations or mono-centered dilatations appeared and were used in \cite[Exp. VIB Ex. 13.3]{SGA3}, \cite{PY06} and \cite{DHdS18}.  The present paper unifies all these constructions, it introduces dilatations of arbitrary rings, schemes and algebraic spaces for arbitrary multi-centers. Allowing multi-centers also leads naturally to the formulation of combinatorial isomorphisms on dilatations and gives birth to refined universal properties.

Beyond rings and algebraic spaces, the concept of dilatations makes sense for other structures and geometric settings. Let us indicate some constructions already available.
Some dilatation constructions in the framework of complex analytic spaces were introduced in \cite{Ka94}. For many other structures than rings, dilatations also make sense (e.g. categories, non-commutative rings, semirings), this is the topic of \cite{Ma23c}.

Recall that dilatations have distinguished  roots, as a consequence, several other terminologies are used to call certain dilatations in literature.  For examples the constructions named \textit{affine blowups}, \textit{affine modifications}, \textit{automatic blowups}, \textit{formal blowups}, \textit{Kaliman-Zaidenberg modifications}, \textit{localizations} and \textit{Néron blowups} are examples of (possibly multi-centered) dilatations.

\subsection{A view on mono-centered dilatations motivating multi-centered dilatations}

Recall that in this paper, we will introduce and study multi-centered dilatations. Here "multi" is related to the cardinality of the set $I$ appearing in Section \ref{Motiv}. For example if $\#I =2$, we talk about double-centered dilatations. If $\#I$ is finite, we talk about finite-centered dilatations. If $\#I =1$, we talk about mono-centered dilatations. As we will see in this paper, there are many combinatorial properties and isomorphisms relating several multi-centered dilatations (Section §\ref{sectionmonopoly}). In particular, we can relate finite-centered dilatations to mono-centered dilatations. Mono-centered dilatations form a special and fundamental case in the theory of algebraic dilatations. As we explained in the previous section, mono-centered dilatations of schemes were introduced earlier by many authors. The most systematic treatment of mono-centered dilatations was carried out by the authors of the Stacks Project in \cite[\href{https://stacks.math.columbia.edu/tag/052P}{Tag 052P}]{stacks-project} for rings and by Mayeux-Richarz-Romagny in \cite{MRR20} for Grothendieck schemes. In this section, we recall the ideas of mono-centered dilatations as a motivation and preliminary to our more general treatment. As our general treatment will cover all cases and is self-contained, we do not recall all the mono-centered definitions formally in the present section. Instead we provide a certain view on \cite[\href{https://stacks.math.columbia.edu/tag/052P}{Tag 052P}]{stacks-project} and \cite{MRR20} sufficient to understand easily and quickly most parts of the present paper.

Let $A$ be a ring (resp. $X$ be a scheme). Let $L$ be an ideal of $A$ (resp. $Z \subset X$ be a closed subscheme). Note that $L$ is denoted by the letter $I$ in \cite{stacks-project}. Let $a \in L$ be an element (resp. $D \subset X$ a locally principal closed subscheme of $X$, in the sense of \cite[\href{https://stacks.math.columbia.edu/tag/01WR}{Tag 01WR}]{stacks-project}, containing $Z$). In particular we get a mono-center in $A$ (resp. $X$) in the sense of Section \ref{Motiv}.

Stacks Project introduces an $A$-algebra $A[\frac{L}{a}]$, called the affine blowup algebra  \cite[\href{https://stacks.math.columbia.edu/tag/052Q}{Tag 052Q}]{stacks-project}. The definition of $A[\frac{L}{a}]$ in Stacks Project is $(\Bl_L(A))_{(a_{(1)})}$, the degree zero part of the localization of the Rees algebra where $a$ is placed in degree $1$. Then \cite[\href{https://stacks.math.columbia.edu/tag/07Z3}{Tag 07Z3}]{stacks-project} says that (the image of) $a$ is a non-zero-divisor in $A[\frac{L}{a}]$ and that  $L A[\frac{L}{a}] = a A[\frac{L}{a}] $: this precisely says that one has a "fraction"  $l/a \in A[\frac{L}{a}] $ for each $l \in L$. In the multi-centered case, we will mainly use the point of view of fractions instead of Rees algebras. Also we will use the terminology dilatations of rings instead of affine blowup algebras (which comes from Rees algebras and related projective blowups).

The paper \cite{MRR20} introduces an $X$-affine scheme $\Bl _Z^D X$, called the dilatation of $X$ with center $Z,D$ (in \cite{MRR20}, the terminology "affine blowup" is also often used). The definition of $\Bl _Z^D X$ involves a quasi-coherent Rees algebra and is a global counterpart of the definition in \cite[\href{https://stacks.math.columbia.edu/tag/052Q}{Tag 052Q}]{stacks-project}. Then \cite[Lemma 2.4]{MRR20} says that $\Bl_Z^D X \times _X D$ is a Cartier divisor in $\Bl _Z ^D X$ and that $\Bl_Z^D X \times _X D = \Bl_Z^D X \times _X Z$; these are global counterparts of the assertions of \cite[\href{https://stacks.math.columbia.edu/tag/07Z3}{Tag 07Z3}]{stacks-project}.

So we have attached to our center a certain object, the dilatation, satisfying two remarkable properties \cite[\href{https://stacks.math.columbia.edu/tag/07Z3}{Tag 07Z3}]{stacks-project} (resp. \cite[Lemma 2.4]{MRR20}). The next step in the theory of \cite{MRR20} is that in fact the dilatation is universal, in a precise sense,  among all objects satisfying two similar properties. This is formalized by saying that the dilatation represents a certain functor (\cite[Proposition 2.6]{MRR20}). To state this precisely, we need the category $Sch_X^{D\text{-reg}}$ defined in \cite{MRR20} as follows. Let $f:T \to X$ be an $X$-scheme, by  \cite[\href{https://stacks.math.columbia.edu/tag/053P}{Tag 053P}]{stacks-project} we always have that $f^{-1}(D)$ (i.e. $T \times _X D $) is locally principal in $T$. 
The category $Sch_X^{D\text{-reg}}$ is defined as the full subcategory of schemes over $X$ whose objects are morphisms $f:T \to X $ such that the pullback $f^{-1}(D)$ is a Cartier divisor in $T$. We now state the universal property of mono-centered dilatations. 

\prop (\cite[Proposition 2.6]{MRR20}) 
The dilatation $\Bl_Z^DX\to X$ represents the contravariant functor $Sch_X^{D\text{-}\reg}\to Set$ given by
\[
(f\co T\to X) \;\longmapsto\; \begin{cases}\{*\}, \; \text{if $f|_{T\x_XD}$ factors through $Z\subset X$;}\\ \varnothing,\;\text{else.}\end{cases}
\]
\xprop 
In particular, given $f:T \to X $ in $Sch_X^{D\text{-reg}}$, we have that $\# \Hom _{X} ( T , \Bl _Z^D X ) \leq 1$, and we know when this equals $0$ or $1$. 

Then \cite{MRR20} studies several facets of mono-centered dilatations (functoriality, base change, exceptional divisor, preservation of group schemes, group structure on the exceptional divisor, iterated dilatations, flatness and smoothness, etc) many of them will be generalized to multi-centered dilatations in the present paper. However, some results of \cite{MRR20} are not generalized to multi-centered dilatations and are not stated in the present paper (cf. e.g. \cite[Proposition 2.9]{MRR20}, \cite[Theorem 3.5]{MRR20} and \cite[Lemma 3.8]{MRR20}).
In the present paper, we prove some results on multi-centered dilatations which do not make sense in the framework of mono-centered dilatations (cf. e.g. most of Section \ref{sectionalgebr}, Section \ref{sectionsinglediv} and Section \ref{sectionmonopoly}). We also prove some results on multi-centered dilatations which make sense also for mono-centered dilatations but do not appear (in the mono-centered case) in  \cite{MRR20} (cf. e.g. Proposition \ref{defined}, Proposition \ref{LieLie} and Proposition \ref{normalizes}).
For a survey covering, among other things, most of \cite{MRR20} and the present paper, we refer to \cite{DMdS23}.

Recall that, before \cite{MRR20}, mono-centered dilatations of schemes were introduced and studied in specific contexts in e.g. \cite{BLR90}, \cite{KZ99} and \cite{Du05}, but the restricted conventions and frameworks made these pioneering treatments less flexible or natural. Note that \cite{MRR20} cites \cite{BLR90} but cites neither \cite{KZ99} nor \cite{Du05}. Note also that \cite{Du05} cites \cite{KZ99} but does not cite \cite{BLR90}. Finally, note that \cite{KZ99} does not cite \cite{BLR90}. Recall that the terminology "dilatation" appeared in \cite{BLR90} and is now standard. In \cite{KZ99} and \cite{Du05}, the word "dilatation" does not appear. The authors of \cite{KZ99} and \cite{Du05} used the compound word "affine modification". 

To conclude this section, we saw that mono-centered dilatations of rings can be viewed as certains operations adding or imposing some fractions and that dilatations of schemes form a global counterpart of this notion. Localization is also an operation adding or imposing some fractions to rings. The present paper provides answers to the following questions. Is it possible to unify localizations of rings and mono-centered dilatations of rings \cite[\href{https://stacks.math.columbia.edu/tag/052Q}{Tag 052Q}]{stacks-project} in a natural theory? If yes, is it possible to globalize it on schemes following the paradigm of \cite{MRR20}? Is it possible to globalize it on algebraic spaces? 
The present text answers these questions positively. 
On the road, it deals with more general centers than \cite{stacks-project} and \cite{MRR20}, even in the mono-centered case, leading to even more flexibility (i.e we allow centers $[Y,D]$ where $Y$ is not necessarily assumed to be included in $D$).

\subsection{Results}
Let $A$ be a ring and let $\{[M_i , a_i ]\}_{i \in I}$ be a multi-center. Put $L_i = M_i +(a_i)$ for $i \in I $. Let $Ring_{A}^{a\text{-reg}}$ be the full subcategory of $A$-algebras $f:A \to B $ such that $f(a_i)$ is a non-zero-divisor for all $i \in I$. Dilatations of rings enjoy several properties, the main ones are summarized in the following statement. Many other properties are proved in Section \ref{sectionalgebr}.

\prop The following assertions hold. \begin{enumerate} \item The covariant functor from $Ring_{A}^{a\text{-reg}} $ to $ Set $  given by  \[
(f\co A \to B ) \;\longmapsto\; \begin{cases}\{*\}, \; \text{if } f(M_i) B \subset f(a_i) B  \text{ for $i \in I$;}\\ \varnothing,\;\text{else}~~~~\end{cases} 
\] is representable by an $A$-algebra $A[\frac{M}{a}]$, also denoted $A[\big\{ \frac{M_i}{a_i}\big\}_{i \in I}], $ called the dilatation of $A$ with center $\{[M_i , a_i ]\}_{i \in I}$. $~~~~(\ref{univpropdilarings2}, \ref{remaequiano})$
 \item The image of $a_i$ in $A[\frac{M}{a}]$ is a non-zero-divisor for any $i \in I$, moreover \[ a_iA[\frac{M}{a}] = L_i A[\frac{M}{a}]. ~~~~(\ref{faitstakalg1}, \ref{faitstakalg2}) \]

 \item If $A$ is a domain and $a_i \ne 0$ for all $i$, then $A[\frac{M}{a}]$ is a domain. ~~~~(\ref{domain})
 \item If $A$ is reduced, then $A[\frac{M}{a}]$ is reduced. ~~~~(\ref{reduced})
 
 \item Assume that $I = \{1 , \ldots , k \}$ is finite. Then we have a canonical identification of $A$-algebras
 \begin{center}$A[\big\{ \frac{M_i}{a_i}\big\}_{i \in I}]= A[\frac{\sum _{i \in I} ( M_i \cdot \prod _{j \in I \setminus {\{i\}}} a_j )}{a_1 \cdots a_k}] .$~~~~(\ref{multifinimonoring})\end{center}
 
\item If $M_i= A$ for all $i$, then $A[\big\{ \frac{M_i}{a_i}\big\}_{i \in I}]$ identifies with the localization $S^{-1}A$ where $S$ is the multiplicative subset of $A$ generated by $\{ a_i \} _{i \in I}$. ~~~~(\ref{localis})

\item The $A$-algebra  $A[\big\{ \frac{M_i}{a_i}\big\}_{i \in I}]$ identifies with the sub-$A$-algebra of $S^{-1}A$ generated by $\big\{ \frac{M_i}{a_i}\big\}_{i \in I}$, where $S$ is the multiplicative subset of $A$ generated by $\{ a_i \} _{i \in I}$. ~~~~(\ref{dilaloc})

\item Any sub-$A$-algebra of a localization $S^{-1} A$ can be obtained as a multi-centered dilatation.  ~~~~(\ref{factdiladansloc})

\item Let $K \subset I$ be a subset. Assume that for all $i \in I \setminus K$, there exists $k(i) $ in $K$ such that $a_{k(i)} \in L_i $ and $L_i \subset L_{k(i)}$.
Then we have a canonical identification\begin{center}$ A[\big\{ \frac{M_i}{a_i}\big\}_{i \in I}] = A[\big\{ \frac{M_k}{a_k}\big\}_{k \in K}]_{\{ \frac{a_i}{a_{k(i)}}\}_{i \in I \setminus K}}$. ~~~~ (\ref{LiLiprime}) \end{center}

 \item In regular cases (in the sense of \cite[\href{https://stacks.math.columbia.edu/tag/07CV}{Tag 07CV}]{stacks-project}), one has an explicit description of $A[\frac{M}{a}]$ as a quotient of a polynomial algebra. ~~~~(cf. \ref{regupoly} for precise statement)
\end{enumerate}
\xprop 
 
 Now let $X$ be a scheme and let $\{[Y_i, D_i]\}_{i \in I } $ be a multi-center. Put $Z_i = Y_i \cap D_i $ for all $i \in I$. Let $Sch_{X}^{D\text{-reg}}$ be the full subcategory of $X$-schemes $f:T \to X  $ such that $f^{-1}(D_i)$ is a Cartier divisor for all $i \in I$. Our main results on dilatations of schemes are summarized in the following theorem. Note that we work with general algebraic spaces later in the paper. 

\theo The following assertions hold.  \label{theoschintro}
 \begin{enumerate} \item The contravariant functor from $Sch ^{D\text{-reg}}_X$ to $Set$ given by  \[
(f\co T\to X) \;\longmapsto\; \begin{cases}\{*\}, \; \text{if $f|_{T\x_XD_i ^{}}$ factors through $Y_i \subset X$ for $i \in I$;}\\ \varnothing,\;\text{else}\end{cases}
\] is representable by an $X$-affine scheme $\Bl_Y^DX$, also denoted $\Bl \big\{{}^{D_i}_{Y_i} \big\}_{i \in I } X $ and $\Bl _{\{Y_i\}_{i \in I}}^{\{D_i\}_{i \in I} }X$, called the dilatation of $X$ with multi-center $\{[Y_i, D_i]\}_{i \in I } $.~~~~(\ref{blow.up.rep.prop}) 

\item As closed subschemes of $\Bl_Y^DX$, one has, for all $i \in I$,
\[
\Bl_Y^DX\x_X  Z_i =\Bl_Y^DX\x_X D_i,
\]
which is an effective Cartier divisor on $\Bl_Y^DX$.~~~~(\ref{blow.up.Cartier.lemm})

\item Let $ K ,J$ be subsets of $I$ such that $I = K \sqcup J$. Then
\[ \Bl \big\{_{Y_i} ^{D_i}\big\}_{i \in I} X = \Bl \Big\{_{Y_k \times _X  \Bl \big\{_{Y_i}^{D_i}\big\}_{i \in J} X } ^{
D_k \times _X  \Bl \big\{_{Y_i}^{D_i}\big\}_{i \in J} X  }\Big\} _{k \in K } \Bl \big\{_{Y_i}^{D_i}\big\}_{i \in J} X  .\] In particular, there is a unique $X$-morphism \[\Bl \big\{_{Y_i}^{D_i}\big\}_{i \in I} X  \to \Bl \big\{_{Y_i}^{D_i}\big\}_{i \in J} X . ~~~~(\ref{uniquemorphism},\ref{multietape})\]

 \item Let $K \subset I$ be such that \begin{enumerate} \item  $I \setminus K $ is finite, \item for all $i \in I \setminus K $, there exists $k(i) \in K $ such that $Z_{k(i)} \subset Z_i$ and $Z_i \subset D_{k(i)}$. \end{enumerate} Then the canonical $X$-morphism
 \[ \Bl \big\{_{Y_i}^{D_i}\big\}_{i \in I} X  \to \Bl \big\{_{Y_i}^{D_i}\big\}_{i \in K} X . \] is an open immersion. $~~~(\ref{ZZprimeopen})$

\item 
Assume that $\#I=k$ is finite. We fix an arbitrary bijection $I = \{ 1 ,\ldots , k\}$.  We have a canonical isomorphism of $X$-schemes
\[ \Bl _{\{Y_i\}_{i \in I}}^{\{D_i\}_{i \in I} }X  \cong \Bl ^{(\Bl\cdots ) \times _X D_k}_{( \Bl \cdots ) \times _X Y_k }\Biggl(\cdots \Bl_{(\Bl \cdots ) \times _X Y_3}^{(\Bl\cdots ) \times _X D_3}\biggl( \Bl _{(\Bl_{Y_1}^{D_1}X) \times _{X} Y_2} ^{(\Bl_{Y_1}^{D_1}X) \times _{X} D_2} \bigl(\Bl _{Y_1}^{D_1} X \bigl) \biggl) \Biggl). ~~~~(\ref{multisingle}) \] 

\item (Monopoly isomorphism) Assume that $\#I=k$ is finite. We fix an arbitrary bijection $I = \{ 1 ,\ldots , k\}$. We have a unique isomorphism of $X$-schemes
 \[\Bl _{\{Y_i\}_{i \in I}}^{\{D_i\}_{i \in I} }X  \cong \Bl _{\bigcap _{i \in I }( Y_i + D_1+ \ldots + D_{i-1} + D_{i+1} + \ldots + D_k)} ^{D_1+ \ldots + D_k} X.~~~~(\ref{formulamultimono})\]

\item In many cases, multi-centered dilatations can be iterated in a compatible way with the addition law on closed subschemes.~~~~(\ref{monoidclosed}, \ref{nutheta})

\item Multi-centered dilatations whose multi-center $\{[Y_i, D_i]\}_{i \in I } $ satisfies that $\{ D_i\} _{i \in I }$ are multiples of a single locally principal closed subscheme satisfy additional properties (preservation of flatness, smoothness). ~~~~(cf. e.g. \ref{regupoly}, \ref{flatandsmooth} for precise statements)
\end{enumerate}
\xtheo 

 We also define dilatations in a relative setting as follows. Let $S$ be a scheme and let $X$ be a scheme over $S$. Let $C=\{C_i\}_{i\in I}$ be closed subschemes of $S$ such that, locally, each $C_i$ is principal. Put $D= \{ C_i \times _S X\}_{i \in I}$. Let $Y=\{Y_i\}_{i\in I}$ be closed $S$-subspaces of $X$. We put $\Bl_{Y}^{C^{}} X:= \Bl_{Y}^{D^{}} X$. 
  
\prop \label{ml}  The scheme $\Bl _Y ^{C}X$ represents the contravariant functor from $Sch ^{C\text{-reg}}_S$ to $Set$ given by
\[( f: T \to S ) \mapsto \{x \in \Hom _S (T,X) |   T \times _S C_i^{} \xrightarrow{{x_|}_{C_i^{}}} X \times _S C_i^{} \text{ factors through } Y_i \times _S C_i  ~\forall i \} .\]
\xprop
Proposition \ref{ml}, proved in Proposition \ref{factintrro}, implies that for any $ T \in Sch ^{C\text{-reg}}_S$ (e.g. $T=S$ if each $C_i$ is a Cartier divisor in $S$) we have a canonical inclusion on $T$-points $ \Bl_Y ^{C^{}} X (T ) \subset X(T)$. But in general $\Bl_Y ^{C^{}} X \to X$ is not a monomorphism in the full category of $S$-schemes.

We now describe our results regarding the behaviour of dilatations of schemes endowed with a structure. So assume that $X=G$ is a monoid (resp. group, resp. Lie algebra) scheme (or any structure defined using products) over $S$.
Let $H_i\subset D_i = G|_{C_i} $ be a closed submonoid (resp. subgroup, resp. Lie subalgebra) scheme over $C_i$ for all $i \in I$ and let $H= \{H_i\}_{i \in I}$. (Recall that $C=\{C_i\}_{i\in I}$ are closed subschemes of $S$ such that, locally, each $C_i$ is principal.)
Let $\calG:=\Bl_H^{D}G\to G$ be the associated dilatation.
The structure morphism $\calG\to S$
defines an object in $Sch_S^{C\text{-}\reg}$.

\prop \label{factintro} 
Let $\calG\to S$ be the above dilatations.
\begin{enumerate}
\item Recall that the scheme $\calG\to S$ represents the contravariant functor
$Sch_S^{C\text{-}\reg}\to  Set$ given for $T\to S$ by the set of all $S$-morphisms $T\to G$ such that the induced morphism $T|_{C_i}\to G|_{C_i}$ factors through $H_i\subset G|_{C_i}$ for all $i \in I$.~~~~(\ref{Neron.blow.lemm})

\item Let $T \to S$ be an object in $Sch_S^{C\text{-}\reg}$, then as subsets of $G(T)$ \[ \calG (T) = \bigcap_{i \in I} \big( \Bl_{H_i}^{D_i}G \big)(T).~~~~(\ref{Neron.blow.lemm} )\]

\item The map $\calG\to G$ is affine. 
Its restriction over $C_i$ factors as $\calG_i\to H_i \subset  D_i$ for all $i \in I.$~~~~(\ref{Neron.blow.lemm})

\item Assume the dilatation $\calG\to S$ is flat,
 then $\calG \to S$ is equipped with the structure of a monoid (resp. group, resp. Lie algebra) scheme over $S$ such that $\calG\to G$ is a morphism of $S$-monoid (resp. $S$-group, resp. $S$-Lie algebra) schemes. ~~~~(\ref{kindof})

\item Under flatness assumptions, dilatations commute with the formation of Lie algebra schemes in a natural sense \[ \bbL ie (\Bl_H^{G|_C}G ) \cong \Bl \big\{{}^{\bbL ie (G) \times _S C_i}_{\bbL ie ( H_i )} \big\}_{i \in I} \bbL ie (G) .~~~~(\ref{LieLie} ) \]  
\end{enumerate}
\xprop

As an example of dilatations, let us explain a connexion between dilatations and Yu's famous construction of supercuspidal representations \cite{Yu01} (cf. also \cite[§10]{Yu15}).

 \exam \label{yuex} Assume in this example that $\calO $ is the ring of integers of a non-Archimedean local field $E$ and that $\pi$ is the maximal ideal of $\calO$.
  Let $G$ be a split connected reductive group scheme over $\calO$, i.e. a Demazure group scheme over $\calO$.  Let $\overset{\to}{{G}}=({{G}}^0 \subset {{G}}^1 \subset \ldots \subset {{G}}^d ={G})$ be a sequence of split Levi subgroups of ${{G}}$ over $\calO$. Put ${\mathrm{G}}_i = G_i \times _{\calO} E$ for all $i \in \{ 0 , \ldots , d \}$. Put $\overset{\to}{\mathrm{G}}=({\mathrm{G}}^0 \subset {\mathrm{G}}^1 \subset \ldots \subset {\mathrm{G}}^d =:\mathrm{G})$, a sequence of split Levi subgroups of ${\mathrm{G}}$ over $E$. For $i \in \{0 , \ldots , d\}$, let $x_i$ be the special point in the Bruhat-Tits building of ${\mathrm{G}}_i$  such that $G_i$ corresponds to $x_i$ via Bruhat-Tits theory.  Then $x_i$ comes from $x_0$ via functoriality of buildings as in \cite[§1]{Yu01}.  Let $0 \leq r_0 \leq r_1 \leq \ldots \leq r_d$ be integers. Recall that $e_G$ denotes the trivial closed subgroup scheme of $G$. There is a canonical isomorphism of groups  
  \begin{equation} \label{yuid} \overset{\to}{\mathrm{G}}(E) _{x , \overset{\to}{r}} = \Bl _{e_G, G^0 , G^1 , \ldots ,G^i , \ldots ,G^{d-1} }^{r_0 , r_1 , r_2, \ldots , r_{i+1} , \ldots , r_d} G (\calO) ,\end{equation}
  where $\overset{\to}{\mathrm{G}}(E) _{x , \overset{\to}{r}}$ is defined in \cite[§1 p584]{Yu01} and $\Bl \{^{t_i}_{H_i}\}_{i\in I} G$ denotes $\Bl \{^{\calO /\pi ^t_i}_{H_i}\}_{i\in I} G$ for integers $t_i \geq0$. We sketch the proof of the identification (\ref{yuid}) in Remark \ref{sketchproofexyu}. 
 \xexam

In the paper, we prove the following result. It generalizes the fact that congruence subgroups are normal subgroups. It was also motivated by the fact that the proof of \cite[Lemma 1.4]{Yu01}, related to the group $\overset{\to}{\mathrm{G}}(E) _{x , \overset{\to}{r}}$ appearing in Example \ref{yuex}, is not correct. In Remark \ref{yuproof}, we explain why \cite[Proof of Lemma 1.4]{Yu01} is not correct and sketch a way to correct it in the setting of Example \ref{yuex}.

\prop 
Assume that $C_i$ is a Cartier divisor in $S$ for all $i$. Assume that $G \to S$ is a flat group scheme. Let $\eta : K\to G$ be a morphism of group schemes over $S$ such that $K \to S$ is flat. Assume that $H_i  \subset G $ is a closed subgroup scheme over $S$ such that $H_i \to S$ is flat for all $i$. Assume that $\Bl_{H}^{C} G \to S$ is flat (and in particular a group scheme). Assume that, for all $i$,
 $K_{C_i} $ commutes with ${H_i}_{C_i}$ in the sense that the morphism $K_{C_i} \times _{C_i} {H_i}_{C_i} \to G_{C_i} $, $(k,h) \mapsto \eta(k)h\eta(k)^{-1} $ equals the composition morphism $K_{C_i} \times _{C_i} {H_i}_{C_i} \to {H_i}_{C_i} \subset G_{C_i}$, $(k,h) \mapsto h$. Then $K$ normalizes $\Bl _H^C G$, more precisely the solid composition map \[
 \begin{tikzcd} K \times _S \Bl _{H}^C G \ar[rr, "Id \times \Bl "] \ar[rrd, dashrightarrow] & & K \times _S G \ar[rrr, "k {,} g \mapsto \eta(k)g\eta(k)^{-1} "] & & & G\\ &  &\Bl _H ^C G  \ar[rrru, "\Bl"]&  & \end{tikzcd} \]
 factors uniquely through $\Bl _H^C G $. ~~~~~(\ref{normalizes})
 \xprop

We now discuss some applications of the above theoretical results.
 Our first application is a multi-centered congruent isomorphism.
To explain it, let $(\calO,\pi)$ be a henselian pair where $\pi\subset \calO$ is
an invertible ideal. 
\theo[Multi-centered congruent isomorphism] \label{theointro}
 Let $G$ be a separated and smooth group scheme over $S$. Let $H_0 , H_1 , \ldots , H_k$ be closed subgroup schemes of $G$ such that $H_0=e_G$ is the trivial subgroup and such that $H_i \to S$ is smooth for all $i\in \{0, \ldots , k\}$. Let $s_0 , s_1 , \ldots , s_k$ and $r_0 , r_1 , \ldots , r_k$ be in $\bbN$ such that  
 \begin{enumerate}
 \item $s_i \geq s_0 $ and $r_i \geq r_0 $ for all $i \in \{0, \ldots , k\}$
 \item $ r_i \geq s_i $ and $r_i-s_i \leq s_0$ for all $i \in \{0, \ldots , k \}$.
  \end{enumerate} Assume that $G$ is affine or $\calO $ is local. Assume that a regularity condition (RC) is satisfied (cf. Definition \ref{rc}). Then we have a canonical isomorphism of groups 
\begin{small} \[ \Bl_{H_0, H_1 , \ldots , H_k}^{s_0 ,~ s_1 ,~ \ldots , s_k} G (\calO) /  \Bl_{H_0, H_1 , \ldots , H_k}^{r_0 ,~ r_1 ,~ \ldots , r_k} G  (\calO)  \cong \Lie (\Bl_{H_0, H_1 , \ldots , H_k}^{s_0 ,~ s_1 ,~ \ldots , s_k} G )(\calO) / \Lie ( \Bl_{H_0, H_1 , \ldots , H_k}^{r_0 ,~ r_1 ,~ \ldots , r_k} G   ) (\calO) \] \end{small} where $\Bl_{H_0, \ldots , H_k}^{t_0 ,~ \ldots , t_k} G  $ denotes $ \Bl_{H_0,  ~\ldots~ , H_k}^{\calO/\pi ^{t_0}, \ldots ,  \calO /\pi^{t_k}} G $ for any $t_0 , \ldots , t_k \in \bbN$. ~~~~(\ref{isocongruent})
\xtheo 

\rema 
 Note that \cite[Lemma 1.3]{Yu01} provides a comparable "multi-centered" isomorphism, but in the framework of reductive groups over non-Archimedean local fields (\cite[Lemma 1.3]{Yu01} does not involve dilatations). Note also that Theorem \ref{theointro} extends \cite[Theorem 4.3]{MRR20}. Recall that \cite[Theorem 4.3]{MRR20} is related to the Moy-Prasad isomorphism in the setting of reductive groups over non-Archimedean local fields.  The Moy-Prasad isomorphism for reductive groups is of fundamental importance in representation theory of reductive groups over non-Archimedean local fields. It is at the heart of several constructions of supercuspidal representations (cf. e.g. \cite[Remark 3.3]{MY24}). The proof of the Moy-Prasad isomorphism in the recent reference \cite{KP22} uses \cite[Theorem 4.3]{MRR20},  cf. \cite[Theorem 13.5.1 and its proof, Proposition A.5.19 (3) and its proof]{KP22}.
\xrema 

Our other application is an interpretation of Rost double deformation space in the language of dilatations. Rost double deformation space is a fundamental tool in intersection theory and motivic homotopy theory.
Let $Z \to Y \to X$ be closed immersions of schemes (in \cite{Ro96}, all schemes are assumed to be defined over fields but we work with arbitrary schemes here). Let $\overline{D} ( X,Y,Z)$ be the double deformation space as defined in \cite[§10]{Ro96}.

Let $\bbA^2$ be $\Spec (\bbZ [s,t])$. Let $D_s,  D_{st}$ and $D_{s^2t}$ be the locally principal closed subschemes of $ \bbA ^2$ defined by the ideals $(s),(st)$ and $(s^2 t)$. We now omit the subscript ${}_{\Spec (\bbZ )}$ in fiber products.

\prop  We have a canonical identification: 
\[\overline{D} ( X,Y,Z) \cong \Bl _{ (Y \times \bbA ^2), ~ ~(Z \times \bbA ^2)}^{ (X \times D_{st}), ~(X \times D_s)}  (X \times  \bbA ^2).\]
In other words, Rost double deformation space is naturally interpreted as a double-centered dilatation. ~~~~ (\ref{rost2})
\xprop

\subsection{Structure of the paper}
Section \ref{sectionalgebr} introduces dilatations of rings. Section \ref{333} introduces dilatations of schemes and algebraic spaces. Section \ref{sec4} deals with iterated dilatations. Section \ref{sectionsinglediv} focuses on the case where the multi-center $\{ [ Z_i , D_i]\}_{i \in I}$ satisfies that  $\{ D_i\}$ are given by multiples of a single $D$. Section \ref{sec6} proves some flatness and smoothness results. Section \ref{néron} considers dilatations of monoid, group and Lie algebra schemes. Section \ref{sectioniso} studies congruent isomorphisms. Section \ref{sectionrost} interprets Rost double deformation space as dilatation.

\section{Dilatations of rings} \label{sectionalgebr}

We introduce dilatations of commutative rings. Recall that dilatations of categories also make sense (cf. \cite{Ma23c}), however dilatations of commutative rings behave specifically and it is better to treat them separately. In this paper, rings are assumed to be unital and commutative. As in \cite[\href{https://stacks.math.columbia.edu/tag/00AQ}{Tag 00AQ}]{stacks-project}, the zero ring is a ring. 

\subsection{Definition} \label{defidilaring}
Let $A$ be a unital commutative ring. If $M $ is an ideal of $A$ and $a \in A$ is an element, we say that the pair $[M,a]$ is a center in $A$. Let $I$ be an index set and let  $\{[M_i , a_i] \}_{i \in I}$ be a set of centers indexed by $I$. For $i \in I$, we put $L_i = M_i +(a_i)$, an ideal of $A$. Let $\bbN _I $ be the monoid $\bigoplus _{i \in I} \bbN $. If $\nu = (\nu_1 , \ldots , \nu_i, \ldots ) \in \bbN _I $ we put $L^{\nu}= L_1 ^{\nu_1} \cdots L_i ^{\nu _i} \cdots  $ (product of ideals of $A$) and $a^{\nu}= a_1^{\nu_1 } \cdots a_i ^{\nu _i} \cdots$ (product of elements of $A$). Note that if $\nu \in \bbN _{I}$ is such that $\nu _i =0$ for all $i$, then $L_i^{\nu _i}= L^\nu =A$. We also put $a^{\bbN _I} = \{ a^\nu | \nu \in \bbN _I \} \subset A$.

\depr \label{generaldila} The dilatation of $A$ with multi-center $\{[M_i , a_i] \}_{i \in I}$ is the unital commutative ring $A[\big\{ \frac{M_i}{a_i}\big\}_{i \in I}]$ defined as follows:

$\bullet$  The underlying set of $A[\big\{ \frac{M_i}{a_i}\big\}_{i \in I}]$ is the set of equivalence classes of symbols $\frac{m}{a^{\nu}} $ where $ \nu \in \bbN _I$ and $m \in L^{\nu}$ under the equivalence relation 
\[ \frac{m}{a^\nu} \equiv \frac{p}{a^{\lambda}} \Leftrightarrow \exists \beta \in \bbN _I \text{ such that } ma^{\beta + \lambda }= p a^{\beta + \nu } \text{ in } A.\] From now on, we abuse notation and denote a class by any of its representative $\frac{m}{a^{\nu}}$ if no confusion is likely.

$\bullet$ The addition law is given by $\frac{m}{a^{\nu}}+ \frac{p}{a^{\beta}}= \frac{m a^{\beta } + p a^{\nu}}{a^{\beta + \nu}}$.

$\bullet$ The multiplication law is given by $\frac{m}{a^{\nu}}\times  \frac{p}{a^{\beta}} = \frac{ mp}{a^{\nu + \beta }}$.

$\bullet$ The additive neutral element is $\frac{0}{1}$ and the multiplicative neutral element is $\frac{1}{1}$.
\begin{flushleft}
We have a canonical morphism of rings $A \to A[\big\{ \frac{M_i}{a_i}\big\}_{i \in I}]$ given by $a \mapsto \frac{a}{1}$. We sometimes use the notations $A[\frac{M}{a}]$ or $A[\big\{ \frac{M_i}{a_i} :i \in I \big\}]$ to denote $A[\big\{ \frac{M_i}{a_i}\big\}_{i \in I}]$.
\end{flushleft}
\xdepr

\pf Let us first prove that the relation is an equivalence relation. Assume \begin{align*}
\frac{m}{a^\nu} &\equiv \frac{p}{a^{\lambda}}  \\  \frac{l}{a^\theta} &\equiv \frac{p}{a^{\lambda}} 
\end{align*} where $\nu, \lambda , \theta \in \bbN _I$, $(m,p,l) \in( L^{\nu }, L^{\lambda }, L^{\theta})$. We want to prove that  \[ \frac{m}{a^\nu} \equiv \frac{l}{a^{\theta}}. \]By definition, there exist $ \beta , \alpha \in \bbN _I$ such that \begin{align*} ma^{\beta + \lambda }&= p a^{\beta + \nu } \\ la^{\alpha + \lambda }&= p a^{\alpha + \theta }.\end{align*} Put $\delta = \beta + \lambda + \alpha $. We get 
\[ma^{\delta + \theta  } = ma^{\beta + \lambda + \alpha + \theta } = pa^{\beta + \nu + \alpha + \theta }= la^{\beta + \nu + \alpha + \lambda } = la^{ \delta+ \nu} ,\] so  $\frac{m}{a^\nu} \equiv \frac{l}{a^{\theta}}$.  The addition and multiplication laws are associative and commutative. The distributivity axiom is satisfied and the additive neutral element is absorbent for the multiplication. So $A[\frac{M}{a}]$ is a unital commutative ring. The formula $a \mapsto \frac{a}{1}$ provides a canonical morphism of rings $A \to A[\frac{M}{a}]$.
\xpf
The element $\frac{a}{1}$  of $ A[\big\{ \frac{M_i}{a_i}\big\}_{i \in I}]$ will sometimes be denoted by $a$ if no confusion is likely.

\rema \label{MNI} Let $\{N_i\}_{i \in I}$ be ideals in $A$ such that  $ N_i +(a_i) = L_i$ for all $i \in I$.  Then we have identifications of $A$-algebras $A[\big\{ \frac{M_i}{a_i}\big\}_{i \in I}]= A[\big\{ \frac{N_i}{a_i}\big\}_{i \in I}]= A[\big\{ \frac{L_i}{a_i}\big\}_{i \in I}].$
\xrema

\rema Note that Def. Prop. \ref{defidilaring} and its proof show that if $A$ is assumed to be just a unital commutative semiring, then $A[\big\{ \frac{M_i}{a_i}\big\}_{i \in I}]$ is a unital commutative semiring and $A \to A[\big\{ \frac{M_i}{a_i}\big\}_{i \in I}]$ is a morphism of semirings. Note that most results of §\ref{defidilaring}-\ref{secpropringring} extend to semirings.
\xrema

\rema Let $\{E_i \}_{i\in I} $ be subsets of $A$, let $P_i$ be the ideal generated by $E_i$ for $ i \in I$. Then one can define $A[\{\frac{E_i}{a_i}\}_{i\in I}]$ as being $A[\{\frac{P_i}{a_i}\}_{i\in I}]$.
\xrema

\defi \label{dilaringmorphidef} Let $f:A \to B$ be a morphism of rings, we say that $f$ is a \textit{dilatation map} or an \textit{affine modification} if there exists a multi-center  $\{[M_i , a_i ]\}_{i \in I}$ in $A$ such that $B \cong  A[\big\{ \frac{M_i}{a_i}\big\}_{i \in I}]$ as $A$-algebras (cf. also Fact \ref{factdiladansloc} for another characterization).
\xdefi 

\subsection{Elementary properties of dilatations} \label{secpropringring}

We proceed with the notation from §\ref{defidilaring}.

\rema\label{enegendréparmi} As an $A$-algebra, $A[\big\{ \frac{M_i}{a_i}\big\}_{i \in I}]$ is generated by $\big\{ \frac{L_i}{a_i}\big\}_{i \in I}$. Since $L_i=M_i+(a_i)$, this implies that  $A[\big\{ \frac{M_i}{a_i}\big\}_{i \in I}]$ is generated by $\big\{ \frac{M_i}{a_i}\big\}_{i \in I}$.
\xrema
\fact \label{equivalencedef}The following assertions are equivalent.\begin{enumerate}
\item  There exists $\nu \in \bbN_I$ such that $a^\nu =0$ in $A$.
\item The ring $A[\big\{ \frac{M_i}{a_i}\big\}_{i \in I}]$ is equal to the zero ring.\end{enumerate}
\xfact
\pf 
Assume (i) holds. Let $\frac{m}{a^\beta} \in A[\frac{M}{a}]$ with $\beta \in \bbN _I$ and $m \in L^\beta$.  Then $a^\nu m=0$ in $A$ and so $\frac{m}{a^\beta}=\frac{0}{1}$ in $A[\big\{ \frac{M_i}{a_i}\big\}_{i \in I}]$. So (ii) holds. Reciprocally, assume (ii) holds. Then $\frac{1}{1}=\frac{0}{1}$ and so there exists $\nu \in \bbN _I$ such that $a^\nu =0$ in $A$. So (i) holds.
\xpf 

\fact \label{domain} Assume that $A$ is a domain and $a_i \ne 0$ for all $i$, then $A[\big\{ \frac{M_i}{a_i}\big\}_{i \in I}]$ is a domain.
\xfact 
\pf
Assume that $\frac{m}{a^\nu} \frac{l}{a^\beta} =0$ in $A[\big\{ \frac{M_i}{a_i}\big\}_{i \in I}]$. Then there exists $\theta \in \bbN _I$ such that $a^{\theta} ml=0$ in $A$. Since $a^\theta \ne 0$ and $A$ is a domain, we get that $m=0$ or $l=0$. This finishes the proof.
\xpf 

\fact \label{reduced} Assume that $A$ is reduced, then $A[\big\{ \frac{M_i}{a_i}\big\}_{i \in I}]$ is reduced.
\xfact 
\pf
Assume that, in $A[\big\{ \frac{M_i}{a_i}\big\}_{i \in I}]$,  $(\frac{m}{a^\nu})^N=0$  for some $N \in \bbN$ then there exists $\beta \in \bbN _I$ such that $a^\beta m^N=0$. We can assume that $\beta = N \theta$ with $\theta \in \bbN _I$. Then $(a^\theta m)^N =0$ in $A$ and so $a^\theta m=0$ since $A$ is reduced. So $\frac{m}{a^\nu}=0 $ in $A[\big\{ \frac{M_i}{a_i}\big\}_{i \in I}]$.
\xpf 

\fact \label{faitstakalg1}   Let $\nu  $ be in $\mathbb{N} _I$. The image of $a ^\nu$ in $A[\big\{ \frac{M_i}{a_i}\big\}_{i \in I}]$ is a non-zero-divisor.
\xfact 
\pf
 Let $b \in A[\big\{ \frac{M_i}{a_i}\big\}_{i \in I}] $ such that $a^\nu b=0$ in $A[\big\{ \frac{M_i}{a_i}\big\}_{i \in I}]$. Write $b = \frac{m}{a^\alpha}$, then we get $\frac{a^\nu m}{a^\alpha}=0$ in $A[\big\{ \frac{M_i}{a_i}\big\}_{i \in I}]$. This implies that there is $\beta \in \mathbb{N}_I$ such that $a^\beta a^\nu m =0$. So $b=0$ in $A[\big\{ \frac{M_i}{a_i}\big\}_{i \in I}]$. So $a^\nu$ is a non-zero-divisor.
 \xpf
 
 \rema
 Note that Fact \ref{faitstakalg1} holds even if  $A[\big\{ \frac{M_i}{a_i}\big\}_{i \in I}]$ is the zero ring. 
 \xrema 

\fact \label{localis}Assume $M_i = A$ for all $i \in I$. Then $A[\big\{ \frac{M_i}{a_i}\big\}_{i \in I}]=(a^{\bbN_I})^{-1} A$ where $(a^{\bbN_I})^{-1} A$ is the localization of $A$ relatively to the multiplicative monoid $a^{\bbN_I}$.
\xfact

\pf For any $\nu \in \bbN _I $, we have $L^\nu = A$. The map $\frac{x}{a^\nu} \mapsto \frac{x}{a^\nu}$ provides an isomorphism of $A$-algebras $A[\big\{ \frac{M_i}{a_i}\big\}_{i \in I}]=(a^{\bbN_I})^{-1} A$.
\xpf 

\rema \label{remalocalis} Dilatations of rings generalize entirely localizations of rings. Indeed, let $A$ be a ring and let $S$ be a multiplicative subset of $A$ (i.e. a submonoid of $A, \times$). Let $I  $ be a set such that $S= \{ s_i \}_{i \in I}$. Then $s^{\bbN _I}= S$ and Fact \ref{localis} says that $S^{-1} A = A[\big\{ \frac{A}{s}\big\}_{s \in S}]$.
\xrema

\fact \label{functosemiano} Let $f:A \to B$ be a morphism of rings. Let $\{[N_i,b_i]\}_{i \in I}$ be centers of $B$ such that $f(M_i) \subset N_i$  and $f(a_i)=b_i$ for all $i\in I$.  Then we have a canonical morphism of $A$-algebras
\begin{center} $\phi: A[\big\{ \frac{M_i}{a_i}\big\}_{i \in I}] \to B[\big\{ \frac{N_i}{b_i}\big\}_{i \in I}]$ .\end{center}
\xfact
\pf
Put $\phi (\frac{m}{a^\nu})=\frac{f(m)}{b^\nu}$.
\xpf

\fact \label{injfa} Let $\{P_i\}_{i\in I}$ be ideals of $A$ such that $P_i \subset M_i$ for all $i \in I$. Then the canonical morphism of $A$-algebras
$A[\big\{ \frac{P_i}{a_i}\big\}_{i \in I}] \to A[\big\{ \frac{M_i}{a_i}\big\}_{i \in I}]  $ is injective.

\xfact
\pf
Clear.
\xpf 

\coro \label{dilaloc}
The $A$-algebra  $A[\big\{ \frac{M_i}{a_i}\big\}_{i \in I}]$ identifies with the sub-$A$-algebra of $A[(a^{\bbN _I})^{-1}]$ generated by $\big\{ \frac{M_i}{a_i}\big\}_{i \in I}$.
\xcoro
\pf
This is an immediate corollary of Fact \ref{injfa}, Fact \ref{localis} and Remark \ref{enegendréparmi}.
\xpf

\fact \label{factdiladansloc}Let $f:A \to B$ be a morphism of rings. The following assertions are equivalent. \begin{enumerate}
\item The morphism $f$ is a dilatation map (cf. Definition \ref{dilaringmorphidef}).
\item There exists a multiplicative subset $S$ of $A$ and a sub-$A$-algebra $C$ of $S^{-1} A$ such that $B \cong C$ as $A$-algebras.
\end{enumerate}
\xfact 
\pf
The first assertion implies the second one by Corollary \ref{dilaloc}. Reciprocally, let $C$ be a sub-$A$-algebra of $S^{-1} A$. Let $I $ be the set defined as 
\[ I := \{ i= (m_i, a_i ) \in A \times S | ~{\tiny \frac{m_i}{a_i}} \text{ belongs to } C \} \subset A \times S. \]
Note that  $S= \{ a_i \} _{i \in I}  $ because $(a_i, a_i ) $ belongs to $I$ for any $a_i \in S$.
Then $C \cong A[\big\{ \frac{(m_i)}{a_i}\big\}_{i \in I}]$, indeed we have a canonical morphism of $A$-algebras $\phi: C \to A[\big\{ \frac{(m_i)}{a_i}\big\}_{i \in I}]$ sending $\frac{m_i}{a_i}$ to $\frac{m_i}{a_i}$. The morphism $\phi$ is clearly injective. The morphism $\phi$ is surjective by Remark \ref{enegendréparmi}. 
\xpf 

\rema Note that the concept of dilatations extend to categories, cf. \cite{Ma23c}. However the analog of Fact \ref{factdiladansloc} fails for categories, cf. \cite{Ma23c} for a detailed explanation.
\xrema 

\fact \label{definedsemiring} Let $c$ be a non-zero-divisor element in $A$. Then $\frac{c}{1}$ is a non-zero-divisor in $A[\big\{ \frac{M_i}{a_i}\big\}_{i \in I}]$.
\xfact 
\pf
Let $\frac{m}{a^\nu} \in A[\big\{ \frac{M_i}{a_i}\big\}_{i \in I}]$ such that $\frac{m}{a^\nu} \frac{c}{1}=0$. Then there exists $\beta \in \bbN _I $ such that $ a^\beta m c=0$ in $A$. Since $c$ is a non-zero-divisor, this implies $a^\beta m=0$ in $A$ and so $\frac{m}{a^\nu} =0$ in $A[\big\{ \frac{M_i}{a_i}\big\}_{i \in I}]$. 
\xpf

\prop \label{functooubli} \label{inutiledzdzdz}Let $K \subset I$ and put $J = I \setminus K$. Then we have a canonical morphism of $A$-algebras
\begin{center} $ \varphi :A[\big\{ \frac{M_i}{a_i}\big\}_{i \in K}] \to A[\big\{ \frac{M_i}{a_i}\big\}_{i \in I}] . $ \end{center} Moreover 
\begin{enumerate}
\item  if $M_i\subset  (a_i)$ for all $i \in J$, then $\varphi$ is surjective, and
\item if $a_i$ is a non-zero-divisor in $A$ for all $i \in J$, then $\varphi$ is injective.
\end{enumerate}
\xprop 
\pf We have a canonical injective morphism of monoids $\bbN _K\to \bbN _I$.
Let $\frac{m}{a^\nu}$ with $\nu \in \bbN _K$ and $m \in L^\nu$, then we put $\varphi (\frac{m}{a^\nu}) = \frac{m}{a^{\nu}} \in A[\big\{ \frac{M_i}{a_i}\big\}_{i \in I}]$.  We now prove the listed properties. \begin{enumerate}
\item It is enough to show that $\frac{M_i}{a_i}$ is in the image of $\varphi$ for all $i \in I$. This is obvious for all $i \in K$. So let $i \in J$ and let $\frac{m_i}{a_i} \in \frac{M_i}{a_i}$. Since $M_i \subset (a_i)$ we write $m_i= a_i x$ with $x \in A$. Then $\frac{m_i}{a_i}=\frac{x}{1}$ belongs to the image of $\varphi$. So $\varphi$ is surjective.
\item  Let $\frac{m}{a^\nu} \in  A[\big\{ \frac{M_i}{a_i}\big\}_{i \in K}]$ with $\nu \in \bbN _K$. Assume that $\frac{m}{a^\nu}=\frac{0}{1} $ in the image $A[\big\{ \frac{M_i}{a_i}\big\}_{i \in I}]$. Then there exists $\beta \in \bbN _I$ such that $ma^{\beta}=0$. Write $\beta = \nu' + \theta $ with $\nu' \in \bbN _K$ and $\theta \in \bbN _J$. Then we have $ma^{\nu'}a^{\theta}=0$ in the ring $A$ and $a^{\theta}$ is a non-zero-divisor, so $ma^{\nu'}=0$ in the ring $A$ and so $\frac{m}{a^{\nu}}=\frac{0}{1}$ in the source. So $\varphi$ is injective. \end{enumerate}
\xpf

\coro \label{1serarien} Let $K \subset I$. Assume that, for all $j \in J= I \setminus K$, the element $a_j $ belongs to the face $A^*$ of invertible elements of the monoid $(A, \times)$, i.e. $a_j$ is invertible for $\times$. Then 
\begin{center} $A[\big\{ \frac{M_i}{a_i}\big\}_{i \in I}] =A[\big\{ \frac{M_i}{a_i}\big\}_{i \in K}].$ \end{center}
\xcoro 
\pf
This follows from Proposition \ref{functooubli}.
\xpf 

\coro \label{petitinutile}Let $J \subset I$ be such that for all $j \in J$, there is $i \in I \setminus J$ satisfying that $a_i=a_j$ and $M_j \subset M_i$, then we have a canonical identification of $A$-algebras
 $ A[\big\{ \frac{M_i}{a_i}\big\}_{i \in I}]=A[\big\{ \frac{M_i}{a_i}\big\}_{i \in I\setminus J}] $.
\xcoro
\pf  Proposition \ref{functooubli} provides a morphism $ \varphi : A[\big\{ \frac{M_i}{a_i}\big\}_{i \in I\setminus J}] \to A[\big\{ \frac{M_i}{a_i}\big\}_{i \in I}]$. This morphism is injective because $a^{\bbN_{I}}= a^{\bbN_{I \setminus J }}$. The assumptions imply that $\{ \frac{M_i}{a_i} \}_{i \in I}$ belongs to the image of $\varphi$. So Remark \ref{enegendréparmi} implies that $\varphi$ is surjective. 
\xpf

\coro \label{maxsemi}  Let $\{d_i\}_{i \in I} $ be positive integers. Let $I = \coprod _{j \in J } I_j $ be a partition of $I$. Assume that, for all $j \in J $, $M_i=M_{i'}=:M_j$ and $a_i=a_{i'}=:a_j$ for all $i,i' \in I_j$.  Assume moreover that, for all $j\in J$, $\max _{i \in I_j} d_i=: d_j$ exists. 
 Then we have a canonical identification 
 \begin{center}
$ A[\big\{ \frac{M_i}{{a_i}^{d_i}}\big\}_{i \in I}] = A[\big\{ \frac{M_j}{{a_j}^{d_j}}\big\}_{j \in J}] $.
 \end{center}
\xcoro 

\pf  This follows from Proposition \ref{functooubli} and elementary arguments.
\xpf

\prop   \label{basechangeetapering}Let $K \subset I$. Then we have a canonical isomorphism of $A[\big\{ \frac{M_i}{a_i}\big\}_{i \in K}]$-algebras
\begin{center}
$A[\big\{ \frac{M_i}{a_i}\big\}_{i \in I}] = A[\big\{ \frac{M_i}{a_i}\big\}_{i \in K}] [\big\{ \frac{A[\big\{ \frac{M_i}{a_i}\big\}_{i \in K}]\frac{M_j}{1}}{\frac{a_j}{1}}\big\}_{j \in I \setminus K}], $
\end{center}
where ${A[\big\{ \frac{M_i}{a_i}\big\}_{i \in K}]\frac{M_j}{1}}$ is the ideal of $A[\big\{ \frac{M_i}{a_i}\big\}_{i \in K}]$ generated by $\frac{M_j}{1} \subset A[\big\{ \frac{M_i}{a_i}\big\}_{i \in I}]$. 
\xprop
\pf
We have a morphism of rings $A[\big\{ \frac{M_i}{a_i}\big\}_{i \in K}] \to A[\big\{ \frac{M_i}{a_i}\big\}_{i \in I}]$ given by Proposition \ref{functooubli}.
 The right-hand side of the equation in the statement of \ref{basechangeetapering} is generated as $ A[\big\{ \frac{M_i}{a_i}\big\}_{i \in K}] $-algebra by $\{ \frac{A[\big\{ \frac{M_i}{a_i}\big\}_{i \in K}]M_j}{a_j}\big\}_{j \in I \setminus K}$. 
 We now define an $ A[\big\{ \frac{M_i}{a_i}\big\}_{i \in K}]$-morphism from the right-hand side to the left-hand side sending $\frac{\frac{m_{\nu}}{a^\nu} m_j}{a_j^k}$ (with $\nu \in I$, $j \in I \setminus K$ and $k \in \bbN$) to $\frac{m_{\nu} m_j }{a^\nu a_j^k}$. This is well-defined and it is easy to check injectivity and surjectivity.
\xpf

\coro  \label{corofds}Let $S$ and $S'$ be the multiplicative monoids in $A  $ and $A[\big\{ \frac{M_i}{a_i}\big\}_{i \in I}]$ given by $\{ a^\nu | \nu \in \bbN _I\}$. Then  $S'^{-1} A[\big\{ \frac{M_i}{a_i}\big\}_{i \in I}] = S^{-1} A$.\xcoro
  \pf
Using Fact \ref{localis}, Proposition \ref{basechangeetapering} and Corollary \ref{petitinutile}, we get {\small { \begin{center}$S'^{-1} A' = A[\big\{ \frac{M_i}{{a_i}}\big\}_{i \in I}][\big\{ \frac{ A[\big\{ \frac{M_i}{{a_i}}\big\}_{i \in I}]}{{a_i}}\big\}_{i \in I}]  =  A[\big\{ \frac{M_i}{{a_i}}\big\}_{i \in I},\big\{ \frac{A}{{a_i}}\big\}_{i \in I}]=  A[\big\{ \frac{A}{{a_i}}\big\}_{i \in I}]= S^{-1} A$ .\end{center}}}\xpf

\prop \label{LiLiprime} Let $K \subset I$ be a subset. Assume that for all $i \in I \setminus K$, there exists $k(i) $ in $K$ such that $a_{k(i)} \in L_i $ and $ L_i \subset L_{k(i)}$.
Then we have a canonical identification\begin{center}$ A[\big\{ \frac{M_i}{a_i}\big\}_{i \in I}] = A[\big\{ \frac{M_k}{a_k}\big\}_{k \in K}]_{\{ \frac{a_i}{a_{k(i)}}\}_{i \in I \setminus K}}$, \end{center} where $A[\{\frac{M_k}{a_k}\}_{k \in K}]_{\{ \frac{a_i}{a_{k(i)}}\}_{i \in I \setminus K}}$ is the localization of $ A[\{\frac{M_k}{a_k}\}_{k \in K }]$ relatively to the multiplicative subset generated by ${\{ \frac{a_i}{a_{k(i)}}\}_{i \in I \setminus K}}$.
\xprop 
\pf
By Proposition \ref{functooubli}, we have a canonical map $\varphi: A[\big\{ \frac{M_k}{a_k}\big\}_{k \in K}] \to  A[\big\{ \frac{M_i}{a_i}\big\}_{i \in I}]$. For all $i \in I \setminus K$, the fraction $\frac{a_i}{a_{k(i)}}$ belongs to $A[\big\{ \frac{M_k}{a_k}\big\}_{k \in K}]$ and its image $\varphi (  \frac{a_i}{a_{k(i)}})$ is invertible in $ A[\big\{ \frac{M_i}{a_i}\big\}_{i \in I}]$, the inverse being $\frac{a_{k(i)}}{a_i}$. So we get a canonical morphism $A[\big\{ \frac{M_k}{a_k}\big\}_{k \in K}]_{\{ \frac{a_i}{a_{k(i)}}\}_{i \in I \setminus K}} \to  A[\big\{ \frac{M_i}{a_i}\big\}_{i \in I}]$. It is enough to prove that it is bijective. Given $i \in I \setminus K $ and $m \in M_i$, the identity $\frac{m}{a_i} = \frac{m}{a_{k(i)}} \frac{a_{k(i)}}{a_i}$ holds in the ring $A[\big\{ \frac{M_i}{a_i}\big\}_{i \in I}]$ and implies surjectivity. Injectivity is also easy. 
\xpf

\prop \label{interesectionalgebrik}Assume that $a_i=a_j=:b$ for all $i,j \in I$, then we have a canonical identification
\begin{center}$A[\big\{ \frac{M_i}{a_i}\big\}_{i \in I}]= A[\frac{\sum _{i\in I}M_i }{b}]$. \end{center}
\xprop
\pf
Clear by Corollary \ref{dilaloc}.
\xpf

\fact \label{faitstakalg2}Let $\nu \in \bbN _I$. We have $L^\nu A[\big\{ \frac{M_i}{a_i}\big\}_{i \in I}] = a^\nu A[\big\{ \frac{M_i}{a_i}\big\}_{i \in I}] $. 
\xfact
\pf 
Obviously $a^\nu A[\big\{ \frac{M_i}{a_i}\big\}_{i \in I}] \subset L^\nu A[\big\{ \frac{M_i}{a_i}\big\}_{i \in I}]$. Let $ y \in  L^\nu $ and $\frac{x}{a^\alpha} \in A[\big\{ \frac{M_i}{a_i}\big\}_{i \in I}]$, the formula \[ y \frac{x}{a^\alpha} = a^\nu \frac{yx}{a^{\alpha + \nu}} \] now shows that $L^\nu A[\big\{ \frac{M_i}{a_i}\big\}_{i \in I}] = a^\nu A[\big\{ \frac{M_i}{a_i}\big\}_{i \in I}] $.
\xpf
 
\prop (Universal property) \label{univpropdilarings}
If $\chi : A \to B$ is a morphism of rings such that $\chi (a_i) $ is a non-zero-divisor and generates $\chi (L_i) B$ for all $i\in I$, then there exists a unique morphism $\chi '$ of $A$-algebras $A[\big\{ \frac{M_i}{a_i}\big\}_{i \in I}] \to B$. The morphism $\chi'$ sends $\frac{l}{a^\nu} $ ($\nu \in \bbN _I , l \in L^\nu$) to the unique element $b \in B $ such that $\chi ( a^\nu) b = \chi (l)$.
\xprop
\pf
The element $b$ in the statement is unique because $\chi (a^\nu)$ is a non-zero-divisor for all $\nu \in \bbN _I$. Clearly, the map $\chi'$ defined in the statement is a morphism of $A$-algebras. Now let $\phi$ be an other morphism of $A$-algebras $A[\big\{ \frac{L_i}{a_i}\big\}_{i \in I}] \to B$. We have 
\[ \chi ( a^\nu ) \phi (\frac{l}{a^\nu})  =\phi (\frac{l}{a^\nu}) \phi ( a^\nu ) = \phi (l) = \chi (l ). \] This implies $\chi ' (\frac{l }{a^\nu}) = \phi (\frac{l}{a^\nu})$.
\xpf

Let $Ring_{A}^{a\text{-reg}}$ be the full subcategory of $A$-algebras $f:A \to B $ such that $f(a_i)$ is a non-zero-divisor for all $i \in I$. 

\rema \label{remaunicite}
Proposition \ref{univpropdilarings} implies that many of the $A$-morphisms previously described in this section are not only canonical, but also unique. For example in Proposition \ref{functooubli}, $\varphi$ is the only $A$-morphism between the considered $A$-algebras (the same comment applies for Fact \ref{functosemiano}). 
\xrema

\coro \label{univpropdilarings2}
The covariant functor from $Ring_{A}^{a\text{-reg}} $ to $ Set $  given by  \[
(f\co A \to B ) \;\longmapsto\; \begin{cases}\{*\}, \; \text{if } f(a_i) B = f (L_i ) B  \text{ for $i \in I$;}\\ \varnothing,\;\text{else}~~~~\end{cases} ~~
\] is representable by $A[\big\{ \frac{M_i}{a_i}\big\}_{i \in I}]$.
\xcoro
\pf Let $F$ be the functor defined in the statement. Let $f: A \to B$ be an $A$-algebra.
If $A[\big\{ \frac{M_i}{a_i}\big\}_{i \in I}] \to B $ is a morphism of $A$-algebras, then 
$f(a_i)B = f(L_i)B $ for all $i$ by Fact \ref{faitstakalg2}. This defines a map 
\[\Hom _{A\text{-alg}} (A[\big\{ \frac{M_i}{a_i}\big\}_{i \in I}], - ) \to F  \] of functors $Ring_{A}^{a\text{-reg}} \to Set $.
We want to show that this map is bijective when evaluated at an object $A \to B $ in $Ring_{A}^{a\text{-reg}}$. This is precisely Proposition \ref{univpropdilarings}.
\xpf 

\rema \label{remaequiano}
In Corollary \ref{univpropdilarings2}, the condition $f(a_i)B=f(L_i)B$ is equivalent to the condition $f(M_i) B \subset f(a_i) B$.
\xrema

\rema The universal property of dilatations generalizes the universal property of localizations. Indeed, let $S$ be a multiplicative subset of $A$ and let $f: A \to B$ be an $A$-algebra such that $f(s)$ is invertible for any $s \in S$. Recall that by \ref{localis}, we have $S^{-1} A = A [\{ \frac{A}{s} \}_{s \in S} ]$. Then obviously $f(s) $ is a non-zero-divisor and $f(s)$ generates $B=f(A)B$, so by the universal property of dilatations there exists a unique morphism $f' $ of $A$-algebras $S^{-1} A  \to B$. So dilatation is a construction that generalizes localization without the need to know localization. Another way to introduce dilatations of rings is to first treat the case of localizations and then to define a dilatation as in Corollary \ref{dilaloc}.
\xrema

\defi \label{defreesalg} The blowup algebra, or the Rees algebra, associated to $A$ and $\{ L_i \} _{i \in I}$ is the $\mathbb{N} _I$-graded $A$-algebra
\[ \Bl _{\{ L_i \} _{i \in I}} A = \bigoplus _{\nu \in \mathbb{N}_I} L^\nu \] where the summand ${L }^{\nu}$ is placed in degree $\nu \in \mathbb{N} _I $.
\xdefi
Let $ (e_i )_{i \in I}$ be the canonical basis of the free $\mathbb{N}$-semimodule $\mathbb{N} _I$ ($e_i$ has value $1$ in place $i$ and $0$ elsewhere). 
Recall that $a_i \in L_i$ for $i \in I$. Denote $a_{i,i}$ the element $a_i$ seen as an element of degree $e_i$  in the Rees algebra $\Bl _{\{ L_i \} _{i \in I}} A$. Let $S$ be the multiplicative subset of $\Bl _{\{ L_i \} _{i \in I}} A$ generated by $\{a_{i,i}\}_{i \in I}$. Let $\big(\Bl _{\{ L_i \} _{i \in I}} A \big)[ S^{-1}]$ be the localization of the multi-Rees algebra relatively to $S$. This $A$-algebra inherits a $\mathbb{Z}_I$-grading given, for any $l \in L^\nu $,  by 
\[ \deg ( \frac{l}{{a_{1,1}}^{\alpha _1} \ldots  {a_{i,i}}^{\alpha _i}}) = \sum _{i\in I} (\nu_i - \alpha _i ) e_i \in \bbZ _I := \bigoplus _{i \in I } \bbZ. \]

\fact \label{expliaffinealgebra} We have a canonical identification of $A$-algebras \begin{center} $A[\big\{ \frac{M_i}{a_i}\big\}_{i \in I}] = \Big[ \big( \Bl _{\{ L_i \}_{i \in I}} A \big) [S^{-1}]\Big]_{\deg=(0, \ldots , 0, \ldots)} $\end{center} where the right-hand part is obtained as degree zero elements in $\big( \Bl _{\{ L_i \}_{i \in I}} A \big) [S^{-1}]$.
\xfact

\pf This is tautological.
\xpf

\prop \label{multifinimonoring}Assume that $I = \{1 , \ldots , k \}$ is finite. Then we have a unique isomorphism of $A$-algebras
 \begin{center}$A[\big\{ \frac{M_i}{a_i}\big\}_{i \in I}]\cong A[\frac{\sum _{i \in I} ( M_i \cdot \prod _{j \in I \setminus {\{i\}}} a_j )}{a_1 \cdots a_k}] .$ \end{center}
\xprop 

\pf
Let us provide a map
\begin{center}$ \phi : A[\frac{\sum _{i \in I} ( M_i \cdot \prod _{j \in I \setminus {\{i\}}} a_j )}{a_1 \cdots a_k}] \to A[\big\{ \frac{M_i}{a_i}\big\}_{i \in I}] .$ \end{center}The ring $ A[\frac{\sum _{i \in I} ( M_i \cdot \prod _{j \in I \setminus {\{i\}}} a_j )}{a_1 \cdots a_k}]$ is generated as $A$-algebra by $\frac{\sum _{i \in I} ( M_i \cdot \prod _{j \in I \setminus {\{i\}}} a_j )}{a_1 \cdots a_k}$, we now define a map $\phi $ via (for $m_i \in M_i$, $i \in I)$:
\[\phi (\frac{\sum _{i \in I} ( m_i \cdot \prod _{j \in I \setminus {\{i\}}} a_j )}{a_1 \cdots a_k}) = \sum _{i \in I} \frac{  m_i }{a_i} .\] This is well-defined and $\phi$ is a morphism of $A$-algebras. It is easy to prove that $\phi$ is injective and surjective. Unicity follows from Proposition \ref{univpropdilarings}. 
\xpf

\rema Assume that $A= \bbZ [X,Y]$. The formal symbol $\frac{X}{3}$ does not make sense in $A[\frac{(2X)+(3Y)}{6}]$. The formal symbol $\frac{2X}{6}$ makes sense and defines an element in $A[\frac{(2X)+(3Y)}{6}]$. The formal symbols $\frac{2X}{6}$ and $\frac{X}{3} $ make sense in  $A[\frac{(X)}{3},\frac{(Y)}{2}]$ and define the same element.  The canonical isomorphism of Proposition \ref{multifinimonoring} sends $\frac{2X}{6}\in A[\frac{(2X)+(3Y)}{6}]$ to $\frac{2X}{6}= \frac{X}{3} \in A[\frac{(X)}{3},\frac{(Y)}{2}]$.

\xrema

\lemm \label{colimmmring} Write $I= \colim _{J\subset I } J $ as a filtered colimit of sets. We have a canonical identification of $A$-algebras
\begin{center}
$A[\big\{ \frac{M_i}{a_i}\big\}_{i \in I}] = \colim _{J \subset I} 
A[\big\{ \frac{M_j}{a_i}\big\}_{i \in J}] $\end{center} where the transition maps are given by Fact \ref{functooubli}.
\xlemm 
\pf
For each $J \subset I$, Fact \ref{functooubli} gives a canonical morphism $A[\big\{ \frac{M_i}{a_i}\big\}_{i \in J}] \to A[\big\{ \frac{M_i}{a_i}\big\}_{i \in I}]$ of $A$-algebras. These morphisms are compatible with transition maps. So we have a canonical $A$-morphism \begin{center}$ \phi :\colim _{J \subset I} 
A[\big\{ \frac{M_i}{a_i}\big\}_{i \in J}] \to A[\big\{ \frac{M_i}{a_i}\big\}_{i \in I}]$.\end{center}
The map $\phi$ is surjective because for any $\nu \in \bbN _I$, there exists a subset $J \subset I$ such that $\nu \in \bbN _J $ (recall that $\bbN _J\subset \bbN _I$). It is easy to check injectivity.
\xpf

\subsection{More properties of dilatations of rings}

We proceed with the notation from §\ref{defidilaring}.

\prop \label{idealsumclosed} Let $T$ be an ideal of $A$. Assume that we have a commutative diagram of $A$-algebras
\begin{tikzcd} & A/T & \\ A[\big\{ \frac{M_i}{a_i}\big\}_{i \in I}] \ar[ru, "\varphi"] & & A \ar[ll, "f"]\ar[lu, "\phi"]\end{tikzcd} where $\phi$ is the quotient map. Assume that $\phi (a_i)$ is a non-zero-divisor for all $i \in I$, i.e. that $\phi$ belongs to the category $Ring_A^{a\text{-reg}}$.
Then \begin{center}$\ker (\varphi) = \sum _{\nu \in \bbN _I} \frac{L^\nu \cap T }{a^\nu}  \subset  A[\big\{ \frac{M_i}{a_i}\big\}_{i \in I}].$\end{center}
\xprop 
\pf Let $\nu \in \bbN _I$ and $\frac{m}{a^\nu}  \in   A[\big\{ \frac{M_i}{a_i}\big\}_{i \in I}]$, i.e. $m \in L^\nu$. We have 
\[ \phi ({a^\nu} ) \varphi (\frac{m}{a^\nu}  )= \varphi (\frac{a^\nu}{1} ) \varphi (\frac{m}{a^\nu}  ) = \varphi ( \frac{m}{1})= \varphi ( f (m)) = \phi (m).\]
Now assume $ \varphi (\frac{m}{a^\nu}  ) =0$, then $\phi (m)=0$ and so $m \in L^{\nu} \cap T$. This shows that $\ker (\varphi) \subset \sum _{\nu \in \bbN _I} \frac{L^\nu \cap T }{a^\nu}$.
 Reciprocally assume $m \in L^\nu \cap T$. Then $\phi (m)=0$. This implies $\varphi (\frac{m}{a^\nu})=0$ because $\phi (a^\nu)$ is a non-zero-divisor by assumption.
\xpf 

\coro \label{coroooooo} Assume that $I=\{i\}$ is a singleton and that $a_i$ is a non-zero divisor in $A/M_i$, so that $\phi : A \to A/M_i $ belongs to $Ring _A^{a_i\text{-reg}}$. Then there exists a unique $A$-morphism $\varphi : A[\frac{M_i}{a_i}] \to A/M_i$, moreover $\ker (\varphi)$ is the ideal of $A[\frac{M_i}{a_i}]$ generated by $\frac{M_i}{a_i}$. If moreover $a_i =b_i^k$ for some $k \in \bbN, b_i \in A$, then  $A[\frac{M_i}{b_i^k}][\frac{\ker(\varphi)}{b_i^d}]= A[\frac{M_i}{b_i^{k+d}}]$ for any $d \in \bbN$.
\xcoro

\pf Since $a_i A/M_i = (a_i +M_i) A/M_i$, Proposition \ref{univpropdilarings} implies the existence and unicity of $\varphi$.
Clearly, $\frac{M_i}{a_i} \subset \ker (\varphi)$, so it is enough to prove that $\ker (\varphi )$ is included in the ideal generated by $\frac{M_i}{a_i}$ (that we denote in this proof by $\langle \frac{M_i}{a_i}\rangle$). 
Let $n \in  \bbN $, by Proposition \ref{idealsumclosed}, it is enough to prove that $ \frac{L_i^n \cap M_i}{a_i^n}$ is included in $\langle \frac{M_i}{a_i}\rangle$. An element $x \in L_i ^n $ can be written as a sum
$x= \sum _{k =0}^{n} m_k a_i^{n-k}  $ with $m_k \in M_i ^k$  (note that, if $x $ belongs to $L_i^n \cap M_i$, then $m_0 a_i^{n}$ also belongs to $M_i$). Then we assume that $x $ belongs to $L_i^n \cap M_i$, it is clear that for $k>0$ the element $\frac{ m_k a_i^{n-k} }{a_i^n} = \frac{ m_k }{a_i^k}$ belongs to $\langle \frac{M_i}{a_i}\rangle$. Now for $k=0$, using that $a_i$ is a non-zero-divisor in $A/M_i$ and that $m_0 a_i^n $ belongs to $M_i$, we get that $m_0 $ belongs to $M_i$ and it is now clear that $\frac{m_0 a_i^n }{a_i^n}$ belongs to $\langle\frac{M_i}{a_i}\rangle $. 
So $\frac{x}{a_i^n}$ belongs to $\langle \frac{M_i}{a_i}\rangle$. \\
Now we deduce the equality  $A[\frac{M_i}{b_i^k}][\frac{\ker(\varphi)}{b_i^d}]= A[\frac{M_i}{b_i^{k+d}}]$ and finish the proof:
\begin{align*}
A[\frac{M_i}{b_i^k}][\frac{\ker(\varphi)}{b_i^d}]&=A[\frac{M_i}{b_i^k}] [ \frac{\frac{M_i}{b_i^k} A[\frac{M_i}{b_i^k}]}{b_i^d} ]\\
&= A[\frac{M_i}{b_i^k}] [ \frac{{M_i} A[\frac{M_i}{b_i^k}]}{b_i^{k+d}} ]\\
\text{ by Proposition \ref{basechangeetapering}}&= A[\frac{M_i}{b_i^k}, \frac{M_i}{b_i^{k+d}}]\\
\text{by Corollary \ref{maxsemi}}&= A[\frac{M_i}{b_i^{k+d}}].
\end{align*}
\xpf

\prop \label{tenseurkernel}
Let $f:A \to B$ be an $A$-algebra. Put $N_i=f(M_i) B$ and $b_i = f(a_i)$ for $i \in I$. Then $B [\big\{ \frac{N_i}{b_i}\big\}_{i \in I}]$ is the quotient of $B \otimes _{A} A[\big\{ \frac{M_i}{a_i}\big\}_{i \in I}]$ by the ideal $T_b$ of elements annihilated by some element in $b^{\mathbb{N} _I }:= \{b^\nu | \nu \in \bbN _I \}$.
\xprop

\pf We can assume that $a_i \in M_i$ and $b_i \in N_i$. Let $B'$ be the quotient of $B \otimes _{A} A[\big\{ \frac{M_i}{a_i}\big\}_{i \in I}]$ by $T_b$. The ring map 
\begin{center} $ B \otimes _{A} A[\big\{ \frac{M_i}{a_i}\big\}_{i \in I}] \to B [\big\{ \frac{N_i}{b_i}\big\}_{i \in I}] $ \end{center} is surjective and annihilates $a^{\bbN _I}$-torsion as elements in $b^{\bbN _I}$ are non-zero-divisors in $B [\big\{ \frac{N_i}{b_i}\big\}_{i \in I}]$. Hence we obtain a surjective map $B' \to B [\big\{ \frac{N_i}{b_i}\big\}_{i \in I}]. $ To see that the kernel is trivial, we construct an inverse map.
 Namely, let $z = \frac{y}{b^{\nu}}$ be an element of $B [\big\{ \frac{N_i}{b_i}\big\}_{i \in I}]$, i.e $y \in N^{\nu}$ for some $\nu \in \bbN _I$.
  Write $y = \sum f(x_i) s_i $ with $x_i \in M^\nu$ and $s_i \in B$. We map $z$ to  the class of $\sum s_i \otimes \frac{x_i}{a^\nu}$ in $B'$. 
  This is well defined because we claim that an element of the kernel of the map $B \otimes _A M^\nu  \to N^\nu $ is annihilated by $a^\nu$
   hence maps to zero in $B'$. We now prove the claim of the previous assertion. Let $\sum_j( s_j \otimes m_j) $ be in the kernel of the map as before ($s_j \in B , m_j \in M^\nu $ for all $j$), so that $\sum_j s_j f(m_j)=0$; we have $a^\nu \sum_j (s_j \otimes m_j)= \sum_j( s_j \otimes m_j a^\nu )= \sum _j ( s_j f(m_j) \otimes a^\nu ) = (\sum _j  s_j f(m_j) ) \otimes a^\nu =0$. This finishes the proof.
\xpf 
\coro \label{flatringbase} We proceed with the notation from Proposition \ref{tenseurkernel} and assume $f:A\to B$ is flat. Then $T_b=0$, in other words we have a canonical isomorphism 
\begin{center}  $B [\{ \frac{N_i}{b_i} \}_{i \in I }] = B \otimes _A A[\{\frac{M_i}{a_i}\}_{i\in I}].$\end{center}
\xcoro
\pf
Since $f$ is flat, the map $\phi: A[\{\frac{M_i}{a_i}\}_{i\in I}] \to  B \otimes _A A[\{\frac{M_i}{a_i}\}_{i\in I}]$ is flat. Since $\phi$ is flat, the image of any non-zero-divisor element in $ A[\{\frac{M_i}{a_i}\}_{i\in I}] $ under the map $\phi$ is a non-zero-divisor. So $T_b=0$ and Proposition \ref{tenseurkernel} finishes the proof.
\xpf

\defi \label{conic} Assume that $a_i$ is a non-zero divisor for all $i \in I$. Let $J_i$ be the ideal of $A$ generated by $a_i$, it is invertible. We consider the $A$-algebra
\[ C^J_L A \defined \bigoplus _{\nu\in \mathbb{N} _I} L^\nu  \otimes J^{- \nu},\] and we call it the associated conic algebra.
\xdefi

\prop \label{conicrhoblowup} We proceed with the notation from Definition \ref{conic}. Let $\zeta$ be the ideal of $C_L^JA$ generated by elements $\rho^\alpha -1$ for $\alpha \in \bbN _I$ where $\rho_i \in C_L^JA$ is the image of $1\in A$ under
 $A \cong J_i \otimes J_i^{-1} \subset L_i \otimes J_i^{-1} \subset C_L^JA$ for any $i \in I$. We have a canonical isomorphism 
 of $A$-algebras \begin{center} $ (C_L^J A)  /\zeta \longto A[\big\{ \frac{L_i}{a_i}\big\}_{i \in I}]$. \end{center}

\xprop

\pf  Let  $t_i=a_i^\vee$ be the generator for $J_i^{-1}$, dual to $a_i$ for all $i \in \bbN_I$. We have a natural morphism of rings  given explicitly by
\[\psi: C_L^J A \longto A[\big\{ \frac{L_i}{a_i}\big\}_{i \in I}], {\sum} _{\nu \in \mathbb{N}_I} l_{\nu} \otimes t^{\nu} \mapsto {\sum}_{\nu \in \mathbb{N}_I} \frac{l_{\nu}}{a^{\nu}}.\] 
The morphism $\psi$ is surjective and $\zeta \subset \ker \psi$. It is enough to prove that $\ker \psi \subset \zeta.$ Let $X= \sum _{\nu \in c} l_{\nu} \otimes t ^{\nu}  \in C_L^J A$  where $l_\nu \in L^{\nu}$,  and $c $ is a finite subset of $\mathbb{N} _I$.
Let $\beta \in \mathbb{N}_I$ defined by $\beta _i = \mathrm{max}_{\nu \in c} \nu _i$ for all $i \in I$. Then we have 
\[ \psi (X) = \sum _{\nu \in c} \frac{l_{\nu}}{a^{\nu}}= \frac{\sum_{\nu \in c} a ^{\beta - \nu } l ^{\nu}}{a^\beta}.\]
 Assume $X \in \ker \psi$, then $\sum _{\nu \in c} a^{\beta - \nu } l ^{\nu} =0$ because $a_i$ are non-zero-divisors for all $i \in I$. So we are allowed to write
  \[ X = \big( \sum _{\nu \in c} l_{\nu} \otimes t ^{\nu} \big) - \big(\sum _{\nu \in c} a^{\beta - \nu } l ^{\nu} \big) \otimes t^{\beta} = \sum _{\nu \in c} \Big[\big(l_{\nu} \otimes t ^{\nu} \big) \big( 1- ( a ^{\beta- \nu} \otimes t ^{\beta - \nu }) \big) \Big].\]
  This finishes the proof since $ a ^{\beta- \nu} \otimes t ^{\beta - \nu } = \rho _1^{\beta _1 - \nu _1} \cdots \rho _i ^{\beta _i - \nu _i }\cdots~$.
\xpf 

\rema \label{rhomoinsunremark} We note that the ideal $\zeta $ appearing in Proposition \ref{conicrhoblowup} is in fact generated by $\{ \rho _i -1\}_{i \in I}$. To see this, use for example that in any ring and for any elements $\rho, \sigma $ in the ring we have $\rho^n-1 = (\rho-1) ( \rho^{n-1} + \ldots + \rho +1)$ and $\rho \sigma -1 = (\rho -1) (\sigma +1) + (\sigma -1) - (\rho -1)$.
\xrema

\fact \label{rfvb}Let $R$ be a ring and assume that $f:R\to A$ is a morphism of rings. Let $\{r_i\}_{i\in I}$ be elements in $R$ and assume that $a_i=f(r_i)$ for all $i\in I$. Let $R_{r}=R[ \big\{ \frac{R}{r_i} \big\}_{i\in I}]$ be the localization of $R$ at $\{r_i\}_{i\in I}$.  For any $i \in I$, let $M_i \otimes r_i^{-1}  \subset A \otimes _R R_r$. Then $
 A[\big\{ \frac{M_i}{a_i}\big\}_{i \in I}] $ identifies with the $A$-subalgebra of $A \otimes _R R_r$ generated by $\{M_i \otimes r_i^{-1} \}_{i\in I}$ and $A$.
\xfact
\pf
By  \cite[\href{https://stacks.math.columbia.edu/tag/00DK}{Tag 00DK}]{stacks-project} and Fact \ref{localis}, we get $A \otimes _R R_r = (A \otimes _R R)_r = A[\big\{ \frac{A}{a_i}\big\}_{i \in I}]$. Moreover $M_i \otimes r_i^{-1} \subset A \otimes _R R_r $ corresponds to $\frac{M_i}{a_i} \subset A[\big\{ \frac{A}{a_i}\big\}_{i \in I}]$. Now  \ref{enegendréparmi} and \ref{injfa} finish the proof.
\xpf

\rema \label{relkj}We discussed before that dilatations of rings provide a formalism unifying localizations of rings and affine blowups of rings (recall that affine blowups are studied in \cite[\href{https://stacks.math.columbia.edu/tag/052P}{Tag 052P}]{stacks-project}). It is easy to check on examples that dilatations strictly generalize localizations and affine blowups. In fact such examples already appeared in the literature in the framework of flat group schemes of finite type over discrete valuation rings, cf. e.g. \cite[Exp. VIB Ex. 13.3]{SGA3}, \cite[§7, after Proposition 7.3]{PY06}, \cite[Definitions 5.1 and 5.5]{DHdS18}. 
\xrema

\subsection{Dilatations of modules} \label{sectionmodulealg} Recall that dilatations of rings generalize localizations of rings. Recall also that localizations of $A$-modules make sense. In this remark we explain that dilatations of $A$-modules also make sense. Let $\mathrm{M}$ be an $A$-module.
The dilatation of $\mathrm{M}$ with multi-center $\{[M_i , a_i] \}_{i \in I}$ is the $A[\big\{ \frac{M_i}{a_i}\big\}_{i \in I}]$-module $\mathrm{M}[\big\{ \frac{M_i}{a_i}\big\}_{i \in I}]$ defined as follows:

$\bullet$  The underlying set of $\mathrm{M}[\big\{ \frac{M_i}{a_i}\big\}_{i \in I}]$ is the set of equivalence classes of symbols $\frac{l m}{a^{\nu}} $ where $ \nu \in \bbN _I$, $m\in \mathrm{M}$ and $l \in L^{\nu}$ under the equivalence relation 
\[ \frac{l m}{a^\nu} \equiv \frac{hp}{a^{\lambda}} \Leftrightarrow \exists \beta \in \bbN _I \text{ such that } a^{\beta + \lambda }l m=  a^{\beta + \nu }hp \text{ in } M.\] From now on, we abuse notation and denote a class by any of its representative $\frac{lm}{a^{\nu}}$ if no confusion is likely.

$\bullet$ The addition law is given by $\frac{lm}{a^{\nu}}+ \frac{hp}{a^{\beta}}= \frac{a^{\beta}lm   + a^{\nu}hp }{a^{\beta + \nu}}$.

$\bullet$ The action law is given by $\frac{l}{a^\nu} \frac{hp}{a^\beta} = \frac{lhp}{a^{\nu + \beta}}$.

$\bullet$ The additive neutral element is $\frac{0}{1}$.
\begin{flushleft}

We have a canonical morphism of $A$-modules from $\mathrm{M}$ to $\mathrm{M}[\big\{ \frac{M_i}{a_i}\big\}_{i \in I}]$ given by $m \mapsto \frac{m}{1}$. We now put $\mathrm{M} '= \mathrm{M}[\big\{ \frac{M_i}{a_i}\big\}_{i \in I}]$.
\end{flushleft}

\prop Let $\nu \in \bbN _I$. \begin{enumerate} \item Let $m \in \mathrm{M}'$. If $a^\nu m=0$, then $m=0$. 
\item  We have $a^\nu \mathrm{M}' = L^{\nu }\mathrm{M}'$.
\end{enumerate}
\xprop
\pf
\begin{enumerate}
\item Write $m = \frac{hp}{a^\beta}$. There exists $\gamma$ such that $a^\gamma a ^\nu hp=0$ in $\mathrm{M}$, so $m=0$ in $\mathrm{M}'$.
\item It is enough to prove that $a^\nu \mathrm{M}' \supset L^{\nu }\mathrm{M}'$. This follows from the identity, for $l \in L^\nu$: \[l \frac{hp}{a^\beta} = a^\nu \frac{lhp}{a^{\beta+\nu}}.\]
\end{enumerate}
\xpf
\prop (Universal Property)
 Let $A$-Mod$_{\mathrm{M}}^{a\text{-reg}}$ be the category whose objects are morphisms of $A$-modules $F: \mathrm{M} \to \mathrm{M'}$ with source $\mathrm{M}$ such that $a_i$ is a non-zero-divisor of $\mathrm{M'}$ for all $i \in I $, then 

\[ \Hom _{A\text{-Mod}_{\mathrm{M}}^{a\text{-reg}}} (\mathrm{M} [\{\frac{M_i}{a_i} \}_{i \in I}] , \mathrm{M} ' ) = \begin{cases}\{*\}, \; \text{if }L_i \mathrm{M'} = a_i \mathrm{M}' \text{ for } i \in I;\\ \varnothing,\;\text{else.}\end{cases}\]
\xprop
\pf
Similar to the proof of Proposition \ref{univpropdilarings2}.
\xpf  

Dilatations of modules enjoy similar properties than dilatations of rings, we do not list all of them here.

\rema In general $\mathrm{M} [\{\frac{M_i}{a_i}\}_{i \in I}]$ is not equal to $\mathrm{M} \otimes _A A[\{\frac{M_i}{a_i}\}_{i \in I}]$. Indeed, let $A= \bbZ [X]$, $I= \{0\}$, $M_0= (X)$, $a_0=2$ and $\mathrm{M} = A[\frac{M_0}{2}]$. Note that $A[\frac{M_0}{2} ] = \bbZ [\frac{X}{2}] \subset \bbQ [X]$. Then 
$\mathrm{M}[\frac{M_0}{2}]\cong \mathrm M$, however $\mathrm{M} \otimes _{A} A[\frac{M_0}{2} ] \not \cong \mathrm{M}$. Indeed the element $0 \ne T:=(\frac{X}{2} \otimes 1 ) + (-1 \otimes \frac{X}{2} ) \in A[\frac{M_0}{2} ] \otimes _{A} A[\frac{M_0}{2} ]$ satisfies $2T=0$ whereas $\mathrm{M}$ is $2$-torsion free. This contrasts with the case of localizations where we always have $\mathrm{M} \otimes _A S^{-1} A \cong S^{-1} \mathrm{M}$, cf. e.g. \cite[\href{https://stacks.math.columbia.edu/tag/00DK}{Tag 00DK}]{stacks-project}.
\xrema 

\rema
Let $0 \to \mathrm{M}_1 \to \mathrm{M} _2 \to \mathrm{M} _3 \to 0$ be an exact sequence of $A$-modules. Then in general $ \mathrm{M}_1[\{\frac{M_i}{a_i}\}_{i \in I}] \to \mathrm{M} _2[\{\frac{M_i}{a_i}\}_{i \in I}] \to \mathrm{M} _3[\{\frac{M_i}{a_i}\}_{i \in I}] $ is not exact. For example, take $A= \bbZ [X]$,  $I= \{0\}$, $M_0= (X)$, $a_0=2$ and consider the exact sequence $0 \to \bbZ[X] \xrightarrow{m \mapsto 2m} \bbZ[X] \to \bbZ /2\bbZ [X] \to 0$. Then $\bbZ /2\bbZ [X] [\frac{X}{2} ]=0$, however $\bbZ[\frac{X}{2}] \xrightarrow{m \mapsto 2m} \bbZ [\frac{X}{2}]$ is not surjective.  This contrasts with the case of localizations where we always have preservation of exact sequences, cf. e.g. \cite[\href{https://stacks.math.columbia.edu/tag/00CS}{Tag 00CS}]{stacks-project}.
\xrema

\section{Multi-centered dilatations in the absolute setting}  \label{333}

In this section, we define multi-centered dilatations and prove some properties. 

\subsection{Definitions}\label{blow.up.definition.sec}
Let $S$ be a scheme. An $S$-space is an $S$-algebraic space. Let us fix an $S$-space $X$. For the convenience of the reader, we recall some basic notations and well-known facts.
\nota \label{monoidclosed} \label{morphism}
Let $Clo (X)$ be the set of closed $S$-subspaces of $X$. Recall that $Clo(X)$ corresponds to quasi-coherent ideals of $\calO _X$ via \cite[\href{https://stacks.math.columbia.edu/tag/03MB}{Tag 03MB}]{stacks-project}. Let $IQCoh(\calO _X)$ denote the set of quasi-coherent ideals of $\calO _X$. It is clear that $(IQCoh(\calO _X), +, \times , 0 , \calO _X)$ is a semiring. So we obtain a semiring structure on $Clo (X)$, usually denoted by $(Clo(X), \cap , + , X, \varnothing)$. For clarity, we now recall directly operations on $Clo(X)$.
Given two closed subspaces $Y_1,Y_2$ given by ideals $\calJ _1 , \calJ_2$, their sum $Y_1+Y_2$ is defined as the closed subspace given by the ideal $\calJ_1 \calJ_2$. Moreover, if  $n \in \bbN$, we denote by $nY_1$ the $n$-th multiple of $Y_1$.  The set of locally principal closed subspaces of $X$ (cf. \cite[\href{https://stacks.math.columbia.edu/tag/083B}{Tag 083B}]{stacks-project}), denoted $Pri(X)$, forms a submonoid of $(Clo(X),+)$. Effective Cartier divisors of $X$, denoted $Car(X)$, form a submonoid of $(Pri (X),+)$. Note that $Car(X)$ is a face of $Pri(X)$. We have another monoid structure on $Clo (X)$ given by intersection, this law is denoted $\cap$. The operation $\cap$ corresponds to the sum of quasi-coherent sheaves of ideals. The set $Clo(X)$ endowed with $ \cap, +$ is a semiring whose neutral element for $+$ is $\varnothing$ and whose neutral element for $\cap$ is $X$. Let $C \in Car(X) $, a non-zero-divisor (for $+$) in the semiring $Clo(X)$. Let $Y,Y' \in Clo (X)$. If $C+Y$ is a closed subspace of $C+Y'$, then $Y $ is a closed subspace of $Y'$. Moreover if $C+Y=C+Y'$, then $Y=Y'$. Let $f:X' \to X$ be a morphism of $S$-spaces, then $f$ induces a morphism of semirings $Clo(f) : Clo (X) \to Clo (X'), Y \mapsto Y \times _X X'$, moreover $Clo(f) $ restricted to $(Pri (X),+)$ factors through $(Pri(X'),+)$, this morphism of monoids is denoted $Pri (f)$. Let $Y_1,Y_2 \in Clo(X)$, we write $Y_1 \subset Y_2$ if $Y_1 $ is a closed subspace of $Y_2$. We obtain a poset $(Clo(X), \subset)$. Let $Y_1,Y_2,Y_3 \in Clo(X)$, if $Y_1 \subset Y_2 $ and $Y_1 \subset Y_3$ then $Y_1 \subset Y_2 \cap Y_3$. Let $Y_1,Y_2 \in Clo (X)$, then $(Y_1 \cap Y_2) \subset Y_1$ and $Y_1 \subset ( Y_1+Y_2)$.
Finally, if $Y=\{ Y_e \} _{e \in E}$ is a subset of $Clo (X)$ and if $\nu \in \bbN ^{E}$, we put $Y^\nu = \{\nu _e Y_e  \}_{e \in E}$ and if moreover $\nu \in \bbN _E$, we put $\nu Y = \sum _{e \in E } \nu _e Y_e$.\xnota
\rema
Be careful that the operation $+$ on $Clo(X)$ is not the operation $\cup$ of \cite[\href{https://stacks.math.columbia.edu/tag/0C4H}{Tag 0C4H}]{stacks-project}. Recall that $+$ corresponds to multiplication of ideals whereas $\cup$ corresponds to intersection of ideals.
\xrema
\rema We proceed with the notation from Notation \ref{morphism}. In general the image of the map $Pri(f) |_{Car(X)}$ is not included in $Car (X').$ 
\xrema

\defi \label{defcatreg} Let $D= \{ D_i \} _{i \in I}$ be a subset of $Clo (X)$.  \begin{enumerate} \item Let $\Spac _X ^{D\text{-reg}}$ be the category of $S$-algebraic spaces $f:T \to X $ over $X$ such that for any $i \in I$, $T \times _X D_i $ is a Cartier divisor in $T$. 
\item If $X=S$ is a scheme, let $Sch _X ^{D\text{-reg}}$ be the category of $X$-schemes $f:T \to X $  such that for any $i \in I$, $T \times _X D_i $ is a Cartier divisor in $T$. 
\end{enumerate}
\xdefi

If $T'\to T$ is flat and $T\to X$ is an object in $\Spac_X^{D\text{-}\reg}$ or $Sch _X ^{D\text{-reg}}$, so is the composition $T'\to T\to X$ by  \cite[\href{https://stacks.math.columbia.edu/tag/083Z}{Tag 083Z}]{stacks-project} and \cite[\href{https://stacks.math.columbia.edu/tag/02OO}{Tag 02OO}]{stacks-project}.

\fact \label{fact3.3} \label{product}  Let $D= \{ D_i \} _{i \in I}$ be a subset of $Clo (X)$. 
 \begin{enumerate} \item Let $f:T\to X$ be an object in $\Spac _X ^{D\text{-reg}}$. Then for any $\nu \in \bbN _I$, the space $T \times _X \nu D$ is a Cartier divisor in $T$, namely $\nu ( T \times _X D)$.
\item Assume $\{ D_i \} _{i \in I}$ is made of some multiples of a finite number of locally principal closed subschemes $D_1', \ldots, D_k'$ (e.g. $\#I$ is finite). Then $\Spac _X ^{\{ D_i \} _{i \in I}\text{-reg}}$ equals $\Spac _X^{D_1'+ \ldots + D_k'}$.  Moreover finite products exist in $\Spac_X^{\{ D_i \} _{i \in I}\text{-reg}}$ (and preserve schemes).
\end{enumerate}
\xfact

\pf \begin{enumerate} \item 
This follows from the fact that $(Clo(f),+)$ is a morphism of monoids and the fact that $Car (T)$ is a submonoid of $Clo(T)$ (e.g. cf. the discussion in Notation \ref{morphism}).
\item The first assertion follows from the fact that $Car(X)$ is a face of the monoid $Pri (X)$ (e.g. cf. the discussion in Notation \ref{morphism}). 
 In case of schemes, the second assertion now follows from \cite[§3. p741]{MRR20} (in case of algebraic spaces use the same argument and  \cite[\href{https://stacks.math.columbia.edu/tag/085U}{Tag 085U}]{stacks-project}). \end{enumerate}
\xpf

\defi \label{multi}
 A multi-center in $X$ is a set $\{ [Y_i , D_i ]\} _{i \in I} $ such that \begin{enumerate}
 \item  $Y_i $ and $ D_i$ belong to $ Clo (X) $,
 \item there exists an affine étale covering $\{ U_{\gamma } \to X \}_{\gamma \in \Gamma} $ of $X$ such that $D_i |_{U_\gamma}$ is principal for all $i \in I $ and $\gamma \in \Gamma $ (in particular $D_i $ belongs to $Pri (X)$ for all $i$). \end{enumerate}
 In other words a multi-center  $\{ [Y_i , D_i ]\} _{i \in I} $ is a set of pairs of closed $S$-spaces such that locally each $D_i$ is principal.
\xdefi

\rema Let $\{ Y_i , D_i \} _{i \in I} $ such that   $Y_i \in Clo (X) $ and $D_i \in Pri (X)$ for any $i \in I$. Assume that $I$ is finite, then $\{ [Y_i , D_i ]\} _{i \in I} $ is a multi-center in $X$, i.e. the second condition in Definition \ref{multi} is satisfied.
\xrema 

We now fix a multi-center $\{ [Y_i , D_i ]\} _{i \in I} $ in $X$.
 Denote by
$\calM_i ,  \calJ_i$ the quasi-coherent sheaves
of ideals in $\calO_X$ so that $Y_i=V(\calM _i), V(\calJ _i)=D_i$. We put $Z_i = Y_i \cap D_i$ and $\calL_i = \calM _i + \calJ _i $ so that $Z_i = V ( \calL _i )$ for any $i \in I $. We put $Y= \{Y_i \}_{i \in I}$, $D= \{D_i \}_{i \in I}$ and $Z= \{ Z_i \}_{i \in I}$.
We now introduce dilatations $\calO _X$-algebras by glueing (cf. \cite[\href{https://stacks.math.columbia.edu/tag/04TP}{Tag 04TP}]{stacks-project}). 

\depr \label{defiqcoalg} \begin{sloppypar}
The dilatation of $\calO_X$ with multi-center $\{[\calM _i , \calJ _i ] \}_{i \in I}$ is the quasi-coherent $\calO _X$-algebra  $\calO _X \Big[\Big\{ \frac{\calM_i}{\calJ_i}\Big\}_{i \in I} \Big]$ obtained by glueing as follows.  The quasi-coherent $\calO _X$-algebra $\calO _X \Big[\Big\{ \frac{\calM_i}{\calJ_i}\Big\}_{i \in I} \Big]$ is characterized by the fact that its restriction, on any étale $S$-morphism $\varphi: U \to X$ such that $U$ is an affine scheme and each $D_i$ is principal on $U$ and generated by $a_{iU}$, is given by
\[ \Big( \calO _X \Big[\Big\{ \frac{\calM_i}{\calJ_i}\Big\}_{i \in I} \Big] \Big) _{\big|_{U_{}}}  =\widetilde{ \Gamma (U , \calO _X ) \Big[\Big\{ \frac{\Gamma (U, \calM _i) }{a_{i U} }\Big\}_{i \in I}} \Big] \]
where $\widetilde{~~~~~~~~~}$ is given by \cite[\href{https://stacks.math.columbia.edu/tag/01I7}{Tag 01I7}]{stacks-project} and \cite[\href{https://stacks.math.columbia.edu/tag/03DT}{Tag 03DT}]{stacks-project} (we work with small étale sites). \end{sloppypar}
\xdepr
\pf By Definition \ref{multi}, the affine schemes $U \to X$ satisfying the conditions in the statement form an étale covering of $X$.
Now Proposition \ref{defiqcoalg} follows from \cite[\href{https://stacks.math.columbia.edu/tag/03M0}{Tag 03M0}]{stacks-project}, \cite[\href{https://stacks.math.columbia.edu/tag/04TR}{Tag 04TR}]{stacks-project} and Corollary \ref{flatringbase}.
\xpf
 \begin{sloppypar} 
Let $\Bl_{\{  \calL_ i \} _{i \in I}}\calO_X= \bigoplus _{ \nu \in \mathbb{N} _I}\calL ^\nu $ denote the multi-Rees algebra, it is a quasi-coherent $\bbN _I$-graded $\calO_X$-algebra. 
 By localization, we get a quasi-coherent $\calO _X$-algebra  $\big(\Bl_{\{  \calL_ i \} _{i \in I}}\calO_X \big) [\{ \calJ _i^{-1}\}_{i \in I} ]$ (locally, we invert  a generator of $\calJ _i$, for each $i \in I$). This $\calO _X$-algebra inherits a grading giving local generators of $\calJ _i$ degree $e_i$. \end{sloppypar} 
 \fact \label{blow.up.algebra.defi} \begin{sloppypar}
  We have a canonical identification of quasi-coherent $\calO _X$-algebras 
\[ \calO _X \Big[\Big\{ \frac{\calM_i}{\calJ_i}\Big\}_{i \in I} \Big] = \Big[ \big(\Bl_{\{  \calL_ i \} _{i \in I}}\calO_X \big) [\{ \calJ _i^{-1}\}_{i \in I} ]\Big]_{\deg=(0,\ldots,0, \ldots)},
\] where the right-hand side is
obtained as the subsheaf of degree zero elements in $ \big(\Bl_{\{  \calL_ i \} _{i \in I}}\calO_X \big) [\{ \calJ _i^{-1}\}_{i \in I} ]$. In particular $ \calO _X \Big[\Big\{ \frac{\calM_i}{\calJ_i}\Big\}_{i \in I} \Big]=  \calO _X \Big[\Big\{ \frac{\calL_i}{\calJ_i}\Big\}_{i \in I} \Big]$. \end{sloppypar}
\xfact
\pf
This follows from Fact \ref{expliaffinealgebra}.
\xpf

\defi  \label{defimultidilaalgsp}
The dilatation of $X$  with multi-center $\{[Y_i , D_i ]\} _{i \in I} $ is the $X$-affine algebraic space over $S$
\[
\Bl_Y^DX\defined \Spec _X \big( \calO _X \Big[\Big\{ \frac{\calM_i}{\calJ_i}\Big\}_{i \in I} \Big] \big).
\]
\xdefi

\rema \label{remaYZYZYZ}
Fact \ref{blow.up.algebra.defi} implies that $\Bl _Y ^D X = \Bl _Z ^D X$.
\xrema

\fact \label{spschcase}Assume that $X=S$ is a scheme. Then $\Bl _Y ^D X$ is a scheme.
\xfact

\pf We have an affine morphism $\Bl _Y ^D X \to X$, now the fact follows from \cite[\href{https://stacks.math.columbia.edu/tag/03WG}{Tag 03WG}]{stacks-project}.
\xpf

\nota \label{notabl}We will also use the notation  $ \Bl \big\{ _{Y_i}^{D_i } \big\}_{i\in I} X$ and $\Bl _{ \{Y_i\}_{i\in I}}^{\{D_i\}_{i\in I}} X $ to denote $\Bl _Y ^D X $. If $I= \{i\}$ is a singleton  we also use the notation $\Bl _{Y_i} ^{D_i} X $. If $K \subset I$, we sometimes use the notation $\Bl _{ \{Y_i\}_{i\in K}, \{Y_i\}_{i\in I \setminus K}}^{\{D_i\}_{i\in K},\{D_i\}_{i\in I \setminus K} } X $. If $I = \{1 , \ldots ,k\}$, we use the notation $\Bl _{Y_1, \ldots , Y_k}^{D_1 , \ldots , D_k} X$. Etc. \xnota 

\defi \label{dilamapdee}Let $X$ be a scheme or an algebraic space over a scheme $S$. We say that a morphism $f:X' \to X$ is a dilatation morphism if $f  $ is equal to $ \Bl \big\{ _{Y_i}^{D_i } \big\}_{i\in I} X \to X$ for some multi-center $\{[Y_i , D_i]\}_{i \in I}.$ The terminologies affine blowups and affine modifications are also used. \xdefi

\fact If $Y_i = \varnothing$ is the empty closed subscheme defined by the ideal $\calO _X$ for all $i \in I $, then we say that  $\Bl \big\{ _{\varnothing}^{D_i } \big\}_{i\in I} X \to X$ is a localization. Moreover if $\# I$ is finite, the dilatation morphism $ \Bl \big\{ _{\varnothing}^{D_i } \big\}_{i\in I} X \to X$ is an open immersion. 
\xfact 
\pf
This is local on $X$ \cite[\href{https://stacks.math.columbia.edu/tag/03M4}{Tag 03M4}]{stacks-project} and follows from Fact \ref{localis}. 
\xpf

\subsection{Exceptional divisors}

We proceed with the notation from \S\ref{blow.up.definition.sec}.

\prop \label{blow.up.Cartier.lemm}
As closed subspaces of $\Bl_Y^DX$, one has, for all $\nu \in \mathbb{N}_I$,
\[
\Bl_Y^DX\x_X \nu Z=\Bl_Y^DX\x_X \nu D,
\]
which is an effective Cartier divisor on $\Bl_Y^DX$.
\xprop
\pf
Our claim is étale local on $X$. 
We reduce to the affine case and apply \ref{faitstakalg1} and \ref{faitstakalg2}.
\xpf

\subsection{Universal property}\label{blow.up.univ.ppty.section}
We proceed with the notation from \S\ref{blow.up.definition.sec}.
As $\Bl_Y^DX\to X$ defines an object in $\Spac_X^{D\text{-}\reg}$ by Proposition \ref{blow.up.Cartier.lemm}, the contravariant functor
\begin{equation}\label{blow.up.represent.eq}
\Spac_X^{D \text{-}\reg}\to Set, \;\;\;(T\to X)\mapsto \Hom_{X\text{-Spaces}}\big(T,\Bl_Y^DX\big)
\end{equation}
together with $\id_{\Bl_Y^DX}$ determines $\Bl_Y^DX\to X$ uniquely up to unique isomorphism.
The next proposition gives the universal property of dilatations.

\prop\label{blow.up.rep.prop}
The dilatation $\Bl_Y^DX\to X$ represents the contravariant functor $\Spac_X^{D\text{-}\reg}\to Set$ given by
\begin{equation}\label{blow.up.iso.eq}
(f\co T\to X) \;\longmapsto\; \begin{cases}\{*\}, \; \text{if $f|_{T\x_XD_i}$ factors through $Y_i\subset X$ for $i \in I$;}\\ \varnothing,\;\text{else.}\end{cases}
\end{equation}
\xprop
\pf Note that  the condition $f|_{T\x_XD_i}$ factors through $Y_i\subset X$ is equivalent to the condition $f|_{T\x_XD_i}$ factors through $Z_i\subset X$, because $Z_i = Y_i \cap D_i$.
Let $F$ be the functor defined by \eqref{blow.up.iso.eq}. 
If $T\to \Bl_Y^DX$ is a map of $X$-spaces, then the structure map $T\to X$ restricted to $T\x_XD_i$ factors through $Z_i\subset X$ by Proposition \ref{blow.up.Cartier.lemm}.
This defines a map
\begin{equation}\label{blow.up.map.sheaves.eq}
\Hom_{X\text{-Spaces}}\big(\,\str\,,\Bl_Y^DX\big)\;\longto\; F
\end{equation}
of contravariant functors $\Spac_X^{D\text{-}\reg}\to Set$.
We want to show that \eqref{blow.up.map.sheaves.eq} is bijective when evaluated at an object $T\to X$ in $\Spac_X^{D\text{-}\reg}$.
As \eqref{blow.up.map.sheaves.eq} is a morphism of étale sheaves, we reduce to the case where both $X$ and  $T$ are affine and $J_i $ is principal for all $i \in I$.
Now Proposition \ref{univpropdilarings} finishes the proof.
\xpf

\prop \label{defined} Put $f : \Bl_Y ^D X \to X$. Then the morphism of monoids $Clo (f) |_{Car(X)} $ factors through $Car( \Bl_Y ^D X )$. In other words, any effective Cartier divisor $C \subset X$ is defined for $f$, i.e. the fiber product $C \times _X \Bl_{Y}^D X \subset \Bl_{Y}^D X $ is an effective Cartier divisor (cf.  \cite[\href{https://stacks.math.columbia.edu/tag/01WV}{Tag 01WV}]{stacks-project}).\xprop

\pf  We reduce to the case where $X= \Spec (A)$ is affine and apply Fact \ref{definedsemiring}.
\xpf

\prop \label{uniquemorphism}Let $J $ be a subset of $I$.
There exists a unique $X$-morphism \[\varphi: \Bl_{\{Y_j\}_{j \in I} }^{\{D_j\}_{j \in I}} X \to \Bl_{\{Y_j\}_{j \in J} }^{\{D_j\}_{j \in J}} X .\]
\xprop   
\pf This follows from
Propositions \ref{blow.up.rep.prop} and  \ref{blow.up.Cartier.lemm}.
\xpf 
\rema Proposition \ref{uniquemorphism} is the spaces version of Proposition \ref{functooubli} (combined with Remark \ref{remaunicite}).
\xrema

\rema Proposition \ref{multietape} will refine Proposition \ref{uniquemorphism} and show that $\varphi $ is in fact a dilatation map (cf. Definition \ref{dilamapdee}).
\xrema 

\prop \label{inutilgeo} Let $K \subset I$ and assume $Z_i = D_i $ is a Cartier divisor in $X$, for all $i \in K$. Then 
\[ \Bl _{\{Z_j\} _{j \in I }}^{ \{D_j\} _{j \in I }} X = \Bl _{ \{Z_j\} _{j \in I \setminus K}}^{ \{D_j\} _{j \in I \setminus K}} X  .\]
\xprop 

\pf 
Both sides belong to $\Spac_X^{D \text{-}\reg}$ by Propositions \ref{defined} and \ref{blow.up.Cartier.lemm}, so it is enough to show that they agree when evaluated  at any $f:T \to X \in \Spac_X^{D\text{-}\reg}$. We have 
\begin{align*}
\Bl _{\{Z_j\} _{j \in I }}^{ \{D_j\} _{j \in I }} X(T) &= \begin{cases}\{*\}, \; \text{if $f|_{T\x_XD_i}$ factors through $Z_i\subset X$ for $i \in I$;}\\ \varnothing,\;\text{else.}\end{cases} \\&= \begin{cases}\{*\}, \; \text{if $f|_{T\x_XD_i}$ factors through $Z_i\subset X$ for $i \in I\setminus K$;}\\ \varnothing,\;\text{else.}\end{cases} \\
&= \Bl _{ \{Z_j\} _{j \in I \setminus K}}^{ \{D_j\} _{j \in I  \setminus K}} X  (T).
\end{align*}
\xpf 

\rema  Proposition \ref{inutilgeo} is the spaces version of Proposition \ref{1serarien}.
\xrema

\prop \label{inter} Assume $D_i=D_j=:D$ for all $i,j \in I$.  Then 
\[ \Bl \big\{^D_{Y_i}\big\}_{i\in I} X = \Bl _{\cap _{i\in I} Y_i }^{D} X .\]
\xprop 
\pf
Both sides belong to $\Spac_X^{D\text{-}\reg}$ by Proposition \ref{blow.up.Cartier.lemm}, so it is enough to show that they agree when evaluated  at any $f:T \to X \in \Spac_X^{D\text{-}\reg}$. We have 
\begin{align*}
\Bl \big\{^D_{Y_i}\big\}_{i\in I} X  (T) &= \begin{cases}\{*\}, \; \text{if $f|_{T\x_XD}$ factors through $Z_i\subset X$ for $i \in I$;}\\ \varnothing,\;\text{else.}\end{cases} \\
&= \begin{cases}\{*\}, \; \text{if $f|_{T\x_XD}$ factors through ${\cap _{i\in I} Z_i }\subset X$ ;}\\ \varnothing,\;\text{else.}\end{cases} \\
&= \Bl _{\cap _{i\in I} Y_i}^D X (T).
\end{align*}
\xpf

\rema  Proposition \ref{inter} is the spaces version of Proposition \ref{interesectionalgebrik}.
\xrema

\fact \label{maxxxxx}   Let $I = \coprod _{j \in J } I_j $ be a partition of $I$. Assume that, for all $j \in J $, $Y_i=Y_{i'}=:Y_j$ and $D_i= D_i'$ for all $i,i' \in I_j$.  Let $\nu $ be in $\bbN ^I$ and assume that for all $j \in J$, the number $\max _{i \in I_j} \nu_i$ exists. Let $\nu  \in \bbN ^J$ defined by $\nu _j=\max _{i \in I_j} \nu_i$. Then 

\[ \Bl \big\{^{\nu _j D_j}_{Y_j}  \big\}_{j \in J } X = \Bl \big\{^{\nu _i D_i}_{Y_i} \big\}_{i \in I } X .\]
\xfact 

\pf We have a canonical morphism $\Bl \big\{^{\nu _j D_j}_{Y_j}  \big\}_{j \in J } X \to \Bl \big\{^{\nu _i D_i}_{Y_i} \big\}_{i \in I } X$ by Proposition \ref{uniquemorphism}, now we reduce to the affine case and apply Corollary \ref{maxsemi}.
\xpf

\rema  Proposition \ref{maxxxxx} is the spaces version of Corollary \ref{maxsemi}.
\xrema

\prop \label{multietape} Let $J $ be a subset of $I$ and put $K = I \setminus J$. Then
\[ \Bl \big\{_{Y_i} ^{D_i}\big\}_{i \in I} X = \Bl \Big\{_{Y_k \times _X  \Bl \big\{_{Y_i}^{D_i}\big\}_{i \in J} X } ^{
D_k \times _X  \Bl \big\{_{Y_i}^{D_i}\big\}_{i \in J} X  }\Big\} _{k \in K } \Bl \big\{_{Y_i}^{D_i}\big\}_{i \in J} X  .\] This in particular gives the unique $X$-morphism \[\Bl \big\{_{Y_i}^{D_i}\big\}_{i \in I} X  \to \Bl \big\{_{Y_i}^{D_i}\big\}_{i \in J} X  \] of Proposition \ref{uniquemorphism}. 
\xprop 

\pf The right hand side is well-defined (e.g. cf. \cite[\href{https://stacks.math.columbia.edu/tag/053P}{Tag 053P}]{stacks-project}). Using Proposition \ref{blow.up.Cartier.lemm} and Proposition \ref{defined}, one obtains that the right hand side is in $\Spac_X^{D \text{-}\reg }$. So it is enough to see that both sides coincide when evaluated at any $f:T\to X \in \Spac_X^{D \text{-}\reg }$. This follows from Proposition \ref{blow.up.rep.prop}. 
\xpf

\prop \label{ZZprimeopen}
 Let $K \subset I$ be such that \begin{enumerate} \item  $I \setminus K $ is finite, \item for all $i \in I \setminus K $, there exists $k(i) \in K $ such that $Z_{k(i)} \subset Z_i$ and $Z_i\subset D_{k(i)}$. \end{enumerate} Then the canonical morphism given by Proposition \ref{uniquemorphism}
 \[ \Bl \big\{_{Y_i}^{D_i}\big\}_{i \in I} X  \to \Bl \big\{_{Y_i}^{D_i}\big\}_{i \in K} X  \] is an open immersion. 
\xprop 

\pf
We reduce to the local case and apply Proposition \ref{LiLiprime}.
\xpf 

\subsection{Universal property in the relative setting} 

In the relative setting, Proposition \ref{blow.up.rep.prop} implies the following statement. Let $S$ be a scheme and let $X$ be a space over $S$. Let $C=\{C_i\}_{i\in I}$ be closed subspaces of $S$ such that, locally, each $C_i$ is principal. Put $D= \{ C_i \times _S X\}_{i \in I}$. Let $Y=\{Y_i\}_{i\in I}$ be closed $S$-subspaces of $X$. We put $\Bl_{Y}^{C^{}} X:= \Bl_{Y}^{D^{}} X$. 
  
\prop  \label{factintrro} The space $\Bl _Y ^{C}X$ represents the contravariant functor from $Sch ^{C\text{-reg}}_S$ to $Set$ given by
\[( f: T \to S ) \mapsto \{x \in \Hom _S (T,X) | T \times _S C_i^{} \xrightarrow{{x_|}_{C_i^{}} } X \times _S C_i^{} \text{ factors through } Y_i \times _S C_i  ~ \forall i\}. \]
\xprop 
\pf
Let $T \to S$ be an object in the category $Sch_S^{C\text{-reg}}$.
Consider the map $\theta: \Hom _S ( T , \Bl_Y ^C X ) \to \Hom _S (T, X)$ which sends a morphism to its composition with the dilatation map $\Bl _Y ^C X \to X $.  Note that any $S$-morphism $T \to X$ belongs to the category $Spaces_X^{\{X_{C_i}\}_{i \in I}\text{-reg}}.$ Using Proposition \ref{blow.up.rep.prop}, we see that $\theta $ is injective and it is easy to prove that the image of $\theta$ is \[\{x \in \Hom _S (T,X) |~ {x_|}_{C_i^{}} : T \times _S C_i^{} \to X \times _S C_i^{} \text{ factors through } Y_i \times _S C_i ~ \forall i \}.\] This finishes the proof. 
\xpf 

\rema \label{remainj} Proposition \ref{factintrro} implies that for any $ T \in Sch ^{C\text{-reg}}_S$ (e.g. $T=S$ if each $C_i$ is a Cartier divisor in $S$) we have a canonical inclusion on $T$-points $ \Bl_Y ^{C^{}} X (T ) \subset X(T)$. But in general $\Bl_Y ^{C^{}} X \to X$ is not a monomorphism in the full category of $S$-spaces.
\xrema

\subsection{Multi-centered dilatations and mono-centered dilatations} \label{sectionmonopoly} We proceed with the notation from \S\ref{blow.up.definition.sec}.

\prop \label{colimitspacegeneral} Write $I = \colim _{J \subset I} J $ as a filtered colimit of sets where transition maps are given by inclusions of subsets. We have a canonical identification 
\begin{center}$\Bl \big\{^{D_i}_{Y_i}  \big\}_{i \in I}  X = \lim_{J \subset I} \Bl \big\{^{D_i}_{Y_i}  \big\}_{i \in J} X $\end{center} where transition maps are described in Propositions \ref{uniquemorphism} and \ref{multietape}. On the right-hand side the direct limit is in the category of $S$-spaces over $X$.
\xprop 
\pf By \cite[\href{https://stacks.math.columbia.edu/tag/07SF}{Tag 07SF}]{stacks-project} the limit exists.
For each $J \subset I$, Propositions \ref{uniquemorphism} and \ref{multietape} give us a $X$-morphism $ \Bl \big\{^{D_i}_{Y_i}  \big\}_{i \in I}  X \to  \Bl \big\{^{D_i}_{Y_i}  \big\}_{i \in J} X$, so we get an $X$-morphism $\phi : \Bl \big\{^{D_i}_{Y_i}  \big\}_{i \in I}  X \to \lim_{J \subset I} \Bl \big\{^{D_i}_{Y_i}  \big\}_{i \in J} X $. To prove that $\phi$ is an isomorphism, we reduce to the affine case where the result follows from Proposition \ref{colimmmring}.
\xpf

\prop \label{multisingle} 
Assume that $\#I=k$ is finite. We fix an arbitrary bijection $I = \{ 1 ,\ldots , k\}$.  We have a canonical isomorphism of $X$-spaces
\[ \Bl _{\{Y_i\}_{i \in I}}^{\{D_i\}_{i \in I} }X  \cong \Bl ^{(\Bl\cdots ) \times _X D_k}_{( \Bl \cdots ) \times _X Y_k }\Biggl(\cdots \Bl_{(\Bl \cdots ) \times _X Y_3}^{(\Bl\cdots ) \times _X D_3}\biggl( \Bl _{(\Bl_{Y_1}^{D_1}X) \times _{X} Y_2} ^{(\Bl_{Y_1}^{D_1}X) \times _{X} D_2} \bigl(\Bl _{Y_1}^{D_1} X \bigl) \biggl) \Biggl). \]
\xprop
\pf  By induction on $k$ using  Proposition \ref{multietape}.
\xpf

\prop[Monopoly isomorphism] \label{formulamultimono} 
Assume that $\#I=k$ is finite. We fix an arbitrary bijection $I = \{ 1 ,\ldots , k\}$. We have a unique isomorphism of $X$-spaces
 \[\Bl _{\{Y_i\}_{i \in I}}^{\{D_i\}_{i \in I} }X  \cong \Bl _{\bigcap _{i \in I }( Y_i + D_1+ \ldots + D_{i-1} + D_{i+1} + \ldots + D_k)} ^{D_1+ \ldots + D_k} X.\]
\xprop 
\pf  \begin{sloppypar}
Since $Car(X)$ is a face of the monoid $Pri(X)$, the right-hand side belongs to $\Spac^{D\text{-reg}}_X$. Let $f  :\Bl _{\bigcap _{i \in I } Z_i + D_1+ \ldots + D_{i-1} + D_{i+1} + \ldots + D_k} ^{D_1+ \ldots + D_k} X \to X$ be the dilatation morphism. Let us prove that $f^{-1} (D_i) \subset f^{-1} (Y_i) $ for all $i \in \{1 , \ldots , k \}$. By Proposition \ref{blow.up.rep.prop}, \[ f^{-1} (D_1 + \ldots + D_k) \subset f^{-1} (\bigcap _{i \in I } Y_i + D_1+ \ldots + D_{i-1} + D_{i+1} + \ldots + D_k)\] is a Cartier divisor in $f^{-1} (X)$.  Moreover by the discussion in Notation \ref{monoidclosed}
\begin{align*}
&f^{-1} (\bigcap _{i \in I } Y_i + D_1+ \ldots + D_{i-1} + D_{i+1} + \ldots + D_k) \\
=& \bigcap _{i \in I }  f^{-1}(Y_i) +  f^{-1}(D_1) + \ldots +  f^{-1}(D_{i-1}) +  f^{-1}(D_{i+1}) + \ldots +  f^{-1}(D_k)\text{, and }
\end{align*}
\[
f^{-1} (D_1 + \ldots + D_k ) = f^{-1} (D_1) + \ldots + f^{-1} (D_k).\]
So for any $i \in \{1 , \ldots , k \}$, we have \begin{small} \[f^{-1} (D_1)+ \ldots + f^{-1} (D_k  ) \subset  f^{-1}(Y_i) +  f^{-1}(D_1) + \ldots +  f^{-1}(D_{i-1}) +  f^{-1}(D_{i+1}) + \ldots +  f^{-1}(D_k).\]\end{small}
Since $f^{-1} (D_l)$ is a Cartier divisor for any $l \in \{1 , \ldots , k\}$, this implies $f^{-1} (D_i ) \subset f^{-1} (Y_i)$. So we obtain a unique $X$-morphism $\phi : \Bl _{\bigcap _{i \in I }( Y_i + D_1+ \ldots + D_{i-1} + D_{i+1} + \ldots + D_k)} ^{D_1+ \ldots + D_k} X \to \Bl _{\{Y_i\}_{i \in I}}^{\{D_i\}_{i \in I} }X  $. To check that $\phi$ is an isomorphism, it is enough to prove that there is an $X$-morphism \[\varphi :  \Bl _{\{Y_i\}_{i \in I}}^{\{D_i\}_{i \in I} }X \to \Bl _{\bigcap _{i \in I }( Y_i + D_1+ \ldots + D_{i-1} + D_{i+1} + \ldots + D_k)} ^{D_1+ \ldots + D_k} X .\] To build $\varphi$, we consider the map $f': \Bl _{\{Y_i\}_{i \in I}}^{\{D_i\}_{i \in I} }X \to X$ and check that $f'^{-1} (D_1+ \ldots +D_k) \subset f'^{-1} ( {\bigcap _{i \in I }( Y_i + D_1+ \ldots + D_{i-1} + D_{i+1} + \ldots + D_k)}) $. This is easy because $f'^{-1} (D_i ) \subset f'^{-1} (Y_i)$ for all $i$. Another method to prove Proposition \ref{formulamultimono} is to build $\phi $ or $\varphi$ and then reduce to the affine case and apply Proposition \ref{multifinimonoring}. \end{sloppypar}
\xpf

\subsection{Functoriality} \label{blow.up.base.functoriality.sec} We proceed with the notation from \S\ref{blow.up.definition.sec}.
Let $X'$ and $\{[Y'_i , D'_i]\} _{i\in I}  $ be another datum as in \S\ref{blow.up.definition.sec}. As usual, put $Z'_i= Y_i' \cap D'_i$.
A morphism $f:X'\to X$ such that, for all $i \in I$, its restriction to $D'_i$ (resp.~${Z'_i}$) factors through $D_i$ (resp.~$Z_i$), and such that $f^{-1} (D_i) = D_i'$, induces a unique morphism  $\Bl_{Y'}^{D'}X'\to \Bl_Y^DX$ such that the following diagram of $S$-spaces
\[
\begin{tikzpicture}[baseline=(current  bounding  box.center)]
\matrix(a)[matrix of math nodes, 
row sep=1.5em, column sep=2em, 
text height=1.5ex, text depth=0.45ex] 
{ 
\Bl_{Y'}^{D'}X'& \Bl_Y^DX \\ 
X'& X\\}; 
\path[->](a-1-1) edge node[above] {} (a-1-2);
\path[->](a-2-1) edge node[above] {} (a-2-2);
\path[->](a-1-1) edge node[right] {} (a-2-1);
\path[->](a-1-2) edge node[right] {} (a-2-2);
\end{tikzpicture}
\] 
commutes. 
This follows directly from Proposition \ref{blow.up.rep.prop}.

\rema
This is the spaces version of Fact \ref{functosemiano}.
\xrema 
\subsection{Base change} \label{blow.up.base.change.sec}
We proceed with the notation from \S\ref{blow.up.definition.sec}.
Let $X'\to X$ be a map of $S$-spaces, and denote by $Y'_i, Z'_i, D'_i \subset X'$ the preimage of $Y_i,Z_i,D_i \subset X$.
Then $D'_i \subset X'$ is locally principal for any $i$ so that the dilatation $\Bl_{Y'}^{D'}X'\to X'$ is well-defined. 
By \S\ref{blow.up.base.functoriality.sec} there is a canonical morphism of $X'$-spaces
\begin{equation}\label{blow.up.base.change.eq}
\Bl_{Y'}^{D'}X'\;\longto\; \Bl_Y^DX\x_{X}X'.
\end{equation}

\lemm\label{blow.up.base.change.lemm}
If $\Bl_Y^DX\x_{X}X'\to X'$ is an object of $\Spac_{X'}^{D'\text{-}\reg}$, then \eqref{blow.up.base.change.eq} is an isomorphism.
\xlemm
\pf
Our claim is étale local on $X$ and $X'$.
We reduce to the case where both $X=\Spec(B)$, $X'=\Spec(B')$ are affine, and $J_i=(b_i)$ is principal for all $i$.
We denote $Z'_i=\Spec(B'/L'_i)$ and $D_i'=\Spec(B'/J'_i)$. 
Then $J_i'=(b'_i)$ is principal as well where $b'_i$ is the image of $b_i$ under $B\to B'$.
We need to show that the map of $B'$-algebras
$
B'\otimes_BB\big[\frac{L}{b}\big]\;\longto\; B'\big[\frac{L'}{b'}\big]
$
is an isomorphism.
However, this map is surjective with kernel the $b'^{\bbN _I}$-torsion elements by Lemma \ref{tenseurkernel}.
As $b'_1, \ldots , b'_i, \ldots $ are non-zero-divisors in $B'\otimes_BB\big[\frac{L}{b}\big]$ by assumption, the lemma follows.
\xpf

\coro \label{blow.up.base.change.cor}
If the morphism $X'\to X$ is flat and satisfies a property $\calP$ which is
stable under base change, then $\Bl_{Y'}^{D'}X'\to\Bl_{Y}^{D}X$ is flat
and satisfies $\calP$.
\xcoro
\pf
Since flatness is stable under base change the projection $p\co\Bl_Y^DX\x_{X}X'\to \Bl_Y^DX$ is flat and has property $\calP$.
By Lemma \ref{blow.up.base.change.lemm}, it is enough to check that the closed subspace $\Bl_Y^DX\x_XD_i'$ defines an effective Cartier divisor on $\Bl_Y^DX\x_{X}X'$ for all $i$.
But this closed subscheme is the preimage of the effective Cartier divisor $\Bl_Y^DX\x_XD_i$ under the flat map $p$, and hence is an effective Cartier divisor as well by \cite[\href{https://stacks.math.columbia.edu/tag/083Z}{Tag 083Z}]{stacks-project}.
\xpf


\subsection{Relation to multi-centered affine projecting cone}
We proceed with the notation from \S\ref{blow.up.definition.sec} and assume that $\{D_i\} _{i \in I}$ belong to $Car(X)$. In this case, we can also realize $\Bl_Y^DX$ as a
closed subspace of the multi-centered affine projecting cone associated to $X,Z$ and $D$. 

\defi 
The {\em affine projecting cone  $\calO_X$-algebra } with multi-center $\{[Z_i=V(\calL_i) , D_i = V (\calJ _i)]\}_{i \in I} $ is
\[
\Cone_{\calL}^{\calJ} \calO_ X\defined\bigoplus_{\nu \in \mathbb{N} _I} \calL^\nu \otimes\calJ^{-\nu}.
\]
The {\em affine projecting cone }  of $X$ with multi-center $\{[Z_i , D_i]\}_{i \in I} $ is \[ \Cone _Z^D X \defined \Spec \big( \Cone_{\calL}^{\calJ} \calO_ X \big).\]
\xdefi
\prop \label{blow.up.closed.in.cone.lemm} The dilatation $\Bl_Z^DX$ is the closed
subspace of the affine projecting cone $\Cone_Z^DX$ defined by the equations $\{\varrho_i-1\}_{i\in I} $,
where  for all $i \in I$, $\varrho _i\in \Cone_{\calL}^{\calJ} \calO_ X$ is the image of $1 \in \calO _X$ under the map
\[\calO_X \cong\calJ _i\otimes\calJ _i ^{-1}\subset \calL _i\otimes\calJ _i^{-1} \subset \Cone_{\calL}^{\calJ} \calO_ X. \]
\xprop
\pf We may
work locally and the proposition follows from Proposition \ref{conicrhoblowup} and Remark \ref{rhomoinsunremark}.
\xpf
\rema
Let $X$ be a scheme and let $D$ be an effective divisor in $Div^+ (X) $ (cf. \cite{Du05}). Let $Z \subset D$ be a closed subscheme of finite type. Then the mono-centered dilatation $\Bl _Z ^D X$ was already defined in \cite[Definition 2.9]{Du05} using the conic point of view.
\xrema

\subsection{Relation to $\Proj$ of multi-graded algebras and multi-centered blowups} 

In \cite[Lemma 2.3]{MRR20}, it is proved that mono-centered dilatations are open subschemes of projective blowups. 
In the multi-centered case, this result can be generalized using multi-centered blowups and Fact \ref{expliaffinealgebra}. Multi-centered blowups are defined using Brenner-Schröer multi-graded Proj schemes. To relate Brenner-Schröer Proj schemes to multi-centered dilatations, it is convenient to use the concept of relevant families in a multi-graded ring as introduced in \cite{MR24}. We refer to \cite{MR24} for relations between multi-graded Proj and multi-centered dilatations.

\subsection{Dilatations of quasi-coherent $\calO_X$-modules}
We proceed with the notation from §\ref{blow.up.definition.sec}.
Let $\calF$ be a quasi-coherent $\calO_X$-module. Working locally as in Definition \ref{defimultidilaalgsp} and using dilatations of modules (cf. §\ref{sectionmodulealg}), we obtain a canonical quasi-coherent sheaf on $\Bl _Y^D X $, denoted $\Bl _Y ^D \calF$ or $\calF \Big[\Big\{ \frac{\calM_i}{\calJ_i}\Big\}_{i \in I} \Big]$. Note that in general $\Bl _Y ^D \calF \not \cong \Bl^* \calF$ where $\Bl: \Bl_Y^D X \to X$ is the dilatation map.
This construction enjoys the following universal property.
 Let $QCoh(X)_{\calF}^{D\text{-reg}}$ be the category whose objects are morphisms of quasi-coherent $\calO_X$-modules $F: \calF \to \calF '$ with source $\calF$ such that, locally, each $\calJ_i$ is defined by a non-zero-divisor of the module $\calF'$, then

\[ \Hom _{QCoh(X)_{\calF}^{D\text{-reg}}} (\calF\Big[\Big\{ \frac{\calM_i}{\calJ_i}\Big\}_{i \in I} \Big], \calF ' ) = \begin{cases}\{*\}, \; \text{if }\calL _i \calF ' = \calJ _i  \calF' \text{ for } i \in I;\\ \varnothing,\;\text{else.}\end{cases}\]

\section{Iterated multi-centered dilatations} \label{sec4}
We proceed with the notation from § \ref{blow.up.definition.sec}. Let $\nu , \theta \in \bbN ^I$ such that $\theta \leq \nu$, i.e. $\theta _i \leq \nu _i $ for all $i\in I$.

\prop \label{prop1iter}
There is a unique $X$-morphism 
\[ \varphi _{\nu , \theta}: \Bl _{Y}^{D^{\nu}} X \to \Bl _{Y}^{D^{\theta}} X  . \]
\xprop
\pf
 Let $\varphi _{\nu} : \Bl _{Y}^{D^{\nu}} X \to X$, $\varphi _{\theta} : \Bl _{Y}^{D^{\theta}} X \to X$  be the dilatation maps.  Let $i \in I$, if $\theta _i =0$, then $\theta _i D_i = X$ and so $\varphi _{\nu}^{-1} (\theta _i D_i )=  \Bl _{Y}^{D^{\nu}} X$ is a Cartier divisor in $ \Bl _{Y}^{D^{\nu}} X$. If $\theta _i > 0$, then $\nu _i >0$ and $\varphi _\nu ^{-1} ( \nu _i D_i) = \nu _i \varphi _\nu ^{-1} (D_i) $ is Cartier and so $ \varphi _\nu ^{-1} (D_i)$ is Cartier because $Car (\Bl _{Y}^{D^{\nu}} X) $ is a face of $Pri (\Bl _{Y}^{D^{\nu}} X)$ (cf. the discussion in Notation \ref{monoidclosed}). Consequently $\varphi _\nu ^{-1} ( \theta _i D_i) = \theta _i \varphi _\nu ^{-1} (D_i) $  is a Cartier divisor. So we proved that $\Bl _{Y}^{D^{\nu} } X$ belongs to $\Spac _X^{D^{\theta}\text{-reg}}$. We now use Proposition \ref{blow.up.rep.prop}. Let $i \in I$, we have $\theta _i D_i \subset \nu _i D_i$. So we have $\varphi_\nu ^{-1} ( \theta _i D_i ) \subset \varphi _\nu^{-1} ( \nu _i D_i ) \subset \varphi _{\nu}^{-1}( Y \cap \nu _i D_i)$.  So \[\varphi _{\nu }^{-1} ( \theta _i D_i) \subset  \big(  \varphi _{\nu}^{-1}( Y \cap \nu _i D_i) \big) \cap \big( \varphi _{\nu }^{-1} ( \theta _i D_i) \big) = \varphi _\nu ^{-1} ( Y \cap \nu _i D_i \cap \theta _i D_i ) = \varphi _\nu ^{-1}  (\theta _i D_i \cap Y ).\] Now we apply Proposition \ref{blow.up.rep.prop} and finish the proof.
\xpf

Assume now moreover that $\nu , \theta \in \bbN _I \subset \bbN ^I.$
We will prove that, under some assumptions, $\varphi _{\nu , \theta}$ is a dilatation morphism with explicit descriptions. We need the following observation.

\prop \label{liftclolem}  Assume that we have a commutative diagram of $S$-spaces \[
\begin{tikzcd}B \ar[rr, "f'"] \ar[dr,"f"]  & & C \ar[dl, "F"] \\  &X \end{tikzcd}.\]
 Assume that $F$ is affine and $f$ is a closed immersion. Then $f'$ is a closed immersion.
\xprop
\pf A closed immersion is affine, so by \cite[\href{https://stacks.math.columbia.edu/tag/08GB}{Tag 08GB}]{stacks-project}, $f'$ is affine. Using \cite[\href{https://stacks.math.columbia.edu/tag/03M4}{Tag 03M4}]{stacks-project}, we reduce to the case where $B,C$ and $X$ are affine (taking an étale covering of $X$ by affine schemes). Now the assertion is clear because closed immersions of affine schemes correspond to surjective morphisms at the level of rings.
\xpf

\coro \label{liftclo}  Assume that we have a commutative diagram of $S$-spaces \[
\begin{tikzcd}B \ar[rr, "f'"] \ar[dr,"f"]  & & \Bl_Y^D X \ar[dl, ] \\  &X \end{tikzcd}\] where the right-hand side morphism is the dilatation map.
 Assume that $f$ is a closed immersion. Then $f'$ is a closed immersion.
\xcoro
\pf
Clear by Proposition \ref{liftclolem}.
\xpf

We now assume that $Z_i \subset Y_i$ is a Cartier divisor inclusion for all $i \in I$.
 Let $\mathcal{D}_i $ be the canonical diagram of closed immersions \[\begin{tikzcd} Y_i \ar[r] \arrow[rd, "\square",phantom] &  \Bl _{Y_i}^{\nu_i D_i } X   \\ Z_i \ar[u] \ar[r] & \ar[u] D_i \end{tikzcd} \] obtained by Propositions \ref{blow.up.rep.prop} and \ref{liftclo}.

\lemm \label{iterationlemm1}
Assume $I= \{i\}$ and let $n_i \in \bbN$. We have an identification of $X$-spaces \[\Bl _{Y_i}^{(\nu _i +n_i )D_i} X= \Bl _{Y_i}^{n_i D_i } \Bl _{Y_i}^{\nu _i D_i}  X.\]
\xlemm
\pf 
Both sides belong to $Spaces_X ^{D_i\text{-reg}}$. 
Let $\Upsilon :  \Bl _{Y_i}^{n_i D_i } \Bl _{Y_i}^{\nu _i D_i}  X \to  \Bl _{Y_i}^{\nu _i D_i}  X$ and $h :\Bl _{Y_i}^{\nu _i D_i}  X \to X$ be the dilatation maps. 
Corollary \ref{coroooooo} implies that in $Clo ( \Bl _{Y_i}^{\nu _i D_i} X)$, we have $Y_i + h^{-1} (\nu _i D_i ) = h^{-1} (Y_i)$.
So we have 
\begin{align*}
(h \circ \Upsilon) ^{-1} \Big((\nu _i + n_i )D_i \Big)&= \Upsilon ^{-1} \Big( h^{-1} (\nu _i D_i )\Big)+ \Upsilon ^{-1} \Big(h^{-1} (n_i D_i) \Big)\\ 
&= \Upsilon ^{-1} \Big( h^{-1} (\nu _i D_i )\Big)+ \Upsilon ^{-1} \Big(h^{-1} (n_i D_i) \cap Y_i  \Big)\\
&= \Upsilon ^{-1} \Big( h^{-1} (\nu _i D_i )+  (h^{-1} (n_i D_i) \cap Y_i ) \Big)\\
&= \Upsilon ^{-1} \Big( ( h^{-1} (\nu _i D_i )+  h^{-1} (n_i D_i) )  \cap ( h^{-1} (\nu _i D_i )+ Y_i ) \Big)\\
& =  \Upsilon ^{-1} \Big( ( h^{-1} ((\nu _i+n_i) D_i )  \cap  h^{-1} ( Y_i ) \Big)\\
 & = (h \circ \Upsilon) ^{-1} \Big( \big((\nu _i + n_i )D_i\big) \cap Y_i \Big).
\end{align*}
Consequently, by the universal property of dilatations, we have a unique morphism of $X$-spaces $\Bl _{Y_i}^{n_i D_i } \Bl _{Y_i}^{\nu _i D_i}  X \to \Bl _{Y_i}^{(\nu _i +n_i )D_i} X$.
 To prove that it is an isomorphism, we reduce to the affine case and apply Corollary \ref{coroooooo}.
\xpf

Let $f_i$ be the canonical morphism (e.g. cf. \ref{uniquemorphism}, \ref{multietape} or \ref{prop1iter})
\[ \Bl _{Y}^{D^{\nu}} X \to \Bl _{Y_i}^{{\nu_i}D_i } X  .\] 
 We denote by $Y_i \times _{\Bl _{Y_i}^{{\nu_i}D_i } X}\Bl _{Y}^{D^{\nu}}X$ the fiber product obtained via the arrows given by $f_i$ and $\mathcal{D}_i$. We use similarly the notation $D_i \times _{\Bl _{Y_i}^{{\nu_i}D_i } X}\Bl _{Y}^{D^{\nu}}X$.

\lemm \label{nueiei} Let $i \in I$ and let $n _i \in \bbN$. Let $\gamma _i \in \bbN _I$ be $(0,\ldots ,0, n_i, 0 , \ldots )$ where $n_i$ is in place $i$. We have an identification
 \[\Bl _{Y}^{D^{\nu+ \gamma_i}}X = \Bl ^{n_iD_i \times _{\Bl _{Y_i}^{{\nu_i}D_i } X}\Bl _{Y}^{D^{\nu}}X} _{Y_i \times _{\Bl _{Y_i}^{{\nu_i}D_i } X}\Bl _{Y}^{D^{\nu}}X} \Bl _{Y}^{D^{\nu}}X  .\] In particular we have a canonical dilatation morphism \[\varphi _{\nu +\gamma_i , \nu}:\Bl _{Y}^{D^{\nu+ \gamma_i}}X \to \Bl _{Y}^{D^{\nu}}X.\]
\xlemm
\pf  We have
 \begin{align*}& \Bl ^{{n_iD_i \times _{\Bl _{Y_i}^{{\nu_i}D_i } X}\Bl _{Y}^{D^\nu} X}} _{{Y \times _{\Bl _{Y_i}^{{\nu_i}D_i } X}\Bl _{Y}^{D^\nu} X}} \Bl _{Y}^{{D^\nu}}X \\
\text{by Proposition \ref{multietape}} =~&  \Bl ^{{n_iD_i \times _{\Bl _{Y_i}^{{\nu_i}D_i } X}\Bl _{Y}^{D^\nu} X}} _{{Y \times _{\Bl _{Y_i}^{{\nu_i}D_i } X}\Bl _{Y}^{D^\nu} X}} \Bl _{\{Y_j \times _X \Bl _{Y_i}^{{\nu _i D_i}}X\}_{j \in I \setminus \{i\}}}^{{\{\nu_j D_j \times _X \Bl _{Y_i}^{{\nu _i D_i}}X\}_{j \in I \setminus \{i\}}}} \Bl _{Y_i}^{{\nu _i D_i}}X  \\
\text{by Proposition \ref{multietape}} 
=~ & \Bl _{Y_i,  \{Y_j \times _X \Bl _{Y_i} ^{{{\nu _i}D_i} } X \}_{j \in I \setminus \{i\} } }^{n_iD_i ,  \{\nu_j D_j  \times _X \Bl _{Y_i} ^{{\nu _i D_i^{}} } X \}_{j \in I \setminus \{i\} }} \Bl _{Y_i}^{{{\nu _i}D_i}}X \\
\text{by Prop. \ref{multietape} and Lem. \ref{iterationlemm1}} =~ & \Bl _{\{Y_j \times _X \Bl _{Y_i} ^{{{(\nu _i +n_i)}D_i} } X \}_{j \in I \setminus \{i\} }} ^{\{\nu _jD_j \times _X \Bl _{Y_i}^{(\nu_i +n_i)D_i }X\}_{j \in I \setminus \{i\}}} \Bl _{Y_i}^{({\nu_i +{n _i}})D_i }X \\
 \text{by Proposition \ref{multietape}}=~& \Bl _{Y}^{D^{\nu+ \gamma_i}}X.
 \end{align*}
\xpf
\prop \label{nutheta} Recall that $\theta \leq \nu $. Put $\gamma = \nu - \theta$. Put $K = \{ i \in I | \gamma _i >0 \}.$ We have an identification
 \[\Bl _{Y}^{D^{\nu}} X = \Bl ^{\{{{\gamma _i}D_i \times _{\Bl _{Y_i}^{\theta _i D_i } X}\Bl _{Y}^{D^{\theta}} X\}_{i \in K}}} _{\{Y_i   \times _{\Bl _{Y_i}^{\theta _i D_i } X}\Bl _{Y}^{D^{\theta}} X \}_{i \in K}} \Bl _{Y}^{D^{\theta}} X .\]
  In particular the unique $X$-morphism 
\[ \varphi _{\nu , \theta}: \Bl _{Y}^{D^{\nu}} X \to \Bl _{Y}^{D^{\theta}} X   \] of Proposition \ref{prop1iter} is a dilatation map.
\xprop

\pf We prove the first assertion by induction on $\# \{i \in I |\nu _i >0 \}$. If $k=1$ the assertion follows from Lemma \ref{iterationlemm1}. The passage from $k-1$ to $k$ follows from Lemma \ref{nueiei} and Proposition \ref{multietape}. 
\xpf 

It is now natural to introduce the following terminology.

\defi 
 For any $\nu \in \bbN _I, $ let us consider
\[ \Bl _{Y}^{D^{\nu}} X= \Bl \big\{_{Y_i}^{\nu _i D_i}\big\}_{i \in I} X\] and call it the $\nu$-th iterated dilatation of $X$ with multi-center $\{[Y_i,D_i]\}_{i\in I}$.
\xdefi

\section{Multi-centered dilatations along multiples of a single divisor} \label{sectionsinglediv}

Let $X$ be an $S$-space. We fix a locally principal closed subscheme $D\subset X$. Let $Y_0 , Y_1 ,\ldots , Y_i, \ldots, Y_k$ be closed $S$-subspaces of $X$ such that $D \cap Y_i \subset Y_i$ is a Cartier divisor for all $i$.  We assume moreover that $Y_0 \subset Y_i $ for $i \in \{1 , \ldots , k \}$. Let $s_0 , s_1 , \ldots , s_k \in \bbN$ be integers.
We claim that we have a canonical closed immersion $Y_0 \to \Bl _{Y_0 ,~~ \ldots~ , Y_k }^{s_0D,  \ldots , {s_k}D} X$. This follows from Propositions \ref{blow.up.rep.prop} and \ref{liftclo} observing that the map $Y_0 \to X$ restricted to $s_i D$ factors through $Y_i \cap s_i D$ for any $i$. We now use the notation $\Bl _{Y_0 ,\ldots , Y_k }^{s_0,  \ldots , {s_k}} X$ to denote $\Bl _{Y_0 ,~~ \ldots~ , Y_k }^{s_0D,  \ldots , {s_k}D} X$. The following fact is a direct generalization of the first assertion of Corollary \ref{coroooooo} to the present situation.

\fact \label{idealL}Assume that $X= \Spec (A)$ is affine. Assume that $D= \Spec (A/(a))$ and $Y_i = \Spec (A/M_i)$ for $i \in I$. Then the ideal $Q$ of $A':=A[\frac{M_0}{a^{s_0}}, \frac{M_1}{a^{s_1}}\cdots \frac{M_k}{a^{s_k}}]$ corresponding to the canonical closed immersion $Y_0  \to \Bl _{Y_0 , Y_1, \ldots , Y_k }^{{s_0}, {s_1} , \ldots , {s_k}} X$ is the ideal $\langle \frac{M_0}{a^{s_0}}, \frac{M_1}{a^{s_1}} , \ldots , \frac{M_k}{a^{s_k}} \rangle$ of $A'$ generated by $\frac{M_0}{a^{s_0}}, \frac{M_1}{a^{s_1}} , \ldots , \frac{M_k}{a^{s_k}}.$
\xfact 

\pf There is no difficulty to adapt the proof of \ref{coroooooo}. We provide details for the convenience of the reader.
Using \ref{idealsumclosed}, it is clear that  $\langle \frac{M_0}{a^{s_0}}, \frac{M_1}{a^{s_1}} , \ldots , \frac{M_k}{a^{s_k}} \rangle \subset Q$, so it is enough to prove that $Q \subset \langle \frac{M_0}{a^{s_0}}, \frac{M_1}{a^{s_1}} , \ldots , \frac{M_k}{a^{s_k}} \rangle$.
So let $\nu \in  \bbN ^k $, we have to prove that $ \frac{L^\nu \cap M_0}{(a^s)^\nu}$ is included in $\langle \frac{M_0}{a^{s_0}}, \frac{M_1}{a^{s_1}} , \ldots , \frac{M_k}{a^{s_k}} \rangle$. An element $x \in L^\nu $ can be written as a sum
$x= \sum _{\nu =\beta + \alpha} m_{\beta} (a^s)^\alpha $ with $m_\beta \in M^\beta$  (note that, if $x $ belongs to $L^\nu  \cap M_0$, then $m_{(0,\ldots, 0)} (a^s)^\nu$ also belongs to $M_0$). Now we assume that $x $ belongs to $L^\nu  \cap M_0$, it is clear that for $\beta \ne (0,\ldots,0)$ the element $\frac{  m_{\beta} (a^s)^\alpha}{(a^s)^\nu} = \frac{ m_\beta }{(a^s)^\beta}$ belongs to $\langle \frac{M_0}{a^{s_0}}, \frac{M_1}{a^{s_1}} , \ldots , \frac{M_k}{a^{s_k}} \rangle$. For $\beta=(0,\ldots,0)$, using that $a_i$ is a non-zero-divisor in $A/M_i$ for all $i$ and that $m_{(0,\ldots, 0)} (a^s)^\nu$ belongs to $M_0$, we get that $m_{(0,\ldots ,0)} $ belongs to $M_0$ and it is now clear that $\frac{m_{(0,\ldots,0)}(a^s)^\nu }{(a^s)^\nu}$ belongs to $\langle \frac{M_0}{a^{s_0}}, \frac{M_1}{a^{s_1}} , \ldots , \frac{M_k}{a^{s_k}} \rangle$.
So $x$ belongs to $\langle \frac{M_0}{a^{s_0}}, \frac{M_1}{a^{s_1}} , \ldots , \frac{M_k}{a^{s_k}} \rangle$. 
\xpf 

\prop \label{clef} Let $0 \leq t \leq s_0$ be an integer. We have a canonical identification of $X$-spaces
\[ \Bl _{Y_0^{}}^{{t}} \Bl _{Y_0 , Y_1, \ldots  , Y_k }^{{s_0}, {s_1} , \ldots , {s_k}} X = \Bl _{Y_0 , ~~Y_1,^{~} \ldots ~,~ Y_k^{~~} }^{{s_0+t}, {s_1+t} , \ldots , {s_k+t}} X .\]
\xprop 

\pf 
Both sides belong to $Spaces_X ^{D\text{-reg}}$. 
Let $\Upsilon : \Bl _{Y_0^{}}^{{t}} \Bl _{Y_0 , Y_1, \ldots  , Y_k }^{{s_0}, {s_1} , \ldots , {s_k}} X \to \Bl _{Y_0 , Y_1, \ldots  , Y_k }^{{s_0}, {s_1} , \ldots , {s_k}} X$ and $h :\Bl _{Y_0 , Y_1, \ldots  , Y_k }^{{s_0}, {s_1} , \ldots , {s_k}} X \to X$ be the dilatation maps. 
Fact \ref{idealL} implies that, in $Clo (\Bl _{Y_0 , Y_1, \ldots  , Y_k }^{{s_0}, {s_1} , \ldots , {s_k}} X)$, we have $Y_0 + h^{-1} (s_i D ) \subset  h^{-1} (Y_i)$ for all $0 \leq i \leq s $.
So we have 
\begin{align*}
(h \circ \Upsilon) ^{-1} \Big((t+ s_i )D \Big)&= \Upsilon ^{-1} \Big( h^{-1} (s_i D )\Big)+ \Upsilon ^{-1} \Big(h^{-1} (t D) \Big)\\ 
&= \Upsilon ^{-1} \Big( h^{-1} (s_i D )\Big)+ \Upsilon ^{-1} \Big(h^{-1} (t D) \cap Y_0  \Big)\\
&= \Upsilon ^{-1} \Big( h^{-1} (s_i D )+  (h^{-1} (t D) \cap Y_0 ) \Big)\\
&= \Upsilon ^{-1} \Big( ( h^{-1} (s_i D )+  h^{-1} (t D) )  \cap ( h^{-1} (s_i D )+ Y_0 ) \Big)\\
& \subset   \Upsilon ^{-1} \Big( ( h^{-1} ((t+s_i) D )  \cap  h^{-1} ( Y_i ) \Big)\\
 & = (h \circ \Upsilon) ^{-1} \Big( \big((t+s_i )D\big) \cap Y_i \Big).
\end{align*}
Consequently, by the universal property of dilatations, we have a unique morphism of $X$-spaces $\Bl _{Y_0^{}}^{{t}} \Bl _{Y_0 , Y_1, \ldots  , Y_k }^{{s_0}, {s_1} , \ldots , {s_k}} X \to \Bl _{Y_0 , ~~Y_1,^{~} \ldots ~,~ Y_k^{~~} }^{{s_0+t}, {s_1+t} , \ldots , {s_k+t}} X$.
 To prove that it is an isomorphism, we reduce to the affine case and use the notation of Fact \ref{idealL}. It is enough to show that we have a canonical identification of rings
\begin{center}$ \big(A[\frac{M_0}{a^{s_0}}, \frac{M_1}{a^{s_1}}\cdots \frac{M_k}{a^{s_k}}]\big)[\frac{Q}{a^t}] =A[\frac{M_0}{a^{s_0 +t}}, \frac{M_1}{a^{s_1 +t}}\cdots \frac{M_k}{a^{s_k +t}}]$.\end{center} For this, as in the  proof of the second assertion of Corollary \ref{coroooooo}, it is enough to apply Proposition  \ref{basechangeetapering} and Corollary \ref{maxsemi}.
\xpf

\fact \label{examployno} Let $A$ be a ring. Let $P = A[T , X_1 , \ldots, X_n]$ be the polynomial algebra in $n+1$ variables. Let $d_1 , \ldots, d_n \in \bbN$. Then there is a canonical identification of $A$-algebras
\begin{center} $ P[x_1 , \ldots , x_n ] / (X_1-T^{d_1} x_1 , \ldots ,X_n-T^{d_n} x_n ) = P[ \frac{(X_1)}{T^{d_1}}, \ldots , \frac{(X_n)}{T^{d_n}}] $.\end{center}
\xfact 
\pf
The map given by $x_i \mapsto \frac{X_i}{T^{d_i}}$ is well-defined and surjective.  The source and target of our map are $T$-torsion free and the map is an isomorphism after inverting $T$ by Corollary \ref{corofds}.
\xpf 

\prop \label{surjectiongipoli}Let $A$ be a ring. Let $a \in A$. Let $g_1 , \ldots , g_n $ be elements in $A$ where $n \in \bbN$. Let $d_1 , \ldots, d_n \in \bbN$. There is a surjection 
\begin{center}$ A[x_1 , \ldots , x_n] / (g_1-a^{d_1} x_1, \ldots , g_n-a^{d_n} x_n ) \to A[\frac{(g_1)}{a^{d_1}}, \ldots , \frac{(g_n)}{a^{d_n}}] $\end{center} whose kernel is the $a$-power torsion in the source.
\xprop 

\pf Consider the map $ P = A [T , X_1 , \ldots , X_n] \to A$ sending $T $ to $a$ and $X_i$ to $g_i$ for $1 \leq i \leq n$. By Fact \ref{examployno}, we have $P[x_1 , \ldots , x_n ] / (X_1-T^{d_1} x_1 , \ldots ,X_n-T^{d_n} x_n ) = P[ \frac{(X_1)}{T^{d_1}}, \ldots , \frac{(X_n)}{T^{d_n}}] $. Now we use Proposition \ref{tenseurkernel} to finish the proof.
\xpf 

\prop \label{regupoly}Let $A$ be a ring. Let $a , g_1 , \ldots , g_n$ be a $H_1$-regular sequence in $A$ (cf. \cite[\href{https://stacks.math.columbia.edu/tag/062E}{Tag 062E}]{stacks-project} for $H_1$-regularity). Let $d_1 , \ldots , d_n$ be positive integers. Then the surjection of Proposition \ref{surjectiongipoli} is an isomorphism. In particular, the dilatation algebra identifies with a quotient of a polynomial algebra as follows 
\begin{center}
$A[\frac{(g_1)}{a^{d_1}}, \ldots , \frac{(g_n)}{a^{d_n}}]=A[x_1 , \ldots , x_n] / (g_1 -a^{d_1} x_1, \ldots , g_n -a^{d_n} x_n )$. 
\end{center}

\xprop 
\pf We can assume that $d_i >0$  for all $i$ by Corollary \ref{1serarien}.
By Proposition \ref{surjectiongipoli}, it is enough to show that the right-hand side is $a$-torsion free. We adapt the proof of  \cite[\href{https://stacks.math.columbia.edu/tag/0BIQ}{Tag 0BIQ}]{stacks-project}. We claim that the sequence $( a, g_1-a^{d_1}x_1, \ldots , g_n -a^{d_n} x_n)$ is $H_1$-regular in $A[x_1, \ldots , x_n]$. Namely, the map \[( a, g_1-a^{d_1}x_1, \ldots , g_n -a^{d_n} x_n):A[x_1 , \ldots ,x_n]^{\oplus (1+n)} \to A[x_1 , \ldots ,x_n] \] used to define the Koszul complex on $a, g_1-a^{d_1}x_1, \ldots , g_n -a^{d_n} x_n $ is isomorphic to the map 
\[(a,g_1,\ldots,g_n): A[x_1 , \ldots ,x_n]^{\oplus (1+n)} \to A[x_1 , \ldots ,x_n] \] used to define the Koszul complex on $a,g_1,\ldots,g_n$ via the isomorphism $\Theta$
\[ A[x_1 , \ldots ,x_n]^{\oplus (1+n)} \to A[x_1 , \ldots ,x_n]^{\oplus (1+n)} \] sending $(P_0 , P_1, \ldots , P_n)$ to 
\[ (P_0 - \sum _{i=1}^{n} a^{d_i-1} x_i P_i, P_1 , P_2 , \ldots , P_n);\] this follows from the identity
\[ aP_0+ \sum _{i=1}^n (g_i-a^{d_i} x_i ) P_i = a\big( P_0 - \sum _{i=1}^n a^{d_i -1} x_i P_i \big) + \sum _{i=1}^n g_i P_i .\] By \cite[\href{https://stacks.math.columbia.edu/tag/0624}{Tag 0624}]{stacks-project} these Koszul complexes are isomorphic. By \cite[\href{https://stacks.math.columbia.edu/tag/0629}{Tag 0629}]{stacks-project} the Koszul complex $K$ on $( a, g_1-a^{d_1}x_1, \ldots , g_n -a^{d_n} x_n)$ is the cone on $a: L \to L$ where $L$ is the Koszul complex on $(g_1-a^{d_1}x_1, \ldots , g_n -a^{d_n} x_n)$, since $H_1 (K)=0$, we conclude that $a: H_0(L) \to H_0(L) $ is injective, so the right-hand side is $a$-torsion free. 
\xpf 
\section{Some flatness and smoothness results}
\label{sec6}
Let $S$ be a scheme and let $C \subset S$ be a Cartier divisor in $S$.
 Let $X$ be a scheme over $S$. Let $D$ be the closed subscheme of $X$ given by $X \times _S C$. Let $X_{j} \subset D$ be closed subchemes for $1 \leq j \leq d$.
 We are now making the following assumption. We assume that locally over $S,X$ the following conditions are satisfied \begin{enumerate}
 \item  $S= \Spec (R)$,  $C = \Spec (R/a)$ and $X= \Spec (A)$,
 \item there exists a sequence $g_1 , \ldots, g_n \in A$  such that $a, g_1 , \ldots, g_n$ is a $H_1$-regular sequence in $A$, 
 \item there exists a sequence  $1 \leq i_1 <   i_2 < \ldots <i_j<  \ldots < i_d =n$ such that 
\begin{align*}
X_1&= \Spec \big( A/ (g_1 , \ldots , g_{i_1})\big)\\
X_2 &= \Spec \big( A / (\overset{\pm}{g_1} , \ldots , \overset{\pm}{g_{i_1}}, g_{i_1+1}, \ldots , g_{i_2}  )\big)\\
\vdots & ~~~~~~~~ \vdots \\
X_j &= \Spec \big ( A / (\overset{\pm}{g_1} , \ldots , \overset{\pm}{g_{i_{j-1}}}, g_{i_{j-1}+1}, \ldots , g_{i_j}) \big)\\
\vdots & ~~~~~~~~ \vdots \\
X_d& = \Spec ( A / (\overset{\pm}{g_1} , \ldots , \overset{\pm}{g_{i_{d-1}}}, g_{i_{d-1}+1}, \ldots , g_{n}) \big)
\end{align*}
where $\pm$ over a symbol means that this symbol possibly appears but not necessarily. We put $U_j = \Spec ( A/ (g_1 ,\ldots, g_{i_j}))$ for $1 \leq i \leq j $.
\end{enumerate}   Let $m_1 \geq  \ldots \geq  m_j \geq  \ldots \geq m_d  \geq 0 \in \bbN $ be integers. 
\prop \label{flatandsmooth} The following assertions hold. 
\begin{enumerate}
\item If $X/S$ is flat and if moreover one of the following holds:
\begin{enumerate}
\item $X_j \to S$ is flat and $S,X$ are locally noetherian for $1 \leq j \leq d$,
\item $  X_j \to S $ is flat and $ X_j \to S$ is locally of finite presentation for all $j$,
\item the local rings of $S$ are valuation rings,
\end{enumerate}
then $\Bl _{X_1, \ldots , X_d}^{m_1 , \ldots , m_d} X  \to S$ is flat.
\item  If $X \to S$ is smooth and (with the local notation of the assumption) $U_j \times _S \Spec ( R/ (a^{m_j})) \to \Spec ( R/ (a^{m_j}))$ is smooth for all $ 1 \leq j \leq d$, then $\Bl _{X_1, \ldots , X_d}^{m_1 , \ldots , m_d} X  \to S$ is smooth.
\end{enumerate}
\xprop 
\pf This is local on $S,X$. We use notations used to state the assumption before the statement.   Corollary \ref{maxsemi} implies that we can and do assume that $X_j = U_j $ for all $1 \leq j \leq d $.\begin{enumerate} \item  We prove the assertions by induction on $d$. If $d=1$ this follows from \cite[Proposition 2.16]{MRR20}. By Proposition \ref{multietape} we have 
$\Bl _{X_1, \ldots , X_d}^{m_1 , \ldots , m_d} X  = \Bl _{X_d'}^{m_d} \Bl _{X_1, \ldots , X_{d-1}}^{m_1 , \ldots , m_{d-1}} X$ where $X'_d=\Bl _{X_1, \ldots , X_{d-1}}^{m_1 , \ldots , m_{d-1}} X \times _{X} X_d$.   
 Proposition \ref{regupoly} implies that $\Bl _{X_1, \ldots , X_{d-1}}^{m_1 , \ldots , m_{d-1}} X $ identifies with the spectrum of $P/ \frakI$ where $ P= A [ x_1 , \ldots , x_{i_{d-1}}]$ and  
\begin{center}
$  \frakI = ( g_1 -a^{m_1} x_1, \ldots , g_{i_1} - a^{m_1} x_{i_1} , g_{i_1+1} - a^{m_2} x_{i_1+1} , \ldots ,  g_{i_2} - a^{m_2} x_{i_2}, \ldots , g_{i_{d-1}} - a^{m_{d-1}} x_{i_{d-1}} ) . $
\end{center}
We claim that the sequence given by $g_{i_{d-1}+1} , \ldots , g_{i_d}$ is $H_1$-regular in $ P/(\frakI +(a^{m_d}))$. Let us prove the claim. 
Since $m_d \leq m_j $ for all $1 \leq j \leq d-1$, the ideal $\frakI + (a^{m_d})$ of $A$ is equal to 
$(a^{m_d},g_1 , \ldots , g _{i_{d-1}} ) 
$. So $P / (\frakI + (a^{m_d})) $ identifies with 
$\big( A/(a^{m_d}, g_1 , \ldots , g _{i_{d-1}} ) \big) [x_1 , \ldots , x_{i _{d-1}} ].$\begin{sloppypar}
Now since $(a, g_1, \ldots , g _n)$ is $H_1$-regular in $A$, we know that $(a^{m_d}, g_1 , \ldots , g_n) $ is $H_1$-regular in $A$ by  \cite[\href{https://stacks.math.columbia.edu/tag/062G}{Tag 062G}]{stacks-project}. So $g_{i_{d-1}+1} , \ldots , g_{i_d}$ is a $H_1$-regular sequence in $A / (a^{m_d}, g_1 , \ldots , g_{i_{d-1}} )$ by \cite[\href{https://stacks.math.columbia.edu/tag/068L}{Tag 068L}]{stacks-project}. This implies that  $g_{i_{d-1}+1} , \ldots , g_{i_d}$ is a $H_1$-regular sequence in $\big( A/(a^{m_d}, g_1 , \ldots , g _{i_{d-1}} ) \big) [x_1 , \ldots , x_{i _{d-1}} ]$. This finishes to prove our claim. We now apply \cite[Proposition 2.16]{MRR20} to deduce the assertions for $d$.\end{sloppypar}
\item   We prove the assertion by induction on $d$. If $d=1$, this is \cite[Proposition 2.16]{MRR20}. We now assume that the assertion is true for $d-1$. Using \cite[Proposition 2.16]{MRR20} and Proposition \ref{multietape}, it is enough to show that $Z_d ':= (X_d \times _X \Bl _{X_1, \ldots , X_d}^{m_1 , \ldots , m_d} X ) \times _S ( \Spec (R / (a^{m_d}) )$ is smooth over $\Spec (R / (a^{m_d})$.   We have \[Z_d'= X_d \times _S \Spec ( R / (a^{m_d})) \times _{\Spec ( R / (a^{m_d}))} \Bl _{X_1, \ldots , X_d}^{m_1 , \ldots , m_d} X  \times _S  \Spec (R / (a^{m_d}) .\]\begin{sloppypar}
We computed $\Bl _{X_1, \ldots , X_d}^{m_1 , \ldots , m_d} X  \times _S  \Spec (R / (a^{m_d})$ in (i) and proved that it is the spectrum of $A/(a^{m_d},g_1 , \ldots, g_{i_{d-1}})[x_1 , \ldots , x_{i_{d-1}}]$. Now $X_d \times _S \Spec ( R / (a^{m_d}))$ is the spectrum of $A/ (a^{d_m}, g_1 , \ldots , g_n)$. Consequently $Z_d'$ is the spectrum of $(A/(a^{d_m}, g_1 , \ldots, g_n ) )[x_1 , \ldots , x_{i_{d-1}}]$. So by assumption $A/(a^{d_m}, g_1 , \ldots, g_n ) $ is smooth over $R/(a^{d_m})$. Moreover any polynomial algebra is smooth over its base ring. A composition of smooth morphisms is smooth. So the composition 
\[ R/(a^{d_m}) \to A/(a^{d_m}, g_1 , \ldots, g_n ) \to (A/(a^{d_m}, g_1 , \ldots, g_n ) )[x_1 , \ldots , x_{i_{d-1}}] \] is smooth. This finishes the proof.
\end{sloppypar}
\end{enumerate}
\xpf 
\section{Dilatations of monoid, group and Lie algebra schemes} \label{néron}

 We study dilatations of schemes endowed with a structure (cf. \cite[Exp. I §2.2]{SGA3}) in this section. We focus on monoid, group and Lie algebra structures.
Let $S$ be a scheme, and let $G\to S$ be a monoid (resp. group, resp. Lie algebra) scheme.
Let $C=\{C_i\} _{i\in I} \subset S$ be locally principal closed subschemes.  Put $D_i= G|_{C_i} = G \times _S C_i$ and $D= \{D_i\} _{i\in I}$.
Let $H_i\subset D_i $ be a closed submonoid (resp. subgroup, resp. Lie subalgebra) scheme over $C_i$ for all $i \in I$ and let $H= \{H_i\}_{i \in I}$. 
Let $\calG:=\Bl_H^{D}G\to G$ be the associated dilatation.
The structure morphism $\calG\to S$
defines an object in $Sch_S^{C\text{-}\reg}$.

\fact\label{Neron.blow.lemm} 
Let $\calG\to S$ be the above dilatations.
\begin{enumerate}
\item[(1)] The scheme $\calG\to S$ represents the contravariant functor
$Sch_S^{C\text{-}\reg}\to  Set$ given for $T\to S$ by the set of all $S$-morphisms $T\to G$ such that the induced morphism $T|_{C_i}\to G|_{C_i}$ factors through $H_i\subset G|_{C_i}$ for all $i \in I$.

\item[(2)] Let $T \to S$ be an object in $Sch_S^{C\text{-}\reg}$, then as subsets of $G(T)$ \[ \calG (T) = \bigcap_{i \in I} \big( \Bl_{H_i}^{D_i}G \big)(T).\]

\item[(3)] The map $\calG\to G$ is affine. 
Its restriction over $S_i$ factors as $\calG_i\to H_i \subset  D_i$ for all $i \in I.$
\end{enumerate}
\xfact
\pf
Part (1) is a reformulation of Proposition \ref{blow.up.rep.prop}. Assertion (2) is immediate using (1). Finally (3) is immediate from Proposition \ref{blow.up.Cartier.lemm}.
\xpf

We now assume that the category $Sch_S^{C \text{-}\reg}$ has products (cf. e.g. Fact \ref{product} for some conditions ensuring this hypothesis).

\prop \label{kindof}  If the dilatation $\calG\to S$ is flat (e.g. using Proposition \ref{flatandsmooth}), then it is equipped with the structure of a monoid (resp. group, resp. Lie algebra) scheme over $S$ such that $\calG\to G$ is a morphism of $S$-monoid (resp. $S$-group, resp. $S$-Lie algebra) schemes.
\xprop 
\pf
By virtue of Fact \ref{Neron.blow.lemm} the (forgetful) map
$\calG\to G$ defines a submonoid (resp. subgroup, resp. Lie subalgebra) functor when restricted to the
category $Sch_S^{C\text{-}\reg}$. As $\calG\to S$ is an object in
$Sch_S^{C \text{-}\reg}$, it is a monoid (resp. group, resp. Lie algebra) object in this category. Now if $X \to S $ and $Y \to S$ are two flat morphisms in $Sch _S^{C\text{-}\reg}$, then the product of $X $ and $Y$ in the category $Sch _S^{C\text{-}\reg}$ exists and is equal to the product of $X$ and $Y$ in the full category of $S$-schemes. 
So $\calG\to S$ is a monoid (resp. group, resp. Lie algebra) object in the full category of $S$-schemes.
\xpf 
Dilatations of group schemes are often called Néron blowups. We note that dilatations preserve similarly structures defined using products and commutative diagrams (cf. \cite[Exp. I §2.2]{SGA3}).

\prop \label{LieLie} Assume that $C_i\subset S$ is a Cartier divisor for all $i$. Assume that $G \to S$ and $H_i \to C_i$ are flat group schemes. Assume that $\Bl _H ^D G \to S$ is flat (and so a group scheme by \ref{kindof}). 
Let $\bbL ie (G)/S$ (resp. $\bbL ie ( H_i )/ C_i$, resp. $\bbL ie (\Bl_H^{D}G )/S$)   be the Lie algebra scheme of $G/S$ (resp. $H_i/C_i$, resp. $\Bl_H^{D}G/S$) (cf. \cite[Exp. II Scholie 4.11.3]{SGA3}). Assume that $\Bl \big\{{}^{\bbL ie (G) \times _S C_i}_{\bbL ie ( H_i )} \big\}_{i \in I} \bbL ie (G) \to S$ is flat (and so a Lie algebra scheme by \ref{kindof}).
Then we have a canonical isomorphism of $S$-Lie algebra schemes:
\[ \bbL ie (\Bl_H^{D}G ) \cong \Bl \big\{{}^{\bbL ie (G) \times _S C_i}_{\bbL ie ( H_i )} \big\}_{i \in I} \bbL ie (G).\] 
\xprop 
\pf We have a morphism of $S$-group schemes $\Bl _{H}^D G \to G $, it induces a morphism of $S$-Lie algebra schemes $\bbL ie (\Bl _{H}^D G ) \to \bbL ie ( G) $. Using the universal property of dilatations of Lie algebras, we obtain a canonical morphism of $S$-Lie algebra schemes 
$\bbL ie (\Bl_H^{D}G ) \to \Bl \big\{{}^{\bbL ie (G) \times _S C_i}_{\bbL ie ( H_i )} \big\}_{i \in I} \bbL ie (G)$. We now show that it is an isomorphism.
Using flatness and the assumptions on divisors $C_i$, we see that both sides belong to $Sch _S ^{C\text{-reg}}$. It is enough (to finish the proof) to evaluate both sides on a test scheme $T \to S$ in $Sch _S ^{C\text{-reg}}$ and obtain an identification of sets. Recall that for any scheme $U$, $I_U = \Spec  (\bbZ[X]/(X^2) ) \times _{\Spec (\bbZ) } U$ denotes the scheme of dual numbers over $U$. The scheme $I_U$ is obviously flat over $U$. In particular $I_T \to S$ belongs to $Sch _S ^{C\text{-reg}}$. We have a canonical morphism $T \to I_T$ induced by  $\bbZ [X]/X^2 \to \bbZ , X \mapsto 0$.
Using \cite[Exp. II Sch. 4.11.3, Cor. 3.9.0.2.]{SGA3} we get
\begin{align*}(  \bbL ie& (\Bl_H^{D}G ) )( T)  = \ker \big( (\Bl_H^{D}G )(I_T) \to (\Bl_H^{D}G )(T) \big) \\ 
&= \big\{ I_T \to  \Bl_H^{D}G  \in  \Hom _S ( I_T , \Bl_H^{D}G ) \vert ~ T \to I_T \to  \Bl_H^{D}G  \text{ is the unit } \big\} \\
& = \big\{ I_T \to G \in \Hom _S ( I_T , G ) \vert ~~ {}^{T \to I_T \to G \text{~is the unit}}_{I_T |_{C_i} \to G|_{C_i} \text{ factors through }H_i \text{ for all } i }\big\} \\
 &= \big\{ T \to \bbL ie (G) \in \Hom _{S} ( T, \bbL ie (G) ) | ~~{}_{ T |_{C_i} \to  \bbL ie (G) |_{C_i} \text{ factors through } \bbL ie ( H_i ) \text{ for all } i }\big\} \\
 &= \big( \Bl \big\{{}^{\bbL ie (G) |_{ C_i}}_{\bbL ie ( H_i )} \big\}_{i \in I} \bbL ie (G) \big)(T).
\end{align*}
This finishes the proof.
\xpf

The following result generalizes the fact that congruence groups are normal subgroups, it is related to Example \ref{yuex} (note that the proof of \cite[Lemma 1.4]{Yu01} is not correct, cf. Remark \ref{yuproof}).

\prop \label{normalizes}
Assume that $C_i$ is a Cartier divisor in $S$ for all $i$. Assume that $G \to S$ is a flat group scheme. Let $\eta : K\to G$ be a morphism of group schemes over $S$ such that $K \to S$ is flat. Assume that $H_i  \subset G $ is a closed subgroup scheme over $S$ such that $H_i \to S$ is flat for all $i$. Assume that $\Bl_{H}^{C} G \to S$ is flat (and in particular a group scheme). Assume that, for all $i$,
 $K_{C_i} $ commutes with ${H_i}_{C_i}$ in the sense that the morphism $K_{C_i} \times _{C_i} {H_i}_{C_i} \to G_{C_i} $, $(k,h) \mapsto \eta(k)h\eta(k)^{-1} $ equals the composition morphism $K_{C_i} \times _{C_i} {H_i}_{C_i} \to {H_i}_{C_i} \subset G_{C_i}$, $(k,h) \mapsto h$. Then $K$ normalizes $\Bl _H^C G$, more precisely the solid composition map \[
 \begin{tikzcd} K \times _S \Bl _{H}^C G \ar[rr, "Id \times \Bl "] \ar[rrd, dashrightarrow] & & K \times _S G \ar[rrr, "k {,} g \mapsto \eta(k)g\eta(k)^{-1} "] & & & G\\ &  &\Bl _H ^C G  \ar[rrru, "\Bl"]&  & \end{tikzcd} \]
 factors uniquely through $\Bl _H^C G $.
 \xprop
 \pf
 Let $\phi$ be the solid composition map, we claim that it belongs to $Sch_G^{G_C \text{-reg}}$. Let us prove this claim. The map $\theta: K \times _S G \xrightarrow{k {,} g \mapsto \eta(k)g\eta(k)^{-1}} G$ is flat. Indeed it is the composition of an isomorphism, namely $K \times _S G \xrightarrow{(k,g) \mapsto (k,\eta (k)g\eta(k)^{-1})} K \times _S G$, with a flat morphism, namely the projection on the second factor $K \times _S G \to G$. So $\theta ^{-1} ( G_{C_i})$ is a Cartier divisor in $K \times _S G$ (note that $G_{C_i}$ is a Cartier divisor in $G$ because $G \to S$ is flat) for all $i$. Lemma \ref{blow.up.base.change.lemm} and \ref{defined} show that $(Id \times \Bl)^{-1} (\theta^{-1} (G_{C_i}))$ is a Cartier divisor. This finishes to prove the claim. Now  by Proposition \ref{blow.up.rep.prop}, $\phi $ factors uniquely through $\Bl _H^C G$ if and only if $\phi|_{G_{C_i}}$ factors through ${H_i}_{C_i}$ for all $i$. The following diagram, obtained using Proposition \ref{blow.up.Cartier.lemm}, finishes the proof \[
 \begin{tikzcd}[row sep=3ex, column sep=-1ex]
(K \times _S \Bl _H ^{C} G ) \times _G G_{C_i} \ar[rr] \ar[d, equal] &  &(K \times _S G) \times _G G_{C_i} \ar[d, equal] \ar[rr]  & & G_{C_i} \ar[d, equal] \\
 K_{C_i} \times _{C_i} ( \Bl _H ^{C} G \times _G G_{C_i} ) \ar[rd, dashrightarrow] \ar[rr] & & K_{C_i} \times _{C_i} G_{C_i}  \ar[rr]  & & G_{C_i} \\
 & K_{C_i} \times _{C_i} {H_i}_{C_i} \ar[ru]\ar[rr, dashrightarrow] & & {H_i}_{C_i}. \ar[ru] &
 \end{tikzcd}\] 
 \xpf 

\rema \label{yuproof} 
The proof of \cite[Lemma 1.4]{Yu01} is not correct. Indeed, the fifth assertion, starting by "Since $G^i(E)_{y,0}$ normalizes [...]" is wrong. Indeed it is wrong in the $SL_2$ case. Take $G^0 \subset G^1$ such that $G^1_E = SL_2/E$ and  $G^0_E = T/E$ where $T/ \bbZ$ is the diagonal split torus. The group $SL_2 (\calO _E) $ does not normalize $T (\calO _E)$.
To obtain a correct proof in the setting of Example \ref{yuex}, one can proceed as follows.
\begin{enumerate}
\item Preserve the first and second sentences of \cite[Proof of Lemma 1.4]{Yu01}.
 
\item Replace the third sentence of \cite[Proof of Lemma 1.4]{Yu01}  by "For $i >0$, $G^i (E)_{y,r_i} \cdots G^d (E) _{y , r_d}$ is a group by induction hypothesis".
 
 \item Replace the fourth and fifth sentences of \cite[Proof of Lemma 1.4]{Yu01} by "Using Proposition \ref{normalizes}, one proves that $G^{i-1} (E) _{y,r_{i-1}}$ normalizes  $G^i (E)_{y,r_i} \cdots G^d (E) _{y , r_d}$, so we see that $G^{i-1} (E)_{y,r_{i-1}} \cdots G^d (E) _{y , r_d}$ is a group".
  
\item  Preserve the sixth sentence of \cite[Proof of Lemma 1.4]{Yu01}.
  \end{enumerate}
\xrema 

\rema \label{sketchproofexyu}
In this remark, we sketch the proof of the identity (\ref{yuid}) appearing in Example \ref{yuex} using \cite{Yu15}. The inclusion $\overset{\to}{\mathrm{G}}(E) _{x , \overset{\to}{r}} \subset  \Bl _{e_G, G^0 , G^1 , \ldots ,G^i , \ldots ,G^{d-1} }^{r_0 , r_1 , r_2, \ldots , r_{i+1} , \ldots , r_d} G (\calO)$ holds since $\overset{\to}{\mathrm{G}}(E) _{x , \overset{\to}{r}}$ is generated by some subgroups belonging to  $\Bl _{e_G, G^0 , G^1 , \ldots ,G^i , \ldots ,G^{d-1} }^{r_0 , r_1 , r_2, \ldots , r_{i+1} , \ldots , r_d} G (\calO)$. To prove the reverse inclusion, we reduce to the case where all $r_i$ are strictly positive, then we use \cite[Lemma §7.4]{Yu15} to write an element in  $\Bl _{e_G, G^0 , G^1 , \ldots ,G^i , \ldots ,G^{d-1} }^{r_0 , r_1 , r_2, \ldots , r_{i+1} , \ldots , r_d} G (\calO)$ (which is a subgroup of $\Bl _{e_G }^{r_0} G (\calO)$) uniquely as product of elements in roots groups. Then we prove that each factor belongs to roots groups considered by Yu, reasoning by contradiction.
\xrema
\rema   
Note that still in relation with Moy-Prasad filtrations, dilatations were also used to study Berkovich's point of view on Bruhat-Tits buildings in positive depth \cite{Ma22}. 
\xrema 
 
\section{Congruent isomorphisms} \label{sectioniso}
Let $(\calO ,\pi)$ be an henselian pair where $\pi \subset \calO $ is an invertible ideal.  Let $S= \Spec (\calO)$ and let $C= \Spec (\calO /\pi)$. If $G\to S$ is a group scheme $\Lie (G) \to S$ denote the underlying group scheme of $\bbL ie (G) \to S$ (cf. §\ref{néron}).

\theo \label{isocongruent}
 Let $G$ be a separated and smooth goup scheme over $S$. Let $H_0 , H_1 , \ldots , H_k$ be closed subgroup schemes of $G$ such that $H_0=e_G$ is the trivial subgroup. Let $s_0 , s_1 , \ldots , s_k$ and $r_0 , r_1 , \ldots , r_k$ be in $\bbN$ such that  
 \begin{enumerate}
 \item $s_i \geq s_0 $ and $r_i \geq r_0 $ for all $i \in \{0, \ldots , k\}$
 \item $ r_i \geq s_i $ and $r_i-s_i \leq s_0$ for all $i \in \{0, \ldots , k \}$.
  \end{enumerate} Assume that $G$ is affine or $\calO $ is local.  Assume that the regularity condition (RC) introduced below is satisfied (cf. Definition \ref{rc}). Then we have a canonical isomorphism of groups
 \[ \Bl_{H_0, H_1 , \ldots , H_k}^{s_0 ,~ s_1 ,~ \ldots , s_k} G (\calO) /  \Bl_{H_0, H_1 , \ldots , H_k}^{r_0 ,~ r_1 ,~ \ldots , r_k} G   (\calO) \cong \Lie (\Bl_{H_0, H_1 , \ldots , H_k}^{s_0 ,~ s_1 ,~ \ldots , s_k} G ) (\calO) / \Lie ( \Bl_{H_0, H_1 , \ldots , H_k}^{r_0 ,~ r_1 ,~ \ldots , r_k} G   )(\calO).  \]
\xtheo 

\pf  For $i \in \{0 , \ldots , k \}$, put $t_i= r_i-s_i$. As a first step in our proof, we assume that $t_i=t_j=:t$ for all $i,j\in \{0,\ldots , k\}.$
Proposition \ref{clef} shows that \[\Bl_{H_0, H_1 , \ldots , H_k}^{r_0 ,~ r_1 ,~ \ldots , r_k} G = \Bl_{H_0}^{2t}   \Bl_{H_0, H_1 , \ldots , H_k}^{s_0 - t ,  \ldots , s_k-t} G \text{ and } \Bl_{H_0, H_1 , \ldots , H_k}^{s_0 ,~ s_1 ,~ \ldots , s_k} G = \Bl_{H_0}^{t}   \Bl_{H_0, H_1 , \ldots , H_k}^{s_0 - t ,  \ldots , s_k-t} G.\] Put $G' =  \Bl_{H_0, H_1 , \ldots , H_k}^{s_0 - t ,  \ldots , s_k-t} G$. The scheme $G'$ is smooth over $S$ by (RC). By \cite[Theorem 4.3]{MRR20}, we have a canonical isomorphism
\[ \Bl_{H_0}^{t} G' (\calO) / \Bl_{H_0}^{2t}  G' (\calO) \cong \Lie ( \Bl_{H_0}^{t}   G' ) (\calO)/ \Lie (\Bl_{H_0}^{2t} G') (\calO) .\] This finishes the proof of the case where $t_i=t_j $ for all $i,j \in \{1 , \ldots , k \}$. Now we prove the general case. Put $t_m= \max _{i \in \{0 , \ldots ,k\}} t_i$. The isomorphism \cite[Theorem 4.3]{MRR20} is functorial in $G$. Applying this functoriality to the morphism $ \Bl_{\{H_i\}_{0 \leq i \leq k}}^{\{s_i-t_m\}_{0 \leq i \leq k}} G \to \Bl_{\{H_i\}_{0 \leq i \leq k}}^{\{r_i-2t_{m}\}_{0 \leq i \leq k}} G$ and with the integers given by the inequality $0 \leq \frac{2t_m}{2} \leq t_m \leq 2t_m $, we get a canonical commutative diagram
\begin{scriptsize}\[
\begin{tikzcd}[row sep=6ex, column sep=-24.5ex]
 & \Bl_{\{H_i\}_{0 \leq i \leq k}}^{\{s_i\}_{0 \leq i \leq k}} G(\calO)/\Bl_{\{H_i\}_{0 \leq i \leq k}}^{\{r_i\}_{0 \leq i \leq k}} G (\calO)
\ar[rd,hookrightarrow] & \\
  \Bl_{\{H_i\}_{0 \leq i \leq k}}^{\{s_i\}_{0 \leq i \leq k}} G(\calO)/\Bl_{\{H_i\}_{0 \leq i \leq k}}^{\{s_i +t_{m}\}_{0 \leq i \leq k}} G(\calO) \ar[ru, two heads] \ar[rr]\ar[d, equal]& &
\Bl_{\{H_i\}_{0 \leq i \leq k}}^{\{r_i-t_{m}\}_{0 \leq i \leq k}} G(\calO)/\Bl_{\{H_i\}_{0 \leq i \leq k}}^{\{r_i\}_{0 \leq i \leq k}} G (\calO)\ar[d,equal]\\ 
   \Lie (\Bl_{\{H_i\}_{0 \leq i \leq k}}^{\{s_i\}_{0 \leq i \leq k}} G)(\calO)/ \Lie (\Bl_{\{H_i\}_{0 \leq i \leq k}}^{\{s_i +t_{m}\}_{0 \leq i \leq k}} G)(\calO) \ar[rr] \ar[rd,two heads] & &
 \Lie (\Bl_{\{H_i\}_{0 \leq i \leq k}}^{\{r_i-t_{m}\}_{0 \leq i \leq k}} G)(\calO)/ \Lie (\Bl_{\{H_i\}_{0 \leq i \leq k}}^{\{r_i\}_{0 \leq i \leq k}} G)(\calO) \\
  &\Lie ( \Bl_{\{H_i\}_{0 \leq i \leq k}}^{\{s_i\}_{0 \leq i \leq k}} G) (\calO)/ \Lie (\Bl_{\{H_i\}_{0 \leq i \leq k}}^{\{r_i\}_{0 \leq i \leq k}} G )(\calO) .\ar[ru,hookrightarrow]
 & 
\end{tikzcd}
\]\end{scriptsize}  The injectivity of the two hookarrows follows from Remark \ref{remainj}. This identifies \[ \Bl_{\{H_i\}_{0 \leq i \leq k}}^{\{s_i\}_{0 \leq i \leq k}} G(\calO)/\Bl_{\{H_i\}_{0 \leq i \leq k}}^{\{r_i\}_{0 \leq i \leq k}} G (\calO)\] and \[\Lie ( \Bl_{\{H_i\}_{0 \leq i \leq k}}^{\{s_i\}_{0 \leq i \leq k}} G) (\calO)/ \Lie (\Bl_{\{H_i\}_{0 \leq i \leq k}}^{\{r_i\}_{0 \leq i \leq k}} G )(\calO)\] inside the right part of the diagram.
\xpf 

\defi \label{rc} Let $G , \{H_i,s_i,r_i\}_{0 \leq i \leq k }$ be as in Theorem \ref{isocongruent}. Put $t_m = \max _{i \in \{0 , \ldots , k\}} (r_i-s_i)$. We introduce the following regularity condition
\begin{center}
(RC) $  \Bl_{\{H_i\}_{0 \leq i \leq k}}^{\{s_i-t_m\}_{0 \leq i \leq k}} G  $ and $ \Bl_{\{H_i\}_{0 \leq i \leq k}}^{\{r_i-2t_{m}\}_{0 \leq i \leq k}} G$ are smooth over $S$.
\end{center}
\xdefi
We recall that Proposition \ref{flatandsmooth} offers a way to check (RC) in many cases. We finish with the following result.
\coro \label{isocorocongru}
Let $G$ be a separated and smooth goup scheme over $S$. Let $H_0 \subset H_1 \subset \ldots \subset  H_k$ be closed subgroup schemes of $G$ such that $H_i$ is smooth over $S$ for $0 \leq i \leq d$ and $H_0=e_G$. Let $s_0 , s_1 , \ldots , s_k$ and $r_0 , r_1 , \ldots , r_k$ be in $\bbN$ such that  
 \begin{enumerate}
 \item $s_i \geq s_0 $ and $r_i \geq r_0 $ for all $i \in \{0, \ldots , k\}$
 \item $ r_i \geq s_i $ and $r_i-s_i \leq s_0$ for all $i \in \{0, \ldots , k \}$.
  \end{enumerate}  Assume that $G$ is affine or $\calO $ is local.  Then we have a canonical isomorphism of groups
 \begin{small}\[ \Bl_{H_0, H_1 , \ldots , H_k}^{s_0 ,~ s_1 ,~ \ldots , s_k} G (\calO)/  \Bl_{H_0, H_1 , \ldots , H_k}^{r_0 ,~ r_1 ,~ \ldots , r_k} G  (\calO)  \cong \Lie (\Bl_{H_0, H_1 , \ldots , H_k}^{s_0 ,~ s_1 ,~ \ldots , s_k} G ) (\calO)/ \Lie ( \Bl_{H_0, H_1 , \ldots , H_k}^{r_0 ,~ r_1 ,~ \ldots , r_k} G   )(\calO).  \]\end{small}
\xcoro 
\pf By Theorem \ref{isocongruent}, it is enough to check the condition (RC), it follows from \cite[Exp. III, Proposition 4.15]{SGA3} and Proposition \ref{flatandsmooth}.\xpf
\rema Note that the result of the proof of \cite[Exp. III, Proposition 4.15]{SGA3} is stronger than its statement. Indeed the statement uses Koszul-regularity and the proof shows regularity (in the terminology of \cite[\href{https://stacks.math.columbia.edu/tag/063J}{Tag 063J}]{stacks-project}).
\xrema 

\section{Interpretation of Rost double deformation space as dilatation}
\label{sectionrost}
We interpret Rost double deformation space \cite{Ro96} in the language of dilatations. This section emerged after a question of A. Dubouloz.
Let $Z \to Y \to X$ be closed immersions (in \cite{Ro96}, all schemes are assumed to be defined over fields but we work with arbitrary schemes here). Let $\overline{D} ( X,Y,Z)$ be the double deformation space as defined in \cite[§10]{Ro96}.
Let $\bbA^2$ be $\Spec (\bbZ [s,t])$. Let $D_s$, $ D_{st}$ and $D_{s^2t}$ be the locally principal closed subschemes of $ \bbA ^2$ defined by the ideals $(s),(st)$ and $(s^2t)$. We now omit the subscript ${}_{\Spec (\bbZ )}$ in fiber products.

\prop \label{rost2} We have a canonical identification 
\[\overline{D} ( X,Y,Z) \cong \Bl _{ (Y \times \bbA ^2), ~ ~(Z \times \bbA ^2)}^{ (X \times D_{st}), ~(X \times D_s)}  (X \times  \bbA ^2).\]
In other words, Rost double deformation space is canonically interpreted as a double-centered dilatation.
\xprop 
\pf
The definition of $\overline{D}(X,Y,Z)$ is given in \cite[10.5]{Ro96} locally for affine schemes. So we reduce to the case where $X= \Spec (A)$, $Y = \Spec (A/I)$, $Z = \Spec (A /J)$ are affine. Then $\overline{D}(X,Y,Z)$ is defined as the spectrum of the ring $O_{\overline{D}}=\sum _{n, m} I^n J^{m-n}t^{-n} s^{-m} \subset A [t,s,t^{-1}, s^{-1}] \cong A \otimes _{\bbZ} \bbZ[t,s,t^{-1},s^{-1}]$ where $ I^k = J^k = A $ for $k <0$ as in \cite[§10.2]{Ro96}.

We claim that $O _{\overline{D}}$ is equal to the sub-$A[t,s]$-algebra of $A[t,s,t^{-1},s^{-1}]$ generated by $I(ts)^{-1}$ and $Js^{-1}$. Indeed, let $a,b \in \bbN$ and put $m= a+b $ and $n =a$. Then $(I(ts)^{-1})^a( Js^{-1})^b = I^n J ^{m-n} t^{-n} s^{-m}$. So $(A[t,s] ) [ I(ts)^{-1}, Js^{-1} ] $ is included in $O _{\overline{D}}$. Reciprocally, let $n, m \in \bbZ$. Assume firstly that $n <0 $ and put $l=-n >0$, then \[I^n J^{m-n} t^{-n} s^{-m} = J^{m+l} t^l s^{-m} \subset J^m s^{-m} t^l \subset  A[t,s][Js^{-1}] \subset A[t,s][I(ts)^{-1}, Js^{-1}].\]
Assume secondly that $n \geq 0  $ and $m \geq n$ and put $a=n $ and $b =m-n \geq 0$, then 
\[I^n J^{m-n} t^{-n} s^{-m} = (I(ts)^{-1} )^a (Js^{-1})^b \subset A[t,s][I(ts)^{-1}, Js^{-1}]  .\]
Assume thirdly that $n \geq 0$ and $m-n  <0$ and put $c = n-m >0$, then 
\[ I^n J^{m-n} t^{-n} s^{-m} =(I(ts)^{-1})^n s^{c}  \subset  A[t,s][I(ts)^{-1}] \subset A[t,s][I(ts)^{-1}, Js^{-1}] . \]
So in all cases, $ I^n J^{m-n} t^{-n} s^{-m}\subset A[t,s][I(ts)^{-1}, Js^{-1}]$. This finishes to prove our claim.
Now Fact \ref{rfvb} finishes the proof of Proposition \ref{rost2}.
\xpf

\section*{Declarations}

\subsection*{Conflict of interest}
 The author has no conflict of interest to declare that are relevant to this article.

\subsection*{Funding} 
This project has received funding from the European Research Council (ERC) under the European Union’s Horizon 2020 research and innovation programme (grant agreement No 101002592) and the ISF grant 1577/23.

\subsection*{Data Availability Statements}

Data sharing not applicable to this article as no datasets were generated or analysed during the current study.


\begin{thebibliography}{999999}



\bibitem[AZ55]{abhyankar-zariski55} S. S. Abhyankar and O. Zariski, {\it Splitting of valuations in extensions of local domains}. Proc. Natl. Acad. Sci. USA,  41 (1955), 84--90.


\bibitem[Ana73]{Ana73} S.\,Anantharaman: {\it Sch\'emas en groupes: espaces homog\`enes et espaces alg\'ebriques sur une base de dimension $1$}, Sur les groupes alg\'ebriques, Soc. Math. France, Paris (1973), pp. 5-79. Bull. Soc. Math. France, M\'em. \textbf{33}.



\bibitem[Ar69]{artin69}M. Artin,
{\it Algebraic approximation of structures over complete local rings}. 
Inst. Hautes \'Etudes Sci. Publ. Math.   \textbf{36} (1969) 23 -- 58. 




\bibitem[BLR90]{BLR90} S.\,Bosch, W.\,L\"utkebohmert and M.\,Raynaud: {\it N\'eron models},  Ergebnisse der Mathematik und ihrer Grenzgebiete (3) \textbf{21}, Berlin, New York: Springer-Verlag.



\bibitem[Ca58]{Ca58} P.\,Cartier: {\it Questions de rationalit\'e des diviseurs en g\'eom\'etrie alg\'ebrique.}
Bull. Soc. Math. France 86 (1958), 177–251.








\bibitem[Da67]{Da67} E.\, Davis: {\it Ideals of the principal class, R-sequences and a certain monoidal transformation.} Pacific J. Math. 20 (1967), 197–205.




\bibitem[DG70]{SGA3} M.\,Demazure, A.\,Grothendieck: {\it S\'eminaire de G\'eom\'etrie Alg\'ebrique du Bois Marie - 1962-64 - Sch\'emas en groupes - (SGA 3)}, Lecture notes in mathematics \textbf{151} (1970), Berlin; New York: Springer-Verlag. pp. xv+564.


\bibitem[DHdS18]{DHdS18} N. D.\,Duong, P. H.\,Hai, J. P.\, dos Santos:
{\it On the structure of affine flat group schemes over discrete valuation
rings, I}. Ann. Sc. Norm. Super. Pisa Cl. Sci. (5) 18 (2018), no. 3,
977--1032.

\bibitem[DMdS23]{DMdS23} A.\,Dubouloz, A.\,Mayeux, J.\,P.\,dos\,Santos: {A survey on algebraic dilatations}, 2023

\bibitem[Du05]{Du05} A.\,Dubouloz: {\it Quelques remarques sur la notion de modification affine.} 
https://arxiv.org/abs/math/0503142, 2005.






\bibitem[Ka94]{Ka94} S.\, Kaliman: {\it Exotic analytic structures and Eisenman intrinsic measures.} Israel J. Math. 88 (1994), no. 1-3, 411–423.

\bibitem[KP22]{KP22} T. Kaletha, G. Prasad: {\it Bruhat-Tits Theory: A New Approach}, New Mathematical Monographs, Series Number 44, 2022


\bibitem[KZ99]{KZ99}S.\,Kaliman, M.\,Zaidenberg: {\it
Affine modifications and affine hypersurfaces with a very transitive automorphism group.}
Transform. Groups 4 (1999), no. 1, 53–95. 


\bibitem[Ma22]{Ma22} A.\,Mayeux: {\it Bruhat-Tits theory from Berkovich’s point of view. Analytic filtrations}, Ann. H. Lebesgue 5 (2022) 813-839.


\bibitem[Ma23c]{Ma23c} A.\,Mayeux: {\it Dilatations of categories}, https://arxiv.org/abs/2305.11303, 2023



\bibitem[MRR20]{MRR20}  A.\,Mayeux,  T.\,Richarz,  M.\,Romagny:  {\it N\'eron blowups and low-degree cohomological applications}, Algebr. Geom. 10 (6) (2023) 729–753.



\bibitem[MR24]{MR24}A.\,Mayeux, S.\,Riche: {\it On multi-graded proj schemes}, https://arxiv.org/html/2310.13502v3, 2024

\bibitem[MY24]{MY24} A.\,Mayeux, Y.\,Yamamoto: {\it Comparing Bushnell-Kutzko and Sécherre’s constructions of types for $GL_N$ and its inner forms with Yu’s construction},  Bull. Soc. Math. France. 152 (2024).

\bibitem[Ner64]{Ner64} A.\,N\'eron: {\it Mod\`eles minimaux des
vari\'et\'es ab\'eliennes sur les corps locaux et globaux},
Inst. Hautes \'Etudes Sci. Publ. Math. No. 21 (1964), 128 pp.

\bibitem[No1884]{No1884} Max Noether, {\it Rationale Ausf\"uhrung der Operationen in der Theorie der
algebraischen Functionen}. Math. Ann.  23, 1884,  311--358.



\bibitem[PY06]{PY06} G.\,Prasad, J.\,-K.\,Yu; {\it On quasi-reductive group schemes},
J. Algebraic Geom. \textbf{15} (2006), 507--549. With an appendix by Brian Conrad.


\bibitem[PZ13]{PZ13} G.\,Pappas, X.\,Zhu: {\it Local models of Shimura varieties and a conjecture of Kottwitz}, Invent.\,Math.\,\textbf{194} (2013), 147--254.



\bibitem[Ro96]{Ro96} M. Rost, \emph{Chow groups with coefficients}. Doc. Math. 1 (1996), No. 16, 319-393.





\bibitem[StP]{stacks-project} The {Stacks project authors}: {\it The Stacks Project}.



\bibitem[WW80]{WW80} W.\,C.\,Waterhouse, B.~Weisfeiler: {\it One-dimensional affine group schemes}, J.~Algebra \textbf{66} (1980), 550--568.

\bibitem[Yu01]{Yu01} J.-K.\,Yu: {\it Construction of tame supercuspidal representations},
J. Amer. Math. Soc.~\textbf{14} (2001), no. 3, 579--622.

\bibitem[Yu15]{Yu15} J.-K.\,Yu: {\it Smooth models associated to concave functions in Bruhat-
Tits theory}, Autour des sch\'emas en groupes. Vol. III, 227--258, Panor. Synth\`eses \textbf{47}, Soc. Math. France, Paris, 2015. 


\bibitem[Za43]{Za43} O.\, Zariski: {\it  Foundations of a general theory of birational correspondences,} Trans. Amer. Math. Soc. 53 (1943).



\end{thebibliography}
\end{document}